\pdfoutput=1

\documentclass[a4paper, 11pt, titlepage, oneside]{scrartcl}
\makeatletter

\usepackage[T1]{fontenc}
\usepackage[utf8]{inputenc}

\usepackage[onehalfspacing]{setspace}
\usepackage[headsepline]{scrlayer-scrpage}
\clearpairofpagestyles
\usepackage{layout}
\usepackage{geometry}
\usepackage{adjustbox}
\usepackage{tcolorbox}
\usepackage{authblk}
\usepackage[titletoc,title]{appendix}
\usepackage{comment}

\usepackage[hyphens]{url}

\usepackage{color}

\usepackage{amsmath,amssymb,amsfonts,amsthm, mathtools, bm}
\usepackage{bbold}
\usepackage{thmtools, thm-restate}
\usepackage{nicefrac}
\usepackage{array}

\usepackage[normal]{threeparttable}
\usepackage{threeparttablex}
\usepackage{tabularx}
\usepackage{longtable}
\usepackage{booktabs}
\usepackage{colortbl}
\usepackage{multirow}
\usepackage{makecell} 

\usepackage[acronym]{glossaries}

\usepackage{svg} 
\usepackage[font=footnotesize, labelfont=bf]{caption}
\usepackage{subcaption}
\usepackage{floatrow}
\usepackage{graphicx}
\usepackage{placeins}
\graphicspath{{Figures/}}
\usepackage{wrapfig}

\usepackage{algorithm}
\usepackage{algpseudocode}
\usepackage[noEnd=false]{algpseudocodex}

\usepackage{tikz}
\usetikzlibrary{math}
\usetikzlibrary{shapes}
\usepackage{pgfplots,color}
\pgfplotsset{compat=1.8}
\usepgfplotslibrary{groupplots}
\usepgfplotslibrary{dateplot}
\usepackage{tikzscale}
\usepackage{pbox}

\usepackage{enumerate}
\usepackage{multicol} 
\usepackage{hhline}
\usepackage{etoolbox}
\usepackage{alphalph}

\usepackage{eurosym}
\usepackage{ellipsis}
\usepackage{bm}
\usepackage{bbm}

\usepackage{natbib}

\usepackage{optidef}

\usepackage{import}
\usepackage{scrlayer-scrpage}
\AtBeginDocument{
	\newgeometry{
		textheight=9in,
		textwidth=6.5in,
		top=1in,
		headheight=16.5pt,
		headsep=25pt,
		footskip=30pt
	}
}
\newsavebox{\abstractbox}
\renewenvironment{abstract}
{\begin{lrbox}{0}\begin{minipage}{\textwidth}
			\begin{center}\normalfont\sectfont\abstractname\end{center}\quotation}
		{\endquotation\end{minipage}\end{lrbox}%
	\global\setbox\abstractbox=\box0 }

\makeatletter
\expandafter\patchcmd\csname\string\maketitle\endcsname
{\vskip\z@\@plus3fill}
{\vskip\z@\@plus2fill\box\abstractbox\vskip\z@\@plus1fill}
{}{}
\makeatother

\addtokomafont{captionlabel}{\footnotesize\bfseries} 
\addtokomafont{caption}{\footnotesize\bfseries} 

\bibpunct[, ]{(}{)}{,}{a}{}{,}%
\newtheorem{definition}{Definition}[section]

\newtheorem{result}{Result}

\counterwithin*{equation}{section}

\newenvironment{myfont}{\fontfamily{ccr}\selectfont}{\par}
\DeclareTextFontCommand{\textmyfont}{\myfont}

\AtBeginDocument{%
	\abovedisplayskip=7.5pt plus 0pt minus 0pt
	\abovedisplayshortskip=7.5pt plus 0pt minus 0pt
	\belowdisplayskip=7.5pt plus 0pt minus 0pt
	\belowdisplayshortskip=7.5pt plus 0pt minus 0pt
}

\newcolumntype{L}[1]{>{\raggedright\let\newline\\\arraybackslash\hspace{0pt}}p{#1}}
\newcolumntype{C}[1]{>{\centering\let\newline\\\arraybackslash\hspace{0pt}}p{#1}}
\newcolumntype{R}[1]{>{\raggedleft\let\newline\\\arraybackslash\hspace{0pt}}p{#1}}
\def\checkmark{\tikz\fill[scale=0.4](0,.35) -- (.25,0) -- (1,.7) -- (.25,.15) -- cycle;}

\renewcommand{\emph}[1]{\textit{#1}}

\makeatletter
\newalphalph{\aalphalph}[mult]{\alphalph@alph}{26}
\newcommand{\alphalphval}[1]{%
  \@ifundefined{c@#1}{
    \aalphalph{#1}
  }{%
    \aalphalph{\value{#1}}
  }
}
\makeatother

\AtBeginEnvironment{subequations}{%
}

\clearpairofpagestyles   
\ohead{\pagemark}

\newacronym{acr:ev}{EV}{electric vehicle}
\newacronym{acr:spp}{SPP}{shortest path problem}
\newacronym{acr:rcspp}{RCSPP}{resource constrained shortest path problem}
\newacronym{acr:osm}{OSM}{{O}pen {S}treet {M}ap}
\newacronym{acr:lp}{LP}{linear programming}
\newacronym{acr:vrp}{VRP}{Vehicle Routing Problem}
\newacronym{acr:frlp}{FRLP}{Flow Refueling Location Problem}
\newacronym{acr:mip}{MIP}{Mixed Integer Programming}
\newacronym{acr:bab}{B\&B}{Branch-and-bound}
\newacronym{acr:ils}{ILS}{Iterated Local Search}
\newacronym{acr:soc}{SoC}{state of charge}
\newacronym{acr:evrptw-pr}{EVRPTW-PR}{electric vehicle routing problem with time windows and partial recharging}
\newacronym{acr:olev}{OLEV}{On-Line Electric Vehicle}
\newacronym{acr:rad}{R\&D}{Research \& Development}
\newacronym{acr:milas}{MILAS}{\textit{Modulare intelligente induktive Ladesysteme für autonome Shuttles}}

\newacronym{acr:strip-dyn}{\textsc{SDS}}{\textsc{StripDynamicStation}}
\newacronym{acr:remove-stat}{\textsc{RSS}}{\textsc{RemoveStaticStation}}
\newacronym{acr:swap-to-dyn}{\textsc{S2D}}{\textsc{Swap2DynamicStation}}
\newacronym{acr:swap-to-stat}{\textsc{S2S}}{\textsc{Swap2StaticStation}}
\newacronym{acr:add-dyn}{\textsc{ADS}}{\textsc{AddDynamicStation}}
\newacronym{acr:add-stat}{\textsc{ASS}}{\textsc{AddStaticStation}}

\graphicspath{{./figures/}}

\begin{document}
\emergencystretch 3em

\newcommand{\ConsumptionFunction}{C}
\newcommand{\Graph}{G}
\newcommand{\VehicleBound}{L^{\Vehicle}}
\newcommand{\Path}{P}
\newcommand{\RechargeFunction}{R}
\newcommand{\Solution}{S}
\newcommand{\TravelTimeFunction}{T}
\newcommand{\Cost}{c}
\newcommand{\Depot}{d}
\newcommand{\TwBegin}{e}
\newcommand{\Station}{f}
\newcommand{\HeuristicComponent}{g}
\newcommand{\StationChargingRate}{h}
\newcommand{\Line}{k}
\newcommand{\TwEnd}{l}
\newcommand{\EnergyPrice}{p}
\newcommand{\Consumption}{q}
\newcommand{\Recharge}{r}
\newcommand{\Stop}{s}
\newcommand{\TravelTime}{t}
\newcommand{\Vertex}{v}
\newcommand{\TransformerVar}{w}
\newcommand{\ArcVar}{x}
\newcommand{\StaticChargerVar}{y}
\newcommand{\DynamicChargerVar}{z}
\newcommand{\SetOfArcs}{\mathcal{A}}
\newcommand{\StopsOfLine}{\mathcal{C}^{\Vehicle}}
\newcommand{\SetOfStations}{\mathcal{F}}
\newcommand{\SetOfSegmentStarts}{\mathcal{I}}
\newcommand{\SetOfSegmentEnds}{\mathcal{J}}
\newcommand{\SetOfLines}{\mathcal{K}}
\newcommand{\SetOfLabels}{\mathcal{L}}
\newcommand{\SetOfSegments}{\mathcal{O}}
\newcommand{\Queue}{\mathcal{Q}}
\newcommand{\SetOfStops}{\mathcal{S}}
\newcommand{\SetOfTransformers}{\mathcal{T}}
\newcommand{\SetOfVertices}{\mathcal{V}}
\newcommand{\ExploredLabels}{\mathcal{X}}
\newcommand{\SetOfStationaryChargerReps}{\mathcal{Y}^{\Vehicle}}
\newcommand{\SegmentStart}{\alpha}
\newcommand{\LabellingIterations}{\beta}
\newcommand{\LabelCost}{\gamma}
\newcommand{\NeighborhoodSize}{\zeta}
\newcommand{\HeuristicSortingKey}{\eta}
\newcommand{\SegmentIndex}{\theta}
\newcommand{\ILSInitStepSize}{\kappa}
\newcommand{\Label}{\lambda}
\newcommand{\Requirement}{\mu}
\newcommand{\ILSIterations}{\xi}
\newcommand{\MaxChargingTimeSteps}{\pi}
\newcommand{\StateOfChargeVar}{\rho}
\newcommand{\Headway}{\sigma}
\newcommand{\ArrivalTimeVar}{\tau}
\newcommand{\PerturbationStrength}{\phi}
\newcommand{\TimeToCharge}{\chi}
\newcommand{\StopOfLine}{\psi}
\newcommand{\SegmentEnd}{\omega}

\newcommand{\ChargingTimeVar}{\Delta\tau}
\newcommand{\RechargeVar}{\Delta\rho}

\newcommand{\StartDepot}{\Depot^{+}}
\newcommand{\EndDepot}{\Depot^{-}}
\newcommand{\SetOfVehicles}{\SetOfLines}
\newcommand{\Vehicle}{\Line}
\newcommand{\SetOfStaticStations}{\SetOfStations^s}
\newcommand{\SetOfDynamicStations}{\SetOfStations^d}
\newcommand{\MinSoC}{\Consumption_{\min}}
\newcommand{\MaxSoC}{\Consumption_{\max}}
\newcommand{\UnbracedArc}{i, j}
\newcommand{\Arc}{(\UnbracedArc)}
\newcommand{\ArcTravelTime}{\TravelTime}
\newcommand{\ArcConsumption}{\Consumption}
\newcommand{\UnbracedVehicleArc}{i, j, \TravelTime, \Consumption, \Recharge}
\newcommand{\VehicleArc}{(\UnbracedVehicleArc)}

\newcommand{\Segment}{(\SegmentStart, \SegmentEnd)}
\newcommand{\UnbracedSegment}{\SegmentStart, \SegmentEnd}

\newcommand{\InvestmentCost}{\Cost}
\newcommand{\TransformerCost}{\Cost}
\newcommand{\InfrastructureCost}{\Cost_{\texttt{I}}}
\newcommand{\RoutingCost}{\Cost_{\texttt{R}}}
\newcommand{\InitialSoC}{\Consumption_{\text{init}}}
\newcommand{\TimeStep}{\delta} 
\newcommand{\ForwardLabel}{\overrightarrow{\Label}}
\newcommand{\BackwardLabel}{\overleftarrow{\Label}}

\newcommand{\VertexEnergyBound}{\underline{\Vertex}}
\newcommand{\SetOfForwardLabels}{\SetOfLabels^{\texttt{F}}}
\newcommand{\SetOfBackwardLabels}{\SetOfLabels^{\texttt{B}}}

\newcommand{\ExploredFwLabels}{\overrightarrow{\ExploredLabels}}
\newcommand{\ExploredBwLabels}{\overleftarrow{\ExploredLabels}}
\newcommand{\FwQueue}{\overrightarrow{\Queue}}
\newcommand{\BwQueue}{\overleftarrow{\Queue}}

\newcommand{\IncomingArcs}{\delta^{+}}
\newcommand{\OutgoingArcs}{\delta^{-}}
\newcommand{\AllVertices}{\forall\Vertex \in \SetOfVertices}
\newcommand{\AllArcs}{\forall\Arc \in \SetOfArcs}
\newcommand{\AllLines}{\forall\Line \in \SetOfLines}
\newcommand{\AllStops}{\forall\Stop \in \SetOfStops}

\title{\large Infrastructure Planning for Inductive Charging in Electrified Shuttle Systems}

\author[1]{\normalsize Paul Bischoff}
\author[2]{\normalsize Salma Hammani}
\author[3]{\normalsize Maximilian Schiffer}

\affil[1]{\small School of Management, Technical University of Munich, Munich, Germany\par
\scriptsize paul.bischoff@tum.de}

\affil[2]{\small École des Ponts ParisTech, Paris, France\par
\scriptsize salma.hammani@outlook.fr}

\affil[3]{\small School of Management \& Munich Data Science Institute, Technical University of Munich, Munich, Germany \par
\scriptsize schiffer@tum.de}

\date{}

\begin{abstract}
\begin{singlespace}
{\small\noindent In response to climate goals, growing environmental awareness, and financial incentives, municipalities increasingly seek to electrify public transportation networks. We study the problem of locating stationary and dynamic inductive charging stations for electric vehicles (EVs), allowing detours from fixed transit routes and schedules. Dynamic charging, which enables energy transfer while driving, reduces space usage in dense urban areas and lowers vehicle idle times. We formulate a cost-minimization problem that considers both infrastructure and operational costs and propose an Iterated Local Search (ILS) algorithm to solve instances of realistic size. Each configuration requires solving a decomposed subproblem comprising multiple resource-constrained shortest-path problems. For this, we employ a bi-directional label-setting algorithm with lazy dominance checks based on local bounds. On adapted benchmark instances, our approach outperforms a commercial solver by up to 60\% in solution quality. We further apply our method to a real-world case study in Hof, Germany. Results indicate that, under current cost structures calibrated from a test track in Bad Staffelstein, dynamic inductive charging is not yet cost-competitive with stationary alternatives. We quantify the value of allowing detours at up to 3.5\% of the total system cost and show that integrating photovoltaics with decentralized energy storage can yield savings exceeding 20\%.
\\
}
{\footnotesize\noindent \textbf{Keywords:} Iterative Local Search, Resource Constrained Shortest Path, Charging Station Location, Wireless Power Transfer, Dynamic Charging}
\end{singlespace}
\end{abstract}

\maketitle

\section{Introduction}
Road-based public transportation systems—such as buses and shuttles—are the backbone of urban mobility. They offer affordable, flexible, and sustainable travel while helping to alleviate congestion and reduce emissions \citep{LiuSong2017}. Compared to rail-based systems, buses require minimal infrastructure investment, making them the dominant mode of travel in urban transit networks \citep[e.g.,][]{MVG2023}. However, most fleets still operate with combustion engines and are significant contributors to CO$_2$ emissions \citep{EUCommission2023}.

The European Commission’s Green Deal mandates that all new city buses be zero-emission by 2030 \citep{EUCommission2023}. Battery-electric vehicles (EVs) are widely considered the most viable path forward due to their low or zero tailpipe emissions \citep{MITClimatePortal2023, MVG2023}. However, their adoption introduces operational challenges: limited driving ranges require intermediate recharging, resulting in downtime, space occupancy, and manual interventions such as plug-in procedures.

Dynamic wireless charging—allowing vehicles to recharge while driving—addresses many of these concerns. It enables smaller, lighter, and more cost-efficient batteries, increasing energy efficiency and vehicle capacity \citep[cf.][]{LiuSong2017}.

This research is motivated by the ongoing \gls{acr:milas} project, which explores the deployment of autonomous, inductively charged shuttles in rural municipalities. As part of a real-world testbed in Bad Staffelstein, Germany, two autonomous electric shuttles are operated and evaluated under realistic service conditions. These shuttles are equipped for both stationary and dynamic inductive charging, enabling them to recharge without manual intervention. This functionality supports fully autonomous operations and extends driving ranges by allowing vehicles to charge during their scheduled trips. Furthermore, the project integrates renewable energy by connecting inductive charging stations to decentralized photovoltaic systems with local energy storage. These storage systems supply energy when solar reserves are sufficient, while fallback to the main grid occurs otherwise, introducing variable, time-dependent energy prices based on local supply availability.

\begin{figure}[t]
    \centering
    \begin{subfigure}[b]{0.48\textwidth}
        \centering
        \includegraphics[width=0.95\textwidth]{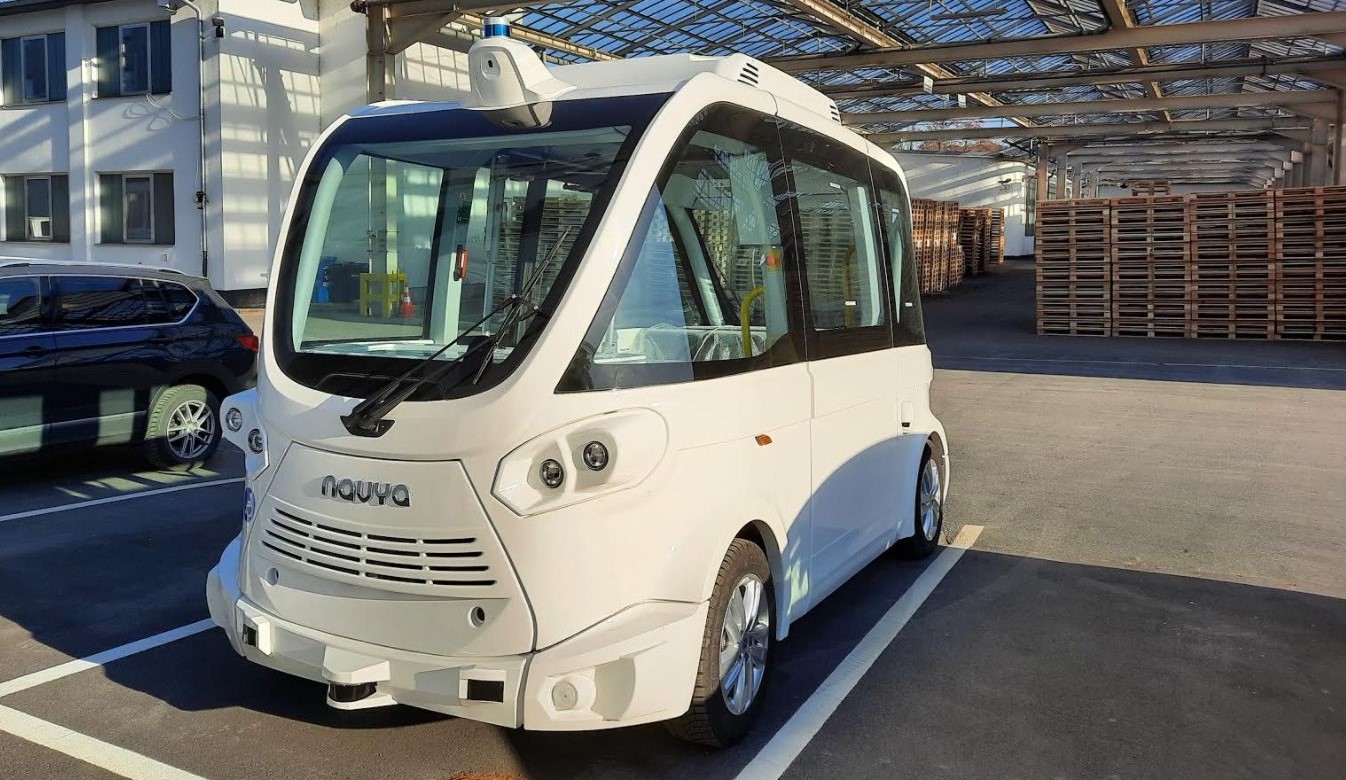}
        \caption{Electric and autonomous bus}
        \label{fig:intro-concepts-left}
    \end{subfigure}\hfill
    \begin{subfigure}[b]{0.48\textwidth}
        \centering
        \includegraphics[width=\textwidth]{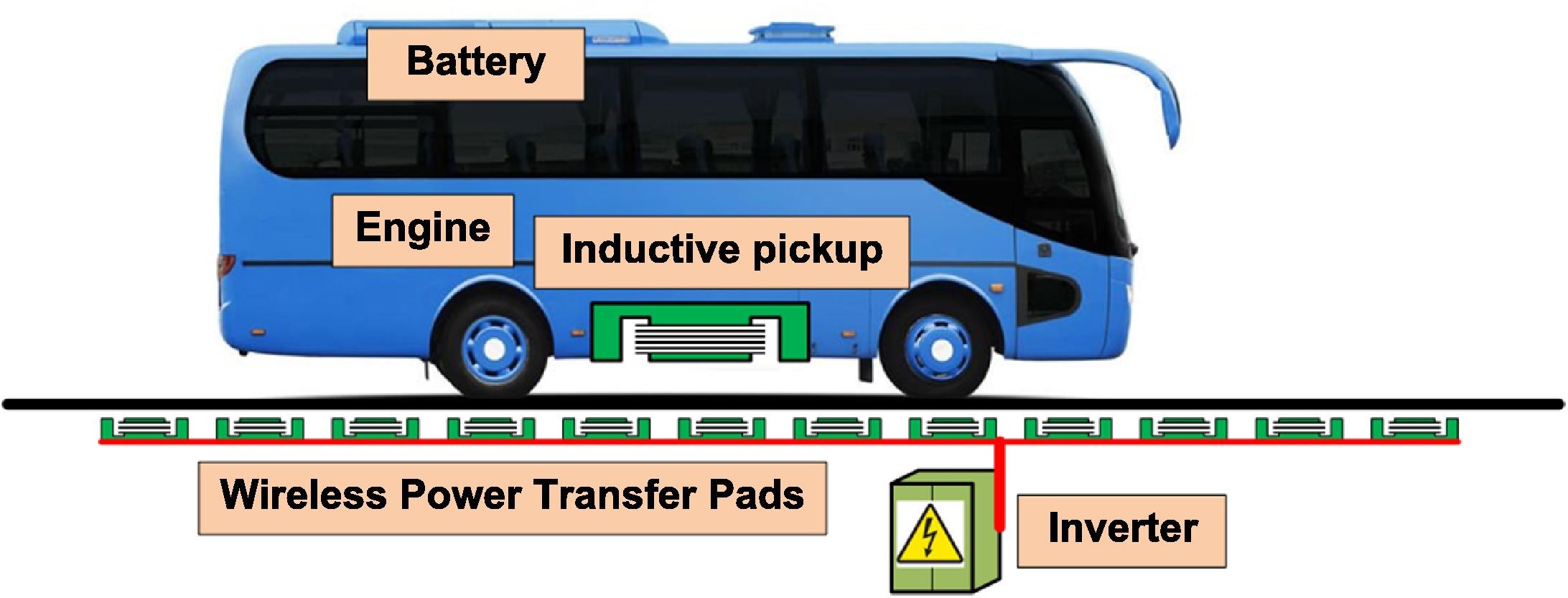}
        \caption{Inductive charging \citep{LiuSong2017}}
        \label{fig:intro-concepts-right}
    \end{subfigure}
    \caption{Inductive charging of autonomous passenger vehicles}
    \label{fig:intro-concepts}
\end{figure}

We study the optimal placement of stationary and dynamic inductive charging infrastructure for such shuttle systems. While stop sequences and schedules are fixed, vehicles may deviate from their shortest paths within time constraints to access charging facilities. Our objective is to minimize total infrastructure and operational costs.

\subsection{State of the art}
The location of charging or refueling infrastructure for electric vehicle fleets has been studied extensively across transportation and energy systems research. Existing approaches can broadly be categorized into three model classes: (i) node-based location models, (ii) path-based location models, and (iii) location-routing problems. Node-based models assume that recharging demand is aggregated at nodes and minimize cost or distance metrics associated with serving this demand \citep[e.g.,][]{TuLiEtAl2016}. Path-based models, by contrast, consider OD-pair-level flows and focus on locating facilities along fixed or flexible routes. The foundational Flow Capturing Location Allocation Model \citep{Hodgson1990} and the subsequent Flow-Refueling Location Model (FRLM) by \cite{KubyLim2005} laid the groundwork by introducing constraints on path lengths and number of stations. These models have since been generalized to allow for arc-based locations \citep{KubyLim2007}, station capacities \citep{UpchurchKubyEtAl2009}, and deviations from the shortest path \citep[e.g.,][]{KimKuby2012, LiHuang2014, YildizArslanEtAl2016}. While early works assumed a fixed set of candidate paths, later formulations such as \cite{YildizArslanEtAl2016} used column generation to endogenously construct deviating paths. A third stream of work focuses on location-routing problems, motivated by applications where operators centrally manage fleets. In such settings, paths are decision variables, and routing and charging placement are integrated. \cite{SchifferWalter2017A} initiated this direction for plug-in EVs, and \cite{SchifferWalter2018} extended it to intra-route charging via Adaptive Large Neighborhood Search.

We focus our review specifically on works involving dynamic charging of EV fleets and highlight research on electrified bus operations using stationary and dynamic inductive charging. We differentiate between \textit{macroscopic approaches}, which model infrastructure at a network or region level, and \textit{microscopic approaches}, which focus on vehicle-specific operations and fleet dispatch. In the macroscopic literature, most works assume fixed OD-pairs and route structures. For example, \cite{MajhiRanjitkarEtAl2022} account for congestion and dynamic state-of-charge evolution along fixed routes, while \cite{SongCheng2024} remove the fixed-path assumption and develop a column generation framework to determine the optimal location of dynamic charging segments. \cite{MubarakUesterEtAl2021} propose a Benders decomposition for infrastructure placement in multi-class EV networks. Other works focus on user equilibrium behavior. \cite{ChenHeEtAl2016}, \cite{RiemannWangEtAl2015}, and \cite{LiuWangEtAl2017} incorporate stochastic user routing and vehicle choice models to study station demand under charging preferences. \cite{SunChenEtAl2020} embed these choices in a bi-level optimization model. Recent advances leverage time-expanded or space-time network models to capture heterogeneity in vehicle types and infrastructure. Notably, \cite{ZengXieEtAl2024} apply dynamic programming for vehicle-specific decisions, while \cite{SongChengEtAl2023} propose a Lagrangian relaxation with decomposition by OD pair to jointly deploy stationary and dynamic charging systems.

In contrast, microscopic approaches focus on the operational integration of inductive charging technologies into public transportation systems. This literature builds heavily on early research around the Online Electric Vehicle (OLEV) system developed by KAIST. \cite{KoJang2013} initiated this line of work by optimizing the design of a dynamic charging bus system. Follow-up studies developed optimization-based methods for charging infrastructure design and battery sizing \citep{KoJangEtAl2015}, extended the models to overlapping routes \citep{LiuSongEtAl2017, HwangJangEtAl2018}, and examined scheduling under vehicle homogeneity versus heterogeneity \citep{AlwesabiWangEtAl2020}. Some works introduce uncertainty by applying robust optimization techniques to travel times and energy consumption \citep{LiuSong2017}. More recent studies evaluate dynamic wireless charging feasibility in other use cases, such as airport bus operations \citep{HelberBroihanEtAl2018} and logistics fleets \citep{SchwerdfegerBockEtAl2022}, often using simulation-based analysis. However, many of these models assume fixed vehicle routes, do not account for state-of-charge dynamics, or ignore charging opportunities via detours. \cite{ChenYinEtAl2018} compare dynamic, stationary, and battery-swapping approaches in cost analyses for single-route networks, while \cite{BroihanNozinskiEtAl2022} introduce vehicle heterogeneity in simulation-based cost assessments. A more integrated view is taken by \cite{IliopoulouKepaptsoglou2019}, who co-optimize transit route design and charging infrastructure using a bi-level model, albeit without variable energy pricing or idling for charge. To the best of our knowledge, \cite{LiHeEtAl2024} offer the only existing microscopic model that integrates variable energy prices with dynamic charging and vehicle scheduling decisions. However, their focus remains on simplified network structures and does not support combined deployment of stationary and dynamic infrastructure.

\subsection{Contribution}
We propose a novel planning problem for determining the cost-optimal deployment of heterogeneous charging infrastructure—comprising both stationary and dynamic inductive charging—for electric shuttle fleets operating on fixed stop sequences and timetables. Our model is the first to integrate variable electricity prices induced by local photovoltaic storage systems and to allow detours from shortest paths to enable infrastructure sharing across partially overlapping routes. This combination captures important real-world trade-offs between cost, routing flexibility, and energy autonomy.

We formulate the problem as a mixed-integer program and develop a scalable metaheuristic based on Iterated Local Search (ILS). The core subproblem, which determines vehicle-level charging and routing decisions, is solved via dynamic programming using a bi-directional label-setting algorithm with lazy dominance checks.

We demonstrate that our approach improves over a warm-started commercial solver by up to 60\% on benchmark instances. In a case study for the city of Hof, Germany, we show that dynamic inductive charging is technically feasible but not economically viable under current cost conditions. Nonetheless, permitting detours reduces system costs by up to 3.5\%, and integrating decentralized solar-powered energy storage reduces costs by more than 20\%.

\def\checkmark{\tikz\fill[scale=0.4](0,.35) -- (.25,0) -- (1,.7) -- (.25,.15) -- cycle;} 
\newcommand{\rt}[1]{\rotatebox{90}{#1}}
\newcommand{\lb}[1]{\pbox{12mm}{#1}}
\begin{table}[!tb]
\footnotesize
\begin{tabular}{lccccccc}
\hline
\textbf{} & \textbf{{[}1{]}} & \textbf{{[}2{]}} & \textbf{{[}3{]}}  & \textbf{{[}4{]}} & \textbf{{[}5{]}} & \textbf{{[}6{]}} & \textbf{\begin{tabular}[c]{@{}c@{}}Our \\ work\end{tabular}} \\ \hline
Stationary and Dynamic Charging & - & - & $\checkmark$ & - & - & -  & $\checkmark$ \\
Route Deviations & - & -  & $\checkmark$ & - & - & - &  $\checkmark$\\
Multiple Routes & - & $\checkmark$ & $\checkmark$ & $\checkmark$ & - & $\checkmark$ & $\checkmark$ \\
Variable Prices & - & - & - & - & - & $\checkmark$ & $\checkmark$ \\
\hline
\end{tabular}
\caption{Locating wireless charging facilities for electrified road-based public transportation systems. \textmd{Indices [1] to [6] signify publications as follows: [1] \cite{KoJang2013, KoJangEtAl2015}; [2] \cite{LiuSong2017,LiuSongEtAl2017, HwangJangEtAl2018, AlwesabiWangEtAl2020}; [3] \cite{IliopoulouKepaptsoglou2019}; [4] \cite{HelberBroihanEtAl2018,BroihanNozinskiEtAl2022}; [5] \cite{ChenYinEtAl2018}; [6] \cite{LiHeEtAl2024}}}
\label{tbl:contribution}
\end{table}

The remainder of this paper is structured as follows. Section~\ref{sec:prob-description} details the problem formulation. Section~\ref{sec:alg-framework} describes our ILS-based solution method. The experimental design is presented in Section~\ref{sec:exp-design}, followed by results in Section~\ref{sec:results}. We conclude in Section~\ref{sec:conclusion}.

\section{Problem Description} \label{sec:prob-description}
We consider a fleet of electric vehicles $\Line \in \SetOfLines$ that charge inductively via embedded coils in magnetically conductive concrete for wireless power transfer. Inductive battery charging is possible in a stationary setting and in a dynamic setting while driving. Every vehicle $\Line \in \SetOfLines$ operates on its individual route between consecutive stop locations $\Stop \in \SetOfStops$. We assume the fleet to start and end their routes from a central depot that we represent as $\StartDepot$ and $\EndDepot$, to share a homogeneous \gls{acr:soc} $\InitialSoC$ when starting their operations, and to operate within homogeneous \gls{acr:soc} limits denoted by $\MinSoC$ and $\MaxSoC$. We denote the vehicle routes by an ordered sequence of locations $\StopsOfLine := \langle \StartDepot, 1, \dots, n, \EndDepot \rangle$ with $n$ being the number of stops on the vehicle's route. Every location $\StopOfLine \in \StopsOfLine$ on a vehicle's route is associated with a time window $[\TwBegin^{\Line}_{\StopOfLine}, \TwEnd^{\Line}_{\StopOfLine}]$ within which $\Line$ has to serve $\StopOfLine$ in order to meet a timetable that is predetermined. 

\begin{figure}[!tb]
  \centering
  \begin{adjustbox}{center}
    \begin{minipage}{0.55\textwidth}
      \centering
      \def\svgwidth{\linewidth}
      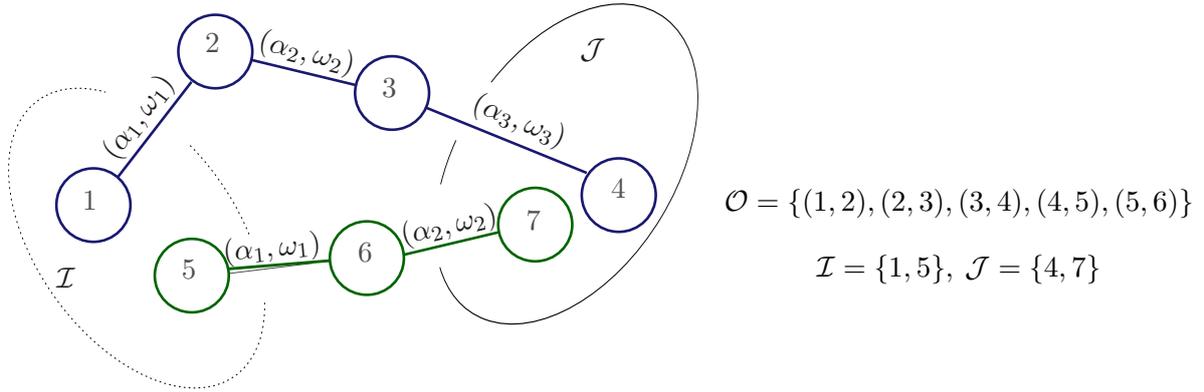
    \end{minipage}
    \hfill
    \begin{minipage}{0.4\textwidth}
      \vspace{0.5cm} 
      \[
        \SetOfSegments = \{ (1, 2), (2, 3), (3, 4), (4, 5), (5, 6)\}
      \]
      \[
        \SetOfSegmentStarts = \{ 1, 5\}, \medspace \SetOfSegmentEnds = \{ 4, 7\}
      \]
    \end{minipage}
  \end{adjustbox}
  \caption{Exemplary setting with two dynamic charging stations (displayed in green and blue) consisting of multiple segments.}
  \label{fig:segment-setting}
\end{figure}

Electric vehicles can recharge inductively either at stationary charging stations located at candidate sites $\Station \in \SetOfStaticStations$, or at dynamic charging stations embedded in the roadway at candidate locations $\Station \in \SetOfDynamicStations$. We define the complete set of candidate charging stations as $\SetOfStations := \SetOfStaticStations \cup \SetOfDynamicStations$. Each station $\Station$ is characterized by a charging rate $\StationChargingRate_{\Station}$ and a construction cost $\InvestmentCost_{\Station}$, both of which may vary across locations to capture heterogeneity in local conditions such as available space or infrastructure requirements. Every dynamic charging station consists of an ordered sequence of consecutive segment indices $\SetOfSegments_\Station := \langle 1, \dots, m \rangle$ with $m \ge 1$. With every index $\SegmentIndex \in \SetOfSegments_{\Station}$ we associate a segment $(\SegmentStart_\SegmentIndex, \SegmentEnd_\SegmentIndex)$ such that the start point of a segment equals the endpoint of its predecessor $\SegmentStart_{\SegmentIndex+1} = \SegmentEnd_{\SegmentIndex}, \medspace \SegmentIndex = 1,\dots, m-1$. In this context, we slightly abuse notation by referring to $(\SegmentStart_\SegmentIndex, \SegmentEnd_\SegmentIndex) \in \SetOfSegments_{\Station}$ when convenient. Furthermore, $\SetOfSegments := \bigcup_{\Station \in \SetOfDynamicStations} \{ (\SegmentStart_\SegmentIndex, \SegmentEnd_\SegmentIndex): \SegmentIndex \in \SetOfSegments_\Station \}$ yields the set of all dynamic charging segments, and accordingly $\SetOfSegmentStarts := \bigcup_{\Station \in \SetOfDynamicStations} \SegmentStart_{1}$, $\SetOfSegmentEnds := \bigcup_{\Station \in \SetOfDynamicStations} \SegmentEnd_{m}$ denote the set of the first segment start points and last segment end points, respectively. Figure~\ref{fig:segment-setting} shows an example with two stations to demonstrate the introduced notation. While the construction cost $\InvestmentCost_{\Station}, \Station \in \SetOfDynamicStations$ is associated with the dynamic charging station and covers fixed costs for power inverters and their infrastructure, every additional segment adds a variable cost term $\InvestmentCost_{\UnbracedSegment}$ that depends linearly on the segment's length. 

We model our problem on a directed graph $\Graph\!=\!(\SetOfVertices,\!\SetOfArcs)$ with a vertex set \mbox{$\SetOfVertices\!:=\!\SetOfStops\!\cup\!\SetOfStaticStations\!\cup\!\SetOfSegmentStarts\!\cup\!\SetOfSegmentEnds\!\cup\!\{\StartDepot, \EndDepot\}$} representing the set of stops, the stationary charger candidate locations, and the depot locations that serve as starting and end point of the vehicles' operations. 

Moreover, we construct arcs $(i,j) \in \SetOfArcs$ as follows. First, we add arcs $\Arc \in \SetOfStops \cup \{\StartDepot\} \times \SetOfStops \cup \{\EndDepot\}$ where $i \neq j$. Second, we connect stationary candidate charging stations to stops. More specifically, we add arcs $\Arc \in \SetOfStops \cup \{\StartDepot\} \times \SetOfStaticStations$ and $\Arc \in \SetOfStaticStations \times \SetOfStops \cup \{\EndDepot\}$. Moreover, we add arcs $\Arc \in \SetOfSegments$ representing the dynamic charging candidate segments. Finally, we connect dynamic charging stations to stops by adding arcs $\Arc \in \SetOfStops \cup \{\StartDepot\} \times \SetOfSegmentStarts$ and $\Arc \in \SetOfSegmentEnds \times \SetOfStops \cup \{\EndDepot\}$. Every arc $\Arc \in \SetOfArcs$ is associated with a  travel time $\TravelTime_{\UnbracedArc}$, and an energy consumption $\ArcConsumption_{\UnbracedArc}$. Furthermore, we denote the set of incoming arcs at a vertex $\Vertex \in \SetOfVertices$ as $\IncomingArcs(\Vertex) \subseteq \SetOfArcs$, and the set of outgoing arcs as $\OutgoingArcs(\Vertex) \subseteq \SetOfArcs$ respectively.

\noindent \textbf{Solution: } A solution $\Solution$ consists of three major components: the network design, the vehicle routes such that all stops are serviced according to their service requirements, and the vehicles' recharging profile in order to keep the fleet operational. 

We encode the design decision as follows: let $\mathbf{\StaticChargerVar} \in \{0,1\}^{|\SetOfStaticStations|}$ be a binary vector where each component corresponds to a stationary candidate charging station. If $\StaticChargerVar_\Station = 1$, station $\Station \in \SetOfStaticStations$ is part of the solution. Furthermore, let $\mathbf{\DynamicChargerVar} \in \{0,1\}^{|\SetOfSegments|}$ be a binary vector where each component corresponds to a dynamic candidate charging segment. Binary vector $\mathbf{\DynamicChargerVar}$ encodes which dynamic charging segments are part of the solution. Lastly, binary vector $\mathbf{\TransformerVar} \in \{0,1\}^{|\SetOfDynamicStations|}$ encodes the construction of power inverters that are connecting the constructed dynamic charging segments of every dynamic station with the power grid.

To encode vehicle routes, we define a binary matrix $\mathbf{X} \in \{0,1\}^{|\SetOfArcs| \times |\SetOfLines|}$, where every column $\mathbf{\ArcVar}^{\Line}$ encodes a vehicle route by $\ArcVar^{\Line}_{\UnbracedArc} = 1$ if vehicle $\Line$ traverses $\Arc$, and $\ArcVar^{\Line}_{\UnbracedArc}=0$ otherwise. Moreover, a solution yields a matrix $\mathbf{T} \in \mathbb{R}^{|\SetOfVertices| \times |\SetOfLines|}_{\ge 0}$ in which every element encodes the vehicle's departure time at the respective vertex. Finally, a solution entails a charging profile represented by the battery \gls{acr:soc} at the time of departure from any vertex $\Vertex \in \SetOfVertices$, modeled as a matrix $\mathbf{R} \in \mathbb{R}^{|\SetOfVertices| \times |\SetOfLines|}_{\ge 0}$ where every column $\StateOfChargeVar^{\Line}$ represents the \gls{acr:soc} profile of the respective vehicle.

\begin{figure}[!t]
    \centering
    \resizebox{\textwidth}{!}{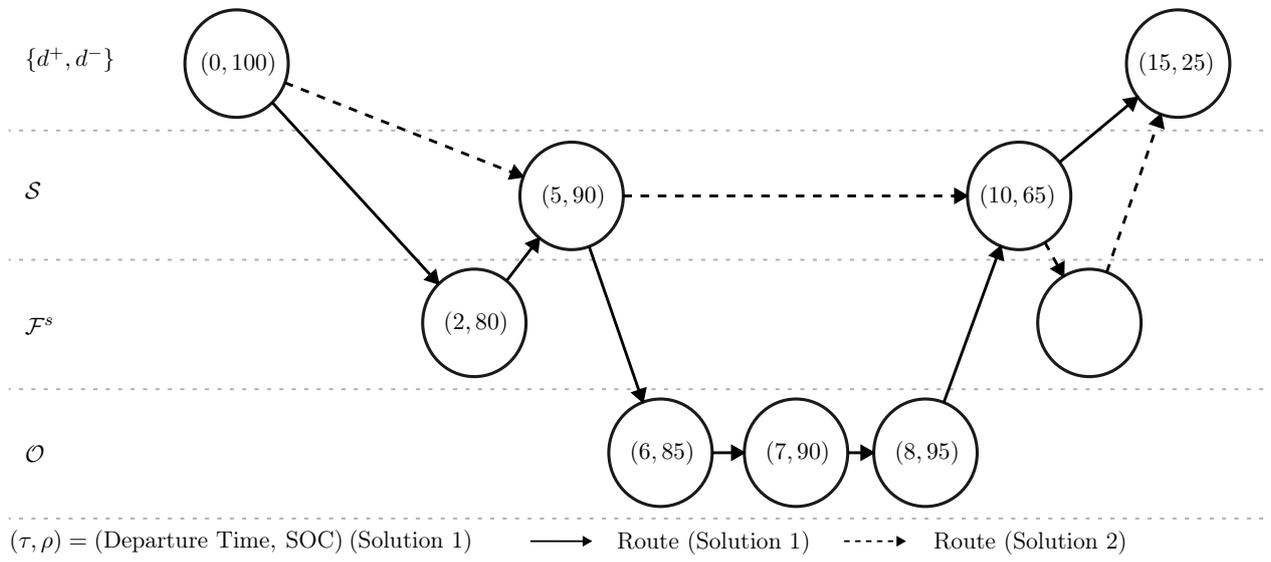}
    \caption{Illustrative example with $|\SetOfVehicles|=1$. Solution $1$ routes the vehicle to a stationary charging station before servicing the first stop. Subsequently the vehicle traverses a dynamic charging station before servicing the next stop and returning to the depot at time $15$ with an \gls{acr:soc} of $25$. The alternative solution $2$ services the stops first and only visits a stationary charging station before returning to the depot. Note that the departure time and \gls{acr:soc} profile for this solution are not displayed.}
    \label{fig:graph-example}
\end{figure}

To summarize, we search for solutions of the form $\Solution = (\mathbf{\StaticChargerVar}, \mathbf{\TransformerVar}, \mathbf{\DynamicChargerVar}, \mathbf{X}, \mathbf{T}, \mathbf{R})$. Figure~\ref{fig:graph-example} shows an example solution for a toy example together with a potential alternative route. For the sake of readability, we do not display departure times and battery \gls{acr:soc} for the alternative solution.

\noindent \textbf{Constraints: } In a feasible solution, all vehicle routes $\mathbf{\ArcVar}^{\Line}$ must service their stop sequences $\StopsOfLine$ in the given order such that $\ArrivalTimeVar^{\Line}_{\StopOfLine} \leq \ArrivalTimeVar^{\Line}_{\StopOfLine+1}, \medspace \Line \in \SetOfLines, \StopOfLine \in \StopsOfLine \setminus \{\EndDepot\}$ and such that they respect the service time windows, i.e., $\TwBegin^{\Line}_{\StopOfLine} \le \ArrivalTimeVar^{\Line}_{\StopOfLine} \le \TwEnd^{\Line}_{\StopOfLine}, \medspace \Line \in \SetOfLines, \StopOfLine \in \StopsOfLine$. Moreover, the \gls{acr:soc} of every vehicle must remain between the given limits. Formally, we require $\MinSoC \leq \StateOfChargeVar^{\Line}_{\Vertex} \leq \MaxSoC, \medspace \Line \in \SetOfLines, \Vertex \in \SetOfVertices$.

\noindent \textbf{Objective: } We search for a cost-optimal solution such that the sum of infrastructure costs $\InfrastructureCost$ and operational costs $\RoutingCost$ is minimal:
\begin{equation*}
    \min\limits_{(\mathbf{\StaticChargerVar}, \mathbf{\TransformerVar}, \mathbf{\DynamicChargerVar}, \mathbf{X}, \mathbf{T}, \mathbf{R})} \InfrastructureCost(\mathbf{\StaticChargerVar}, \mathbf{\TransformerVar}, \mathbf{\DynamicChargerVar})  + \RoutingCost(\mathbf{X}, \mathbf{T}, \mathbf{R}).
\end{equation*}

The infrastructure costs contain the fixed costs for constructing stationary and dynamic charging stations as well as the variable cost per dynamic charging segment.
\begin{equation} \label{eq:infrastructure-cost}
    \InfrastructureCost(\mathbf{\StaticChargerVar}, \mathbf{\TransformerVar}, \mathbf{\DynamicChargerVar}) = \sum_{\Station \in \SetOfStaticStations} \InvestmentCost_{\Station} \StaticChargerVar_{\Station} + \sum_{\Station \in \SetOfDynamicStations} \InvestmentCost_{\Station} \TransformerVar_{\Station} + \sum_{\Segment \in \SetOfSegments} \InvestmentCost_{\UnbracedSegment} \DynamicChargerVar_{\UnbracedSegment}
\end{equation}

The operational costs comprise costs for energy consumption and re-charging while operating:
\begin{equation} \label{eq:operational-cost}
    \RoutingCost(\mathbf{X}, \mathbf{T}, \mathbf{R}) = \sum_{\Line \in \SetOfVehicles} \sum_{\Arc \in \SetOfArcs} \EnergyPrice_{\texttt{c}} \, \ArcConsumption_{\UnbracedArc} \, \ArcVar^{\Line}_{\UnbracedArc} + \sum_{\Line \in \SetOfVehicles} \sum_{\Arc \in \SetOfArcs} \EnergyPrice_{\texttt{r}} \, (\StateOfChargeVar^{\Line}_{j} - \StateOfChargeVar^{\Line}_{i} + \ArcConsumption_{\UnbracedArc}) \, \ArcVar^{\Line}_{\UnbracedArc} 
\end{equation}

with energy prices $\EnergyPrice_{\texttt{c}}, \EnergyPrice_{\texttt{r}} > 0$ for consumption and recharging respectively. In ~\eqref{eq:operational-cost}, the first term accumulates the cost for consumed energy while driving, and the second term accumulates the cost for recharged energy derived from the solution's charge profile.

\section{Algorithmic Framework} \label{sec:alg-framework}
We propose an \gls{acr:ils} \citep[cf. e.g.,][]{StützleRuiz2018} to solve the planning problem that we introduced in Section~\ref{sec:prob-description} and elaborate on its components in this Section. We initialize the \gls{acr:ils} with a randomized approach to determine a feasible initial solution $\Solution_{0}$ in which we sample, with equal probabilities, $\ILSInitStepSize$ elements from the set of stations $\SetOfStations$ and add them to the configuration. 

Formally, we set $\StaticChargerVar_{\Station}=1$ for a drawn $\Station \in \SetOfStaticStations$, and $\DynamicChargerVar_{\UnbracedSegment} = 1, \medspace \forall \Segment \in \SetOfSegments_{\Station}$ for a drawn $\Station \in \SetOfDynamicStations$. We repeat until we can prove feasibility for the resulting configuration (see Section~\ref{sec:subproblem}). 

The remainder of this section details our algorithmic framework. Specifically, Subsection~\ref{sec:ils} explains the components of the \gls{acr:ils} and Section~\ref{sec:subproblem} yields a description of the subproblem decomposing by vehicle into multiple \glspl{acr:rcspp}. We use the subproblem to evaluate configurations by determining the operational costs \eqref{eq:operational-cost}, or label the configurations as infeasible alternatively. 

\subsection{Iterated Local Search} \label{sec:ils}
Algorithm~\ref{alg:ils} shows the pseudocode of our \gls{acr:ils} framework. The algorithm starts by initializing the maximum perturbation strength $\PerturbationStrength_{\texttt{max}}$ to the number of stationary charging stations and dynamic charging segments in the initial solution~(l.~\ref{l:init}). Afterwards it performs an initial call to a sample-based local search~(l.~\ref{l:init-ls}). The algorithm keeps track of the number of performed iterations without an update to the incumbent $\ILSIterations$ and steers the amount of exploration with the strength parameter $\PerturbationStrength^{\prime}$ accordingly. We initialize these parameters at their starting values (l.~\ref{l:init-paras}). The \gls{acr:ils} terminates when the algorithm reaches a time limit~(l.~\ref{l:time-limit}). In every iteration, the algorithm adjusts the strength parameter $\PerturbationStrength^{\prime} \geq 1$ that determines the extend of the randomized exploration in the procedure. The algorithm increases the strength by $\PerturbationStrength$ in every iteration after more than $\ILSIterations_{\texttt{max}}$ iterations without updates to the incumbent have passed~(l.~\ref{l:increase-strength}). The \gls{acr:ils} iterates between perturbing the best found solution~(l.~\ref{l:perturb}) and exploring the stochastic neighborhood of the perturbed solution via the sample-based local search~(l.~\ref{l:loop-ls}). The algorithm's acceptance criterion compares objective values~(l.~\ref{l:accept}) and accepts solutions that improve the current best found solution with a small tolerance $\epsilon \ll 1$. Once the algorithm finds a new incumbent, it resets the perturbation strength to its baseline value $\PerturbationStrength$~(l.~\ref{l:reset-paras}), and increases the counter of iterations without update to the incumbent otherwise~(l.~\ref{l:cnt-iters}).

\begin{algorithm}[!t]
\caption{\textsc{IteratedLocalSearch ($\Solution_0:=(\mathbf{\StaticChargerVar}_{0}, \mathbf{\TransformerVar}_{0}, \mathbf{\DynamicChargerVar}_{0}, \mathbf{X}_{0}, \mathbf{T}_{0}, \mathbf{R}_{0})$, $\PerturbationStrength$, $\ILSIterations_{\texttt{max}}$)}}
\label{alg:ils}
\begin{algorithmic}[1]
    \State $\PerturbationStrength_{\texttt{max}} \gets \mathbf{1}^\text{T} \mathbf{\StaticChargerVar_{0}} + \mathbf{1}^\text{T} \mathbf{\DynamicChargerVar_{0}}$ \label{l:init}
    \State $\Solution^* \gets \textsc{LocalSearch}(\Solution_0)$ \label{l:init-ls}
    \State $\ILSIterations \gets 0, \medspace \PerturbationStrength^{\prime} \gets \PerturbationStrength$ \label{l:init-paras} \Comment{Initialize iteration counter at zero and strength at baseline} 
    \While{$\textbf{not} \medspace \textsc{TimeOut}$} \label{l:time-limit}
        \If{$\ILSIterations \ge \ILSIterations_{\texttt{max}}$} \Comment{Iterations w/o new incumbent exceeds threshold}
            \State $\PerturbationStrength^{\prime} \gets \min \{\PerturbationStrength^{\prime} + \PerturbationStrength, \PerturbationStrength_{\texttt{max}} \}$ \label{l:increase-strength} \Comment{Increase strength for more exploration}
        \EndIf
        \State $\Solution^{\prime} \gets \textsc{Perturbation}(\Solution^*, \PerturbationStrength^{\prime})$ \label{l:perturb} \Comment{Contains repeated calls to the subproblem}
        \State $\Solution^{\prime \prime} \gets \textsc{LocalSearch}(\Solution^{\prime})$ \label{l:loop-ls} \Comment{Contains repeated calls to the subproblem}
        \If{$\textsc{EvalSubproblem}(\Solution^{\prime \prime}) \leq \textsc{EvalSubproblem}(\Solution^*) + \epsilon$} \label{l:accept}
            \State $\Solution^* \gets \Solution^{\prime \prime}$ \Comment{New incumbent found}
            \State $\ILSIterations \gets 0, \medspace \PerturbationStrength^{\prime} \gets \PerturbationStrength$ \label{l:reset-paras} \Comment{Reset parameters to reset exploration}
        \Else
            \State $\ILSIterations \gets +1$ \label{l:cnt-iters}
        \EndIf
    \EndWhile
    \State \Return $\Solution^{*}$
\end{algorithmic}
\end{algorithm}

\noindent \textbf{Sample-based local search:} We apply a sample-based local search to find local minima within randomized neighborhoods of a given solution. Algorithm~\ref{alg:local-search} shows a pseudocode describing the sample-based local search. We initialize the upper bound on the local minimum with the cost of the solution $\Solution := (\mathbf{\StaticChargerVar}, \mathbf{\TransformerVar}, \mathbf{\DynamicChargerVar}, \mathbf{X}, \mathbf{T}, \mathbf{R})$ that we pass to the local search~(l.~\ref{l:init-lm}) and add a constant of $2 \, \epsilon$ to ensure that we enter the first iteration of the subsequent while loop. The main body~(l.~\ref{l:while-loop}-\ref{l:end-while-loop}) of the local search applies the different operators sequentially and continues as long as the objective can be further decreased by a small $\epsilon \ll 1$. The operators explore sample-based neighborhoods of a given solution in order to keep the number of evaluated configurations tractable. We choose a homogeneous parameter $\NeighborhoodSize \ge 1$ that controls the size of the sample-based neighborhoods across all operators. Note that if no improving solution exists in a sample-based neighborhood defined by a local search operator, the passed solution itself is a minimum. Consequently, all local search operators return the original solution in this case. 

To define the search operators, we use the following set notation to denote charging stations that are part of a given configuration: $\SetOfStaticStations_{\mathbf{y}=1} := \{ \Station \in \SetOfStaticStations: \StaticChargerVar_{\Station} = 1 \}$  denotes the stationary charging stations in the current configuration, and similarly $\SetOfDynamicStations_{\mathbf{\TransformerVar}=1} := \{ \Station \in \SetOfDynamicStations: \TransformerVar_{\Station} = 1 \}$  denotes the dynamic charging stations in the current configuration. Thus, $\SetOfStaticStations_{\mathbf{y}=0} = \SetOfStaticStations \setminus \SetOfStaticStations_{\mathbf{y}=1}$ and $\SetOfDynamicStations_{\mathbf{\TransformerVar}=0} = \SetOfDynamicStations \setminus \SetOfDynamicStations_{\mathbf{\TransformerVar}=1}$. Our sample-based local search operators consist of three groups. First, operators \gls{acr:strip-dyn} and \gls{acr:remove-stat} remove random charging stations from the given configuration. We designed these operators to find solutions in which the infrastructure costs~\eqref{eq:infrastructure-cost} outweigh the operational costs~\eqref{eq:operational-cost} or the given configuration contains redundant charging stations. Second, the operators \gls{acr:swap-to-dyn} and \gls{acr:swap-to-stat} exchange random charging stations in the given configuration for dynamic charging stations $\Station \in \SetOfDynamicStations_{\mathbf{\TransformerVar}=0} $, and stationary charging stations $\Station \in \SetOfStaticStations_{\mathbf{y}=0}$ respectively. These operators are capable to evaluate the trade-off between stationary and dynamic charging and shape the solutions accordingly. Third, operators \gls{acr:add-dyn} and \gls{acr:add-stat} add charging stations to the given configuration. We designed these operators to find solutions in which the operational costs~\eqref{eq:operational-cost} outweigh the infrastructure costs~\eqref{eq:infrastructure-cost}. Appendix~\ref{app:sb-ls-operators} provides a detailed description of all operators.

\begin{algorithm}[!t]
    \caption{$\textsc{Sample-based LocalSearch}(\Solution = (\mathbf{\StaticChargerVar}, \mathbf{\TransformerVar}, \mathbf{\DynamicChargerVar}, \mathbf{X}, \mathbf{T}, \mathbf{R}))$}
    \label{alg:local-search}
    \begin{algorithmic}[1]
        \State $\textsc{LocalMinimum} \gets \textsc{EvalSubproblem}(\Solution) + 2\epsilon$ \label{l:init-lm}
        \While{$\textsc{EvalSubproblem}(\Solution) + \epsilon \leq \textsc{LocalMinimum}$} \label{l:while-loop}
            \State $\textsc{LocalMinimum} \gets \textsc{EvalSubproblem}(\Solution)$
            \State $\Solution \gets \glsentrylong{acr:strip-dyn}(\Solution)$
            \State $\Solution \gets \glsentrylong{acr:remove-stat}(\Solution)$
            \State $\Solution \gets \glsentrylong{acr:swap-to-dyn}(\Solution)$
            \State $\Solution \gets \glsentrylong{acr:swap-to-stat}(\Solution)$
            \State $\Solution \gets \glsentrylong{acr:add-dyn}(\Solution)$
            \State $\Solution \gets \glsentrylong{acr:add-stat}(\Solution)$
        \EndWhile \label{l:end-while-loop}
        \State \Return $\Solution$
    \end{algorithmic}
\end{algorithm}

\noindent \textbf{Perturbation:} The perturbation procedure relies on a strength parameter $\PerturbationStrength^{\prime} \geq 1$ that determines the extend of the exploration in the current \gls{acr:ils} iteration. Algorithm~\ref{alg:perturbation} shows a pseudocode that details the perturbation. The perturbation starts by destroying the current solution~(l.~\ref{l:destroy}). For this purpose, the algorithm samples $\min \{\PerturbationStrength^{\prime}, |\SetOfStaticStations_{\mathbf{y}=1}| + |\SetOfDynamicStations_{\mathbf{\TransformerVar}=1}|\}$ from the current configuration and removes the sampled stations from the configuration. The resulting new configuration may be infeasible. Then, the perturbation algorithm enters a loop and adds charging stations to the configuration until the resulting configuration is feasible. In this context, the algorithm samples the stations from the set $\SetOfStaticStations_{\mathbf{y}=0} \cup \SetOfDynamicStations_{\mathbf{\TransformerVar}=0}$ with equal probabilities~(l.~\ref{l:repair-sample}) before adding them to the configuration~(l.~\ref{l:repair-add}). The procedure accepts and returns the resulting configuration~(l.~\ref{l:repair-accept},\ref{l:perturb-return}) if it is feasible. Otherwise, it increases $\PerturbationStrength^{\prime}$ by one~(l.~\ref{l:repair-increase-strengt}) and re-samples corresponding stations in a subsequent iteration.

\noindent \textbf{Discussion:} Three comments on the implementation of our algorithmic framework are in order. First, we avoid to repeatedly solve the same subproblem by storing configurations with their respective subproblem solution values in memory. Moreover, we maintain a cache with information on which configurations the algorithm already evaluated as infeasible. Second, we leverage the fact that removing a charging station from the configuration does not affect the solution of the decomposed subproblem of the respective vehicle if this vehicle does not charge at the respective charging station. We waive the evaluation of such decomposed subproblems and use the solution that we previously stored in memory. Third, we note that the scalability of our algorithmic framework is independent of the number of vehicles $|\SetOfVehicles|$ since the decomposable structure of the subproblem allows the parallel evaluation of the decomposed subproblems.

\begin{algorithm}[!t]
\caption{\textsc{Perturbation ($\Solution$, $\PerturbationStrength^{\prime}$)}}
\label{alg:perturbation}
\begin{algorithmic}[1]
    \State $\Solution^{\prime} \gets \textsc{Destroy}(\PerturbationStrength^{\prime})$ \label{l:destroy}
    \Comment{Remove up to $\PerturbationStrength^{\prime}$ stations from current configuration}
    \While{$\Solution^{\prime} \textsc{ Infeasible}$} \label{l:repair-loop}
        \State $\textsc{NewStations} \gets \textsc{SampleStations}(\PerturbationStrength^{\prime})$ \Comment{With equal probability}\label{l:repair-sample}
        \State $\Solution^{\prime \prime} \gets \textsc{AddStations}(\Solution^{\prime}, \textsc{NewStations})$ \label{l:repair-add}
        \If{$\Solution^{\prime \prime} \textsc{ Feasible}$} 
            \State $\Solution^\prime \gets \Solution^{\prime \prime}$ \label{l:repair-accept} \Comment{Feasible repaired solution found}
        \EndIf
        \State $\PerturbationStrength^{\prime} \gets +1$ \label{l:repair-increase-strengt} \Comment{Increase number of stations to restore feasibility}
    \EndWhile
    \State \Return $\Solution^{\prime}$ \label{l:perturb-return}
\end{algorithmic}
\end{algorithm} 

\subsection{Subproblem} \label{sec:subproblem}
The routing subproblem determines the cost-optimal routes between any consecutive pair of locations $\StopOfLine, \StopOfLine+1, \medspace \StopOfLine \in \StopsOfLine \setminus \{\EndDepot\}$ for every vehicle $\Line \in \SetOfLines$ given limited energy resources and respecting the given time windows. We solve the problem by partially time-expanding graph $\Graph$ into a vehicle-specific graph $\Graph^{\Vehicle}$ and solving a \gls{acr:rcspp} such that a solution to the subproblem can be represented as a path $\Path^{\Vehicle}$ of vertices in $\Graph^{\Vehicle}$. 

Although introduced with a constant energy price $\EnergyPrice_{\texttt{r}} > 0$, we can easily incorporate variable energy prices into our problem formulation by evaluating a given price curve at the respective departure times. Specifically, we define $\EnergyPrice_{\texttt{r}} : \mathbb{R}_{\geq 0} \mapsto \mathbb{R}_{\geq 0}$ where $\EnergyPrice_{\texttt{r}}$ maps each continuous time $\tau$ to an energy price for recharging the fleet. In this context and for the remainder of this section, we interpret $\EnergyPrice_{\texttt{r}}$ as a piecewise constant function instead of as a parameter.

\begin{figure}[!bt]
    \centering
    \fontsize{7pt}{7pt}\selectfont
    \resizebox{\textwidth}{!}{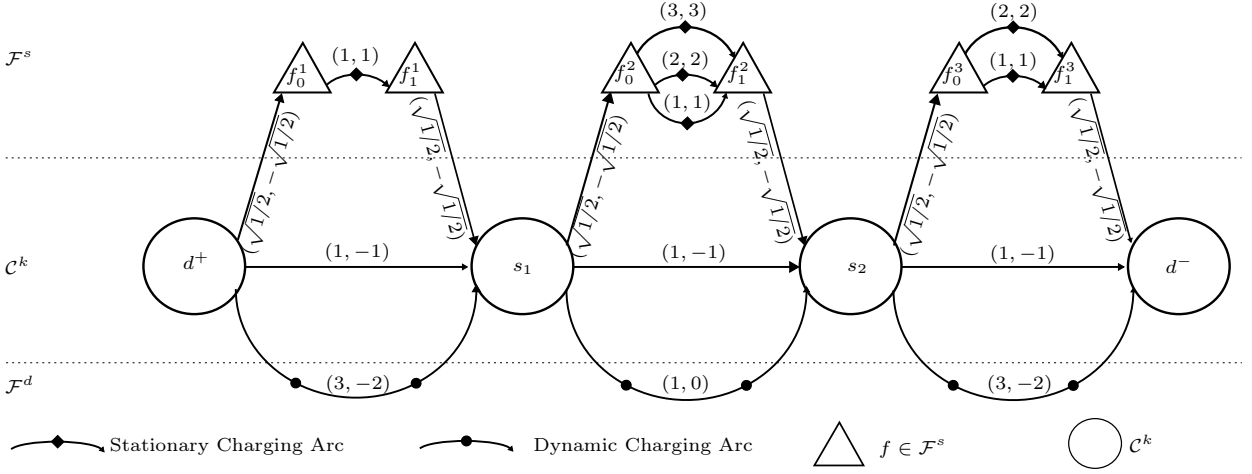}
    \caption{Illustrative example with depot and two stops on the corners of an equilateral triangle with side length one, one stationary charger in the center of the triangle, and a dynamic charging station on the side between $\Stop_1$ and $\Stop_2$. The vehicle travels with speed $1$, we consider a time discretization of $\TimeStep=1$, and all stations charge with $\StationChargingRate_{\Station}=1$. We assume that service time windows are such that the vehicle can charge a maximum of three time steps between servicing the two stops and similarly at other points of the route. Every arc is associated with a tuple reflecting time and net energy balance where a positive value encodes a situation in which charging replenishes more energy than the vehicle consumes.}
    \label{fig:vehicle-graph-example}
\end{figure}

\noindent \textbf{Graph expansion: } The subproblem decomposes by vehicle such that we construct one individual and independent vehicle graph $\Graph^{\Vehicle}=(\SetOfVertices^{\Vehicle}, \SetOfArcs^{\Vehicle})$ for each vehicle $\Vehicle \in \SetOfVehicles$. In the following, we explain the construction of a vehicle graph for an arbitrary vehicle $\Vehicle$ before we elaborate on the pruning steps that we apply to avoid unnecessary graph complexity. Figure~\ref{fig:vehicle-graph-example} provides an intuitive example of a vehicle graph.

Every vehicle route consists of $n+2$ mandatory locations that the vehicle has to service which results in $n+1$ inter-stop paths. Finding every such inter-stop path is a resource constrained shortest path problem on a mini-network in which the available resources depend on the solution of the previous mini-network. 

Let $\SetOfStationaryChargerReps := \{ \Station^{i}_{0}, \Station^{i}_{1}: \Station \in \SetOfStaticStations_{\mathbf{y}=1}, i = 1, \dots, n+1 \}$ be the set of stationary charging station representations such that $\SetOfStationaryChargerReps$ contains two representations per stationary candidate charging station and per inter-stop path of vehicle $\Vehicle$.

Finally, we construct the vertex set $\SetOfVertices^{\Vehicle} := \StopsOfLine \cup \SetOfStationaryChargerReps$ such that $\SetOfVertices^{\Vehicle}$ contains representations for all stops on the vehicle's route and the depot, as well as the stationary charging station representations. In this context, we pass the service time windows to a vertex $\Vertex \in \SetOfVertices^{\Vehicle}$ such that
\[
    (\TwBegin_{\Vertex}, \TwEnd_{\Vertex}) = 
    \begin{cases}
        (\TwBegin^{\Vehicle}_{\StopOfLine}, \TwEnd^{\Vehicle}_{\StopOfLine}) & \text{if } \Vertex = \StopOfLine \in \StopsOfLine \\
        (0, \infty) & \text{else.}
    \end{cases}
\]
Moreover, we construct the arc sets $\SetOfArcs^{\Vehicle}$ as the union of four subsets, which we detail in the following. An arc $\VehicleArc \in \SetOfArcs^{\Vehicle}$ is a quintuple defining their origin $i \in \SetOfVertices^{\Vehicle}$, their target $j \in \SetOfVertices^{\Vehicle}$, the time expenditure $\TravelTime$, the energy consumption $\Consumption$, and the recharged \gls{acr:soc} $\Recharge$ connected to the selection of the arc on the vehicle's route. 

\begin{description}
    \item[\medspace Direct arcs] represent direct connections between consecutive stops on a vehicle route. Formally, we define the set of direct arcs as $\{ (\StopOfLine, \StopOfLine+1, \TravelTime_{\StopOfLine, \StopOfLine+1}, \Consumption_{\StopOfLine, \StopOfLine+1}, 0): \StopOfLine \in \StopsOfLine \setminus \{\EndDepot\} \}$.
    \item[\medspace Stationary charging arcs] represent the stationary charging process. Let $\MaxChargingTimeSteps \in \mathbb{N}$ be any number of periods of duration $\TimeStep$ that the vehicle can spend charging at a stationary charging station $\Station \in \SetOfStaticStations$ while still maintaining feasibility in terms of \gls{acr:soc} limits and departure times. We determine a tight $\MaxChargingTimeSteps$ as part of the graph pruning discussed in more detail in the following. The set of stationary charging arcs is $\{(\Station^{i}_{0}, \Station^{i}_{1}, \TimeStep l, 0, \TimeStep l \StationChargingRate_{\Station}): \Station \in \SetOfStaticStations_{\mathbf{y}=1}, l=1,\dots,\MaxChargingTimeSteps, i=1, \dots, n+1 \}$.
    \item[\medspace Connecting arcs] connect stops with stationary charging stations. Formally, the union of the following sets yields the set of connecting arcs: the set $\{(\StopOfLine, \Station^{\StopOfLine+1}_{0}, \TravelTime_{\StopOfLine, \Station}, \Consumption_{\StopOfLine, \Station}, 0): \Station \in \SetOfStaticStations_{\mathbf{y}=1}, \StopOfLine=0,\dots,n \}$ allows the vehicle to visit a stationary charging station and vice versa $\{(\Station^{\StopOfLine}_{1}, \StopOfLine, \TravelTime_{\Station, \StopOfLine}, \Consumption_{\Station, \StopOfLine}): \Station \in \SetOfStaticStations_{\mathbf{y}=1}, \StopOfLine \in \StopsOfLine \setminus \{\EndDepot\} \}$ allows the vehicle to leave the station and service the next stop on its route.
    \item[\medspace Dynamic charging arcs] represent the option to detour between consecutive stops and charge dynamically at a station $\SetOfDynamicStations_{\mathbf{\TransformerVar} = 1}$. We add the following subset $\{ (\StopOfLine, \StopOfLine+1, \TravelTimeFunction_{\Station}(\StopOfLine, \StopOfLine+1), \ConsumptionFunction_{\Station}(\StopOfLine, \StopOfLine+1), \RechargeFunction_{\Station}): \Station \in \SetOfDynamicStations_{\mathbf{\TransformerVar} = 1}, \StopOfLine \in \StopsOfLine \setminus \{\EndDepot\} \}$ to $\SetOfArcs^{\Vehicle}$ where $\TravelTimeFunction_{\Station}$, and $\ConsumptionFunction_{\Station}$ are bivariate functions that map two vertices $\Vertex, \Vertex^{\prime} \in \SetOfVertices^{\Vehicle}$ to a time expenditure and energy consumption if the vehicle traverses the dynamic charging station $\Station$ between servicing the two given stops; and $\RechargeFunction_{\Station}$ yields the recharged energy when traversing the respective station:
    \begin{align*}
        \TravelTimeFunction_{\Station}(\Vertex, \Vertex^{\prime}) &= \TravelTime_{\Vertex, \SegmentStart_{1}} + \sum_{\Segment \in \SetOfSegments_{\Station}} \TravelTime_{\UnbracedSegment} + \TravelTime_{\SegmentEnd_{m}, \Vertex^{\prime}} \\
        \ConsumptionFunction_{\Station}(\Vertex, \Vertex^{\prime}) &= \Consumption_{\Vertex, \SegmentStart_{1}} + \sum_{\Segment \in \SetOfSegments_{\Station}} \Consumption_{\UnbracedSegment} + \Consumption_{\SegmentEnd_{m}, \Vertex^{\prime}} \\
        \RechargeFunction_{\Station} &= \sum_{\Segment \in \SetOfSegments_{\Station}} \StationChargingRate_{\Station} \TravelTime_{\UnbracedSegment} \mathbb{1}_{\{\DynamicChargerVar_{\UnbracedSegment} = 1\}}
    \end{align*}
\end{description}

    Note that $\SetOfStaticStations \cap \SetOfStops \neq \emptyset$ if candidate stationary charging stations are located at stops, which is why connecting arcs can be associated to zero time expenditure and \gls{acr:soc} balance. Moreover, we store some information on the contracted dynamic charging segments on the respective dynamic charging arcs: first, in order to be able to accurately evaluate the energy price curve $\EnergyPrice_{\texttt{r}}$ at the time when the charging process starts, we maintain variables $\TimeToCharge_{i,j}^{\TravelTime, \Consumption, \Recharge}, \medspace \VehicleArc \in \SetOfArcs^{\Vehicle}$ that yield the time a vehicle traverses the arc until the charging process starts. Thus, $\TimeToCharge_{i,j}^{\TravelTime, \Consumption, \Recharge}$ only takes non-zero values if $\VehicleArc$ is a dynamic charging arc. Second, we store the discrete energy balance along the contracted arcs such that we can discard paths that would violate the \gls{acr:soc} limits along the way. 

\noindent \textbf{Lower bounds: } We introduce lower bounds to the required \gls{acr:soc} resources for every vehicle $\Vehicle \in \SetOfVehicles$ at every vertex in the respective $\Graph^{\Vehicle}$. Let $\Vertex \in \SetOfVertices^{\Vehicle}$ with $\VertexEnergyBound$ being its local energy bound. We compute this bound as follows: we start by assigning a zero bound to the vertex that represents the depot where the vehicle route terminates~(l.~\ref{l:set-bound-end-depot}). We then assign lower bounds to all vertices representing stops or the remaining representation of the depot by backpropagating the local bounds through the graph in reversed order of the stop sequence~(l.~\ref{l:set-bound-stops}). The computed bounds rely on the fact that the stop sequence is fixed and determine the consumption if the vehicle does not detour at all but only services its stop sequence. Finally, we assign lower bounds to the vertices representing stationary chargers by leveraging the already computed bounds at stop representations and correcting the associated bounds by the required consumption to reach the respective stationary station~(l.~\mbox{\ref{l:set-bound-station-0}-\ref{l:set-bound-station-1}}). We save the computed bounds and share them between \gls{acr:ils} iterations by only incrementally updating them for changes in the respective charging station configuration.

\begin{algorithm}[!b]
\caption{\textsc{LocalEnergyBounds($\Graph^{\Vehicle})$}}
\label{alg:local-bounds}
\begin{algorithmic}[1]
    \State $\underline{d}^{-} \leftarrow 0$ \label{l:set-bound-end-depot} \Comment{Assign zero bound to sink node}
    \For{$\StopOfLine$ \text{ in } $\textsc{Reversed}(\StopsOfLine \setminus \{\EndDepot\})$}
        \State Let $\Vertex := \StopOfLine + 1$ and set $\underline{\StopOfLine} \leftarrow \Consumption_{\StopOfLine, \Vertex} + \VertexEnergyBound$ \label{l:set-bound-stops} \Comment{Direct servicing of stop sequence}
    \EndFor
    \For{$\StopOfLine$ \text{ in } $\StopsOfLine \setminus \{\StartDepot\}$}
        \For {$\Station \in \SetOfStaticStations$}
            \State $\underline{\Station}^{\StopOfLine}_{0} \leftarrow - \sum\limits_{\Vertex \in \IncomingArcs(\StopOfLine)} \Consumption_{\Vertex^{\prime}, \Vertex} + \VertexEnergyBound$ \label{l:set-bound-station-0} \Comment{Leveraging the already computed bounds}
            \State $\underline{\Station}^{\StopOfLine}_{1} \leftarrow \sum\limits_{\Vertex \in \OutgoingArcs(\StopOfLine)} \Consumption_{\Vertex, \Vertex^{\prime}} + \VertexEnergyBound$ \label{l:set-bound-station-1}
        \EndFor
    \EndFor
    \State \Return $\texttt{LocalEnergyBounds}$ \label{l:return-local-bounds}
\end{algorithmic}
\end{algorithm}

The local energy bounds additionally allow us to bound the routing cost for a given vehicle using the variable energy prices: 
\begin{equation} \label{eq:lb-routing-cost}
    \VehicleBound := \underline{d}^{+} \, \EnergyPrice_{\texttt{c}} + \max \{\underline{d}^{+} - \InitialSoC + \MinSoC, 0 \} \, \min_{\tilde{\ArrivalTimeVar}} \EnergyPrice_{\texttt{r}}(\tilde{\ArrivalTimeVar}), \medspace \Vehicle \in \SetOfVehicles.
\end{equation}

Equation~\eqref{eq:lb-routing-cost} multiplies the local energy consumption bound in the start representation of the depot with the cost for consumption and adds a positive term representing the minimum energy that has to be recharged at the associated cost.

\noindent \textbf{Graph pruning: } We apply graph pruning in order to avoid unnecessary complexity of the vehicle graphs $\Graph^{\Vehicle}, \Vehicle \in \SetOfVehicles$. Let $\MaxChargingTimeSteps$ denote the maximum number of timesteps that a vehicle can charge at a stationary charging station $\Station \in \SetOfStaticStations$ between servicing stops $(\StopOfLine, \StopOfLine+1), \StopOfLine \in \StopsOfLine \setminus \{\EndDepot\}$. We define $\MaxChargingTimeSteps$ as follows
\begin{equation*}
\MaxChargingTimeSteps := \lfloor \frac{1}{\TimeStep} \min \{ \TwEnd_{\StopOfLine+1} - \TwBegin_{\StopOfLine} - \TravelTime_{\StopOfLine, \Station} - \TravelTime_{\Station, \StopOfLine+1}, \frac{1}{\StationChargingRate_{\Station}} (\MaxSoC - \MinSoC + \Consumption_{\StopOfLine, \Station})\} \rfloor    
\end{equation*}

such that the upper limit on the charging time reflects both spatial and temporal constraints. In this context, we waive the construction of stationary charging arcs belonging to $\Station$ and their respective connecting arcs if 

\begin{itemize}
    \item[i.] The stationary charging station is located at a stop, i.e., $\Station = \StopOfLine \in \StopsOfLine \cap \SetOfStaticStations$ and we are constructing the mini-network that represents the inter-stop path between $\StopOfLine$ and $\StopOfLine+1$. Thus, we restrict the charging operations in this case to happen before the stop is serviced and therefore ensure that the vehicle is adhering to its departure time at the respective stop.
    \item[ii.] The detour to the respective charging station yields a negative net energy balance. Formally, if $\Consumption_{\StopOfLine, \Station} + \Consumption_{\Station, \StopOfLine+1} -\Consumption_{\StopOfLine, \StopOfLine+1} \geq \MaxChargingTimeSteps \TimeStep \StationChargingRate_{\Station}$ holds.
\end{itemize}

After the inclusion of the pruning steps outlined above, we further tighten $\MaxChargingTimeSteps$ by incorporating optimality aspects. Specifically, we tighten 
\[
\MaxChargingTimeSteps \leftarrow \min \{ \MaxChargingTimeSteps, \frac{1}{\TimeStep} \lceil \underline{\StopOfLine} + \Consumption_{\Station, \StopOfLine+1} - \Consumption_{\StopOfLine, \StopOfLine+1} \rceil \}
\]
such that the number of charging arcs at $\Station \in \SetOfStaticStations$ is limited by the local energy bounds at a given vertex $\StopOfLine \in \SetOfVertices^{\Vehicle}$.

With respect to the dynamic charging arcs representing the charging operations at $\Station \in \SetOfDynamicStations$, we apply the following pruning step. We waive the addition of arcs between two consecutive locations $(\StopOfLine, \StopOfLine+1), \StopOfLine \in \StopsOfLine \setminus \{\EndDepot\}$ if $\RechargeFunction_{\Station} \leq \ConsumptionFunction_{\Station}(\StopOfLine, \StopOfLine+1) - \Consumption_{\StopOfLine, \StopOfLine+1}$. Thus, we do not construct a dynamic charging arc if the recharged battery \gls{acr:soc} is smaller than the additional consumption from the detour to the dynamic station.

\noindent \textbf{Bidirectional A*-Algorithm: } We use a customized bidirectional ${A}^{\star}$-Algorithm with a similar structure as \cite{ThomasEtAl2019}. However, we note that our problem is different from the conventional \gls{acr:rcspp} because resources are constrained at every node of the graph.

\begin{algorithm}[!t]
    \caption{Bidirectional ${A}^{\star}$ with lazy dominance checks}
    \label{alg:a-star}
    \footnotesize
    \begin{algorithmic}[1]
        \State $\ExploredFwLabels_{\Vertex} \gets \langle \rangle, \, \ExploredBwLabels_{\Vertex} \gets \langle \rangle, \medspace \forall \Vertex \in \SetOfVertices^{\Vehicle}$ \label{l:init-explored-sets} \Comment{Initialize empty sorted sets of explored labels at all vertices}
        \State $\FwQueue \gets \{ (0, \langle 0, \InitialSoC, \TwBegin_{\StartDepot} \rangle, \StartDepot) \}, \,\BwQueue \gets \{ (0, \langle 0, \MinSoC, \TwEnd_{\EndDepot} \rangle, \EndDepot) \}$ \label{l:init-queues}
        \State $\Queue \gets \FwQueue \cup \BwQueue$ \label{l:init-queue} \Comment{Priority queue sorted by cost approximation $\LabelCost + \HeuristicComponent(\Label)$}
        \State $\LabellingIterations \gets 0$ \label{l:init-iter-count}
        \While{$\FwQueue \neq \emptyset$} \label{l:stop-criterium}
            \State $\LabellingIterations \gets \LabellingIterations + 1$ \label{l:increase-iter-number}
            \State $\_, \Label, \Vertex \gets \Queue.pop()$  \label{l:pop-from-queue}
            \If{$\texttt{IsJointLabel}(\Label) \lor (\Vertex = \EndDepot \land \texttt{IsForwardLabel}(\Label))$} 
                \State \Return $\texttt{ConstructPath}(\Label)$ \label{l:termination}
            \EndIf
            \If{$\texttt{IsForwardLabel}(\Label)$} \label{l:case-fw}
                \If{$\texttt{IsDominated}(\Label, \ExploredFwLabels)$} \Comment{cf. Defintion~\eqref{def:fw-dominance}} 
                    \State \textbf{continue} \label{l:dominance-check-fw}
                \EndIf
                \For{$\Arc \in \OutgoingArcs(\Vertex)$} \label{l:loop-outgoing-arcs}
                    \State $\Label^{\prime} \gets \texttt{PropagateLabel}(\Label, \Arc)$ \Comment{cf. Equations~\eqref{eq:fw-prop-cost}-\eqref{eq:fw-prop-dep-time}} \label{l:propagate-fw}
                    \If{$\texttt{IsValid}(\Label^{\prime})$}  \Comment{cf. Equation~\eqref{eq:fw-prop-feasible}} 
                        \State $\FwQueue.push((\LabelCost_{\Label^{\prime}} + \HeuristicComponent(\Label^{\prime}), \langle \LabelCost_{\Label^{\prime}}, \StateOfChargeVar_{\Label^{\prime}}, \ArrivalTimeVar_{\Label^{\prime}} \rangle, j))$ \label{l:add-fw-label-to-queue}
                    \EndIf
                \EndFor
                \If{$\texttt{WasSuccesfullyPropagated}(\Label)$}
                    \State $\ExploredFwLabels \gets \texttt{AddAndMaintainDominanceBounds}(\Label)$  \label{l:add-fw-label-to-explored}
                \EndIf
                \For{$\BackwardLabel \in \ExploredBwLabels_{\Vertex}$} \label{l:check-fw-merge}
                    \If{$\texttt{CanFormJointLabel}(\Label, \BackwardLabel)$}
                        \State $\Label^{\prime \prime} \gets \texttt{ConstructJointLabel}$ \Comment{cf. Equation~\eqref{eq:joint-label}}
                        \State $\Queue.push((\LabelCost_{\Label} + \LabelCost_{\BackwardLabel}, \Label^{\prime \prime}, \EndDepot))$ \label{l:add-joint-label-to-queue-1}
                        \State \textbf{break}
                    \EndIf
                \EndFor
            \Else
                \If{$\texttt{IsDominated}(\Label, \ExploredBwLabels)$} \Comment{cf. Defintion~\eqref{def:bw-dominance}}
                    \State \textbf{continue} \label{l:dominance-check-bw}
                \EndIf
                \For{$\Arc \in \IncomingArcs(\Vertex)$} \label{l:loop-incoming-arcs}
                    \State $\Label^{\prime} \gets \texttt{PropagateLabel}(\Label, \Arc)$ \label{l:propagate-bw} \Comment{cf. Equations~\eqref{eq:bw-prop-cost}-\eqref{eq:bw-prop-time}}
                    \If{$\texttt{IsValid}(\Label^{\prime})$}  \Comment{cf. Equation~\eqref{eq:bw-prop-feasible}}
                        \State $\BwQueue.push((\LabelCost_{\Label^{\prime}} + \HeuristicComponent(\Label^{\prime}), \langle \LabelCost_{\Label^{\prime}}, \StateOfChargeVar_{\Label^{\prime}}, \ArrivalTimeVar_{\Label^{\prime}} \rangle, i))$ \label{l:add-bw-label-to-queue}
                    \EndIf
                \EndFor
                \If{$\texttt{WasSuccesfullyPropagated}(\Label) \lor \Vertex \in \{\StartDepot, \EndDepot\}$}
                        \State $\ExploredBwLabels \gets \texttt{AddAndMaintainDominanceBounds}(\Label)$ \label{l:add-bw-label-to-explored}
                \EndIf
                \For{$\ForwardLabel \in \ExploredFwLabels_{\Vertex}$} \label{l:check-bw-merge}
                    \If{$\texttt{CanFormJointLabel}(\Label, \ForwardLabel)$}
                        \State $\Label^{\prime \prime} \gets \texttt{ConstructJointLabel}$ \Comment{cf. Equation~\eqref{eq:joint-label}}
                        \State $\Queue.push((\LabelCost_{\Label} + \LabelCost_{\ForwardLabel}, \Label^{\prime \prime}, \EndDepot))$ \label{l:add-joint-label-to-queue-2}
                        \State \textbf{break}
                    \EndIf
                \EndFor
            \EndIf
        \EndWhile
    \State \Return \textbf{infeasible}
    \end{algorithmic}
\end{algorithm}

We detail the bidirectional A*-Algorithm in Algorithm~\ref{alg:a-star}. The algorithm starts by initializing sets of explored forward and backward labels at every node in the vehicle graph~(l.~\ref{l:init-explored-sets}). Moreover, we initialize the priority queue $\Queue$ with a starting forward label at $\StartDepot$, and a starting backward label at $\EndDepot$~(l.~\ref{l:init-queues}-\ref{l:init-queue}), and initialize an iteration count~(l.~\ref{l:init-iter-count}). The algorithm runs until the forward priority queue is empty~(l.~\ref{l:stop-criterium}) or the algorithm found an optimal solution. In every iteration, the algorithm increases the iteration counter by one~(l.~\ref{l:increase-iter-number}) and pops the element from the combined priority queue that yields the label $\Label$ with the lowest term $\LabelCost + \HeuristicComponent(\Label)$~(l.~\ref{l:pop-from-queue}). In case the popped label is a joint label or is associated to $\EndDepot$, the algorithm terminates and returns the respective path~(l.~\ref{l:termination}). Otherwise the algorithm checks if the popped label is a forward or backward label and handles it accordingly. In both cases, the algorithm checks the popped label for dominance~(l.~\ref{l:dominance-check-fw},~\ref{l:dominance-check-bw}) and discards dominated labels before returning to the next iteration. For non-dominated forward labels at a given vertex $\Vertex \in \SetOfVertices^{\Vehicle}$, the algorithm extends the popped label along all outgoing arcs from $\Vertex$~(l.~\ref{l:loop-outgoing-arcs}) and adds the resulting labels that are valid according to~\eqref{eq:fw-prop-feasible} with their associated information to the priority queue~(l.~\ref{l:add-fw-label-to-queue}). Vice versa, if the popped label is a backward label, the algorithm extends it along all incoming arcs to $\Vertex$~(l.~\ref{l:loop-incoming-arcs}), checks for validity according to~\eqref{eq:bw-prop-feasible}, and adds the labels to the respective priority queue~(l.~\ref{l:add-bw-label-to-queue}). If the label represents a feasible path because it resulted in at least one valid extended label, the algorithm adds the label to the set of explored labels and maintains the stored dominance bounds accordingly~(l.~\ref{l:add-fw-label-to-explored},~\ref{l:add-bw-label-to-explored}). Finally, the algorithm determines if the popped label can form a joint label with any label of the opposite direction at the same vertex. For this purpose, it compares a forward label with all explored backward labels~(l.~\ref{l:check-fw-merge}) and vice versa~(l.\ref{l:check-bw-merge}) for the conditions in Equation~\eqref{eq:joint-label}. If two labels can form a joint label, the algorithm adds the joint label to the queue and determines the labels place in the priority queue by the sum of the label costs of the two merged labels~(l.~\ref{l:add-joint-label-to-queue-1}, \ref{l:add-joint-label-to-queue-2}).     

In this context, $\HeuristicComponent(\Label)$ denotes our heuristic distance approximation with which we prioritize promising labels to be extended in the label setting algorithm. More specifically, we compute $\HeuristicComponent(\Label)$ for a label $\Label$ at a vertex $\Vertex \in \SetOfVertices^{\Vehicle}$ as 
\begin{equation} \label{eq:dist-approx}
    \HeuristicComponent(\Label) = 
    \begin{cases}
        \max \{ \MinSoC + \VertexEnergyBound - \StateOfChargeVar_{\Label}, 0 \} \, \min\limits_{\tilde{\ArrivalTimeVar} \geq \ArrivalTimeVar_{\Label}} \EnergyPrice_{\texttt{r}}(\tilde{\ArrivalTimeVar}) + \VertexEnergyBound \, \EnergyPrice_{\texttt{c}}, & \medspace \text{if } \Label \in \SetOfForwardLabels_{\Vertex} \\ 
        \max \{ \StateOfChargeVar_{\Label} + \underline{d}^{+} - \VertexEnergyBound - \InitialSoC, 0 \} \, \min\limits_{\tilde{\ArrivalTimeVar} \leq \ArrivalTimeVar_{\Label}} \EnergyPrice_{\texttt{r}}(\tilde{\ArrivalTimeVar}) + (\underline{d}^{+} -\VertexEnergyBound) \, \EnergyPrice_{\texttt{c}}, & \medspace \text{else}.
    \end{cases}
\end{equation}

When computing the distance approximation with variable prices in Equation~\ref{eq:dist-approx}, we consider the minimum energy price $\EnergyPrice$ that can be achieved after the labels departure time on the price curve, respectively before in the case of a backward label. We apply the resulting price to the term that represents the minimal required recharging amount and the constant cost for consumption to the term that reflects the minimum additional energy requirement at the given label $\Label$. 

\noindent \textbf{Labels and label propagation: } Let $\SetOfForwardLabels_{\Vertex}, \SetOfBackwardLabels_{\Vertex}$ denote all forward labels and backward labels, that are associated with vertex $\Vertex \in \SetOfVertices^{\Vehicle}$. A forward label $\ForwardLabel \in \SetOfForwardLabels$ is a triplet $\ForwardLabel := (\LabelCost_{\ForwardLabel}, \StateOfChargeVar_{\ForwardLabel}, \ArrivalTimeVar_{\ForwardLabel})$ comprising the cost associated with the label, the battery \gls{acr:soc}, and the departure time at the associated vertex of the label. 

The set of backward labels $\SetOfBackwardLabels$ serves as a complement of information that enables our algorithm to yield solutions for real-world instance sizes quickly. For this purpose a backward label $\BackwardLabel := (\LabelCost_{\BackwardLabel}, \StateOfChargeVar_{\BackwardLabel}, \ArrivalTimeVar_{\BackwardLabel}, \Headway_{\BackwardLabel})$ is a quadruple that contains, analogous to a forward label, a cost, a battery \gls{acr:soc}, and a departure time. Backward labels carry an additional value $\Headway$ that propagates the maximum \gls{acr:soc} offset that maintains feasibility when the label forms a joint label with another forward label. 

We propagate forward labels $\ForwardLabel \in \SetOfForwardLabels$ at $i \in \SetOfVertices^{\Vehicle}$ along arcs $(i,j,\TravelTime ,\Consumption, \Recharge) \in \SetOfArcs^{\Vehicle}$ to a new label $\ForwardLabel^{\prime}$ at $j \in \SetOfVertices^{\Vehicle}$ by setting
\begin{align}
    \LabelCost_{\ForwardLabel^{\prime}} &:= \LabelCost_{\ForwardLabel} + \EnergyPrice_{\texttt{r}}(\ArrivalTimeVar_{\ForwardLabel} + \TimeToCharge_{i,j}^{\TravelTime, \Consumption, \Recharge}) \, \Recharge + \EnergyPrice_{\texttt{c}} \, \Consumption \label{eq:fw-prop-cost} \\ 
    \StateOfChargeVar_{\ForwardLabel^{\prime}} &:= \StateOfChargeVar_{\ForwardLabel} - \Consumption + \Recharge \label{eq:fw-prop-soc}\\ 
    \ArrivalTimeVar_{\ForwardLabel^{\prime}} &:= \max\{\ArrivalTimeVar_{\ForwardLabel} + \TravelTime, \TwBegin_{j}\}. \label{eq:fw-prop-dep-time}
\end{align}

The propagated forward label is valid if and only if
\begin{equation} \label{eq:fw-prop-feasible}
    \MinSoC \leq \StateOfChargeVar_{\ForwardLabel^{\prime}} \leq \MaxSoC \medspace \land \medspace \ArrivalTimeVar_{\ForwardLabel^{\prime}} \leq \TwEnd_{j}.
\end{equation}

Three comments are in order. First, we implicitly allow vehicles to idle on arcs without allowing idle time to change the \gls{acr:soc} profile through increased consumption or extended recharging. Second, we do not allow partial recharging in this approach by propagating full charging amounts in Equation~\eqref{eq:fw-prop-soc} and rendering labels infeasible that violate the maximum battery \gls{acr:soc} in~\eqref{eq:fw-prop-feasible}, because initial tests on our field project's test track showed that charging rates are not high enough and dynamic segments are not long enough to make partial re-charging a beneficial edge-case in real-world settings. Third, we note that by emulating the described forward label extension within a standard dynamic programming approach, the resulting algorithm would converge towards an optimal solution. However, in order to solve larger instances with many charging stations and substantial symmetry due to uniform charging rates and charging stations being located at stops, we extended our label setting algorithm to a bidirectional approach.  

In this context, we propagate backward labels $\BackwardLabel \in \SetOfBackwardLabels$ at $j \in \SetOfVertices^{\Vehicle}$ along arcs $(i,j,\TravelTime ,\Consumption,\Recharge) \in \SetOfArcs^{\Vehicle}$ to a new label $\BackwardLabel^{\prime}$ at $i$ by setting
\begin{align}
    \LabelCost_{\BackwardLabel^{\prime}} &:= \LabelCost_{\BackwardLabel} + \EnergyPrice_{\texttt{r}}(\min\{\ArrivalTimeVar_{\BackwardLabel} - \TravelTime, \TwEnd_{i}\} + \TimeToCharge_{i,j}^{\TravelTime, \Consumption, \Recharge}) \, \Recharge + \EnergyPrice_{\texttt{c}} \, \Consumption \label{eq:bw-prop-cost} \\
    \StateOfChargeVar_{\BackwardLabel^{\prime}} &:= \StateOfChargeVar_{\BackwardLabel} + \Consumption - \Recharge \label{eq:bw-prop-soc} \\
    \ArrivalTimeVar_{\BackwardLabel^{\prime}} &:= \min\{\ArrivalTimeVar_{\BackwardLabel} - \TravelTime, \TwEnd_{i}\} \label{eq:bw-prop-time} \\
    \Headway_{\BackwardLabel^{\prime}} &:= \min\{\Headway_{\BackwardLabel},  \MaxSoC - \StateOfChargeVar_{\BackwardLabel^{\prime}} \} \label{eq:bw-prop-headway} 
\end{align}

In Equation~\eqref{eq:bw-prop-cost}, we propagate the label cost and variable energy prices must be evaluated at the new labels departure time in order to be consistent with the forward label propagation. In Equations~\eqref{eq:bw-prop-soc} and~\eqref{eq:bw-prop-time}, we propagate the battery \gls{acr:soc} and the departure time backwards. Additionally, we propagate the headspace of the label's \gls{acr:soc} profile in Equation~\eqref{eq:bw-prop-headway}. In line with our logic when propagating forward labels, we allow the vehicle to idle on arcs. A backward label is valid if and only if
\begin{equation} \label{eq:bw-prop-feasible}
     \MinSoC \leq \StateOfChargeVar_{\BackwardLabel^{\prime}} \leq \MaxSoC \medspace \land \medspace \ArrivalTimeVar_{\BackwardLabel^{\prime}} \geq \TwBegin_{i}.
\end{equation}

A forward label $\ForwardLabel \in \SetOfForwardLabels_{\Vertex}$ and a backward label $\BackwardLabel \in \SetOfBackwardLabels_{\Vertex}$ can form a joint label if
\begin{equation} \label{eq:joint-label}
    0 \leq \StateOfChargeVar_{\ForwardLabel} - \StateOfChargeVar_{\BackwardLabel} \leq \Headway_{\BackwardLabel} \medspace \land \medspace \ArrivalTimeVar_{\ForwardLabel} \leq \ArrivalTimeVar_{\BackwardLabel}.
\end{equation}

Note that forming a joint label and extracting its solution $\Path^{\Vehicle}$ means to adapt the \gls{acr:soc} profile at all vertices in $\Path^{\Vehicle}$ that were explored by a backward label exactly by $\StateOfChargeVar_{\ForwardLabel} - \StateOfChargeVar_{\BackwardLabel}$. To ensure that such an adjustment yields a feasible solution $\Path^{\Vehicle}$, we restrict the difference between $\StateOfChargeVar_{\ForwardLabel}$ and $\StateOfChargeVar_{\BackwardLabel}$ to remain below the backwards propagated headspace $\Headway_{\BackwardLabel}$.

\noindent \textbf{Local dominance bounds: } To accelerate convergence, we apply dominance criteria to both forward and backward labels. We check labels for dominance in a lazy fashion after popping them from the priority queue instead of performing the check when adding them to the queue. This is particularly beneficial for our heuristic approach in which a substantial part of the queue remains unexplored and therefore never invokes a dominance check. 

\begin{definition}[Forward Domination] \label{def:fw-dominance}
    Let $\ForwardLabel_a, \ForwardLabel_b \in \SetOfForwardLabels_{\Vertex}$ be two forward labels associated with the same vertex $\Vertex \in \SetOfVertices^{\Vehicle}$. Then, $\ForwardLabel_a$ dominates $\ForwardLabel_b$ ($\ForwardLabel_a \succcurlyeq \ForwardLabel_b$) if and only if
    \[
        \LabelCost_{\ForwardLabel_a} \leq \LabelCost_{\ForwardLabel_b} \medspace \land \medspace \StateOfChargeVar_{\ForwardLabel_a} \geq \StateOfChargeVar_{\ForwardLabel_b} \medspace \land \medspace \ArrivalTimeVar_{\ForwardLabel_a} \leq \ArrivalTimeVar_{\ForwardLabel_b}. 
    \]
\end{definition}

Definition~\ref{def:fw-dominance} ensures that any forward label $\ForwardLabel_a$ only dominates another forward label $\ForwardLabel_b$ if it achieves at least an equal battery \gls{acr:soc} in at maximum the same time and at at maximum the same cost.

Moreover, we define dominance criteria for backward labels.  

\begin{definition}[Backward Domination] \label{def:bw-dominance}
    Let $\BackwardLabel_a, \BackwardLabel_b \in \SetOfBackwardLabels_{\Vertex}$ be two backward labels associated with the same vertex $\Vertex \in \SetOfVertices^{\Vehicle}$. Then, $\BackwardLabel_a$ dominates $\BackwardLabel_b$ ($\BackwardLabel_a \succcurlyeq \BackwardLabel_b$), if and only if
    \[
        \LabelCost_{\BackwardLabel_a} \leq \LabelCost_{\BackwardLabel_b} \medspace \land \medspace \StateOfChargeVar_{\BackwardLabel_a} \leq \StateOfChargeVar_{\BackwardLabel_b} \medspace \land \medspace \ArrivalTimeVar_{\BackwardLabel_a} \geq \ArrivalTimeVar_{\BackwardLabel_b}.
    \]
\end{definition}

Definition~\ref{def:bw-dominance} reflects the same intuition as Definition~\ref{def:fw-dominance}. Specifically, it states that a backward label $\BackwardLabel_a$ dominates another backward label $\BackwardLabel_b$ if it yields at maximum the same cost for at maximum the same \gls{acr:soc} at at minimum the same time. 

In order to check dominance efficiently, we maintain both local forward dominance bounds as well as local backward dominance bounds at all vertices in $\SetOfVertices^{\Vehicle}$. Specifically, we sort the explored forward labels associated to any given vertex in ascending order by their cost. Let $\langle \ForwardLabel_1, \ForwardLabel_2, \dots, \ForwardLabel_{m} \rangle$ be such a sorted set of forward labels at an arbitrary vertex $\Vertex \in \SetOfVertices^{\Vehicle}$ and let $m=|\SetOfForwardLabels_{\Vertex}|$ be the current number of forward labels at this vertex. We maintain corresponding \gls{acr:soc} bounds $\langle \StateOfChargeVar_{\ForwardLabel_1}, \max\limits_{i=1,2} \StateOfChargeVar_{\ForwardLabel_i}, \dots, \max\limits_{i=1,\dots,m} \StateOfChargeVar_{\ForwardLabel_i} \rangle$ and departure time bounds $\langle \ArrivalTimeVar_{\ForwardLabel_1}, \min\limits_{i=1,2} \ArrivalTimeVar_{\ForwardLabel_i}, \dots, \min\limits_{i=1,\dots,m} \ArrivalTimeVar_{\ForwardLabel_i} \rangle$. We use these local dominance bounds by rejecting dominance for an arbitrary forward label $\ForwardLabel_a$ with $j \in \{2, \dots, m\}$ such that $\LabelCost_{\ForwardLabel_i} \leq \LabelCost_{\ForwardLabel_a}, \medspace \forall i \in \{1, \dots, j-1\}$ if $\max\limits_{i=1,\dots,j} \StateOfChargeVar_{\ForwardLabel_i} \leq \StateOfChargeVar_{\ForwardLabel_a}$ or $\min\limits_{i=1,\dots,j} \ArrivalTimeVar_{\ForwardLabel_i} \geq \ArrivalTimeVar_{\ForwardLabel_a}$ hold. If dominance cannot be rejected based on the local dominance bounds, we check sequentially if labels $\ForwardLabel_j, \dots, \ForwardLabel_1$ dominate $\ForwardLabel_a$. The case in which Label $\ForwardLabel_a$ yields the lowest label cost, i.e., $\LabelCost_{\ForwardLabel_a} < \LabelCost_{\ForwardLabel_1}$, is trivial and $\ForwardLabel_a$ is undominated. Moreover, we maintain similar bounds for backward labels in which the \gls{acr:soc} bounds decrease and the departure time bounds increase from cheaper to more expensive labels.

\noindent \textbf{Heuristic solution: } We interrupt Algorithm~\ref{alg:a-star} if two conditions hold and change the algorithm's focus on finding a heuristic solution. First, we require $\LabellingIterations_{\texttt{min}}$ iterations to be conducted, formally $\LabellingIterations > \LabellingIterations_{\texttt{min}}$. Second, we require that the forward search and the backward search overlap at minimum at one vertex $\Vertex \in \StopsOfLine$, formally $\exists \Vertex \in \StopsOfLine: \ExploredFwLabels_{\Vertex} \neq \langle \rangle \land \ExploredBwLabels_{\Vertex} \neq \langle \rangle$. 

If the two conditions hold, we accept the joint label with the least cost that is already in the queue. If no joint label exists, we re-sort the priority queue based on the label's \gls{acr:soc} deviation from the minimum \gls{acr:soc} required to finish the route $\HeuristicSortingKey_{\Label}$:
\begin{equation} \label{eq:heuristic-sorting}
\HeuristicSortingKey_{\Label} =
\begin{cases}
    |\VertexEnergyBound + \MinSoC - \StateOfChargeVar_{\Label}|, & \medspace \text{if } \Label \in \SetOfForwardLabels_{\Vertex} \\ 
    |\InitialSoC - \underline{d}^{+} + \VertexEnergyBound - \StateOfChargeVar_{\Label}|, & \medspace \text{else}.
\end{cases}
\end{equation}

Additionally, we focus either on forward or on backward labels: if $|\FwQueue| \leq |\BwQueue|$, we only extend forward labels in the order of their associated $\HeuristicSortingKey_{\ForwardLabel}, \medspace \ForwardLabel \in \SetOfForwardLabels$, and otherwise we only extend backward labels. Then, we continue Algorithm~\ref{alg:a-star} with these adaptions until we find the first joint label or either of the two priority queues $\FwQueue$ and $\BwQueue$ is empty. In this context, the sorting by deviation from required \gls{acr:soc} fosters the extension of labels that can form joint labels with already existing labels of the opposite direction. Moreover, the focus on either forward or backward labels accelerates the termination of the algorithm because it focuses on proving infeasibility. However, note that terminating the algorithm when $\BwQueue = \emptyset$, introduces an additional heuristic element if the underlying assumption $\ConsumptionFunction_{\Station}(\StopOfLine, \StopOfLine+1) \geq \RechargeFunction_{\Station}(\StopOfLine, \StopOfLine+1), \medspace \Station \in \SetOfDynamicStations, \StopOfLine \in \StopsOfLine \setminus \{\EndDepot\}$ is violated and optimal solutions require to reach the depot with $\StateOfChargeVar^{\Vehicle}_{\EndDepot} > \MinSoC$.

\section{Experimental Design} \label{sec:exp-design}
We implemented our algorithmic framework with Python 3.10.2 and outsourced Algorithm~\ref{alg:a-star} to a C++ implementation which we integrated via Pybind11. Likewise, we implemented Problem~\ref{mips:deterministic} as a benchmark in IBM ILOG CPLEX Optimization Studio 22.1 using the Python API.

We run all our experiments with the same hyper-parameters. More specifically, we set the perturbation strength $\PerturbationStrength = 2$, and $\ILSIterations_{\texttt{max}} = 10$. The sample-based neighborhoods of our local search operators sample $\NeighborhoodSize = 2$ random elements from the configurations to intensify a given solution, and we use a constant $\epsilon = 1\mathrm{e}-3$. Finally, we run Algorithm~\ref{alg:a-star} for a minimum of $\LabellingIterations_{\texttt{min}} = 2.5\mathrm{e}5$ iterations before we consider falling back to searching for a heuristic solution (cf. Section~\ref{sec:subproblem}).

To benchmark our algorithm and derive managerial insights, we constructed two different instance sets: artificial instances to investigate our algorithms computational performance and real-world instances to evaluate the concept of dynamic inductive charging based on the city of Hof, Germany.

\noindent \textbf{Artificial instances:} We create artificial benchmark instances based on instances from the \gls{acr:evrptw-pr}, which consist of 1 depot, 100 stop locations and 20 stationary charging stations. Although originally motivated by a logistic application, these instances exhibit rather generic features. Moreover, we study a concept that can be extended to use cases beyond the investigated shuttles, e.g., logistic applications or airpron buses. The \gls{acr:evrptw-pr} instances contain three classes: Type C, where the stop locations and charging station locations are clustered; Type R with completely random stop locations; and Type RC which combine both clustered and random stop locations. Within each class, there are two subcategories: instances with small $\MaxSoC$ and low station charging rates, resulting in a larger number of shorter routes (1), as well as instances with large $\MaxSoC$ and high station charging rates resulting in fewer and longer routes (2). We construct the stop sequences on the basis of the routes in the best-known solutions from the literature for every instance, discard stopover times, and consider a time-discretization with $\TimeStep = 1$. 

To evaluate our algorithm's performance based on meaningful instances, we add dynamic candidate charging stations on up to $50\%$ of the arcs $(\StopOfLine, \StopOfLine+1), \medspace \StopOfLine=1, \dots, n-1$ for every vehicle $\Line \in \SetOfVehicles$. We select the longest arcs in these subsets and randomly break the dynamic charging stations up into two or three segments. Furthermore, we add stationary charging stations at up to $2\,|\SetOfVehicles|$ randomly selected intersections between different vehicle routes.

We uniformly generate costs $\InvestmentCost_{\Station} = [15, 25], \Station \in \SetOfStaticStations$, and $\InvestmentCost_{\Station} = [1.5, 2.5], \Station \in \SetOfDynamicStations$. Moreover, we set the variable cost for dynamic charging segments to $1\mathrm{e}-3$, assume $\EnergyPrice_{\texttt{c}}=5\mathrm{e}-2$ and $\EnergyPrice_{\texttt{r}}=5\mathrm{e}-2$, and increase the charging rates of randomly selected dynamic charging stations by a factor of four with a probability of $\frac{2}{3}$.

\begin{figure}[!t]
\centering
\begin{minipage}{0.49\textwidth}
    \centering
    \includegraphics[width=0.95\textwidth]{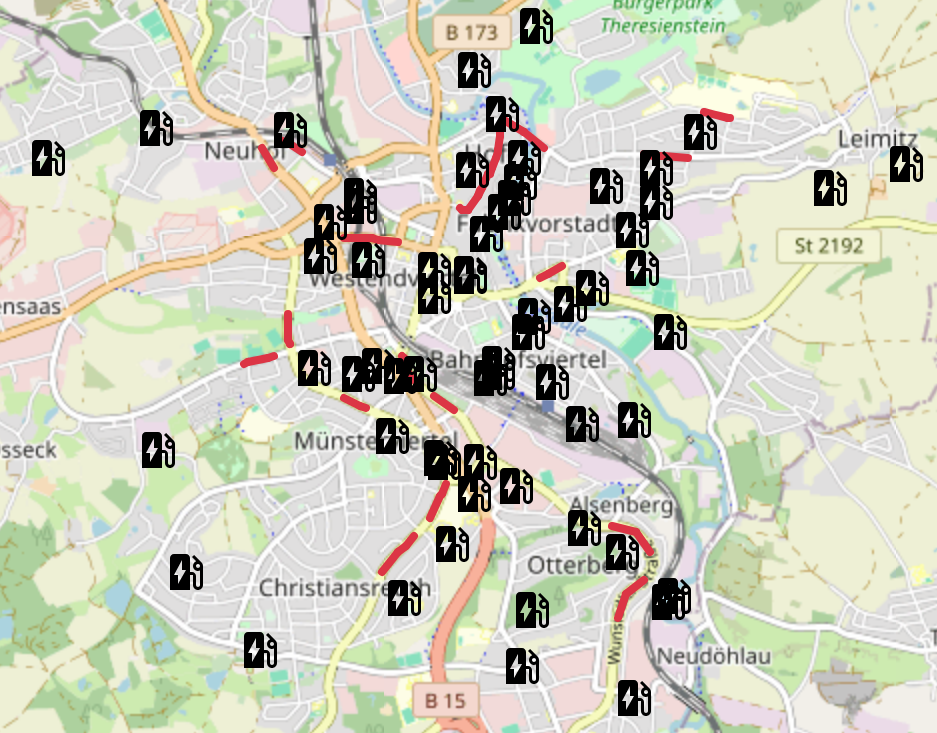}
\end{minipage}%
\begin{minipage}{0.49\textwidth}
    \centering
    \includegraphics[width=0.95\textwidth]{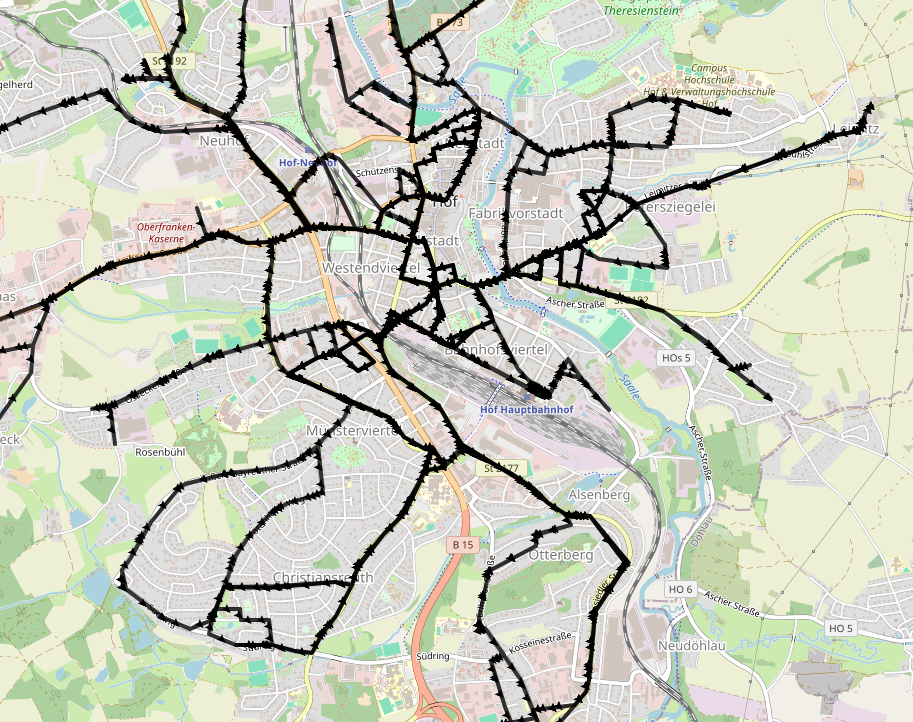} 
\end{minipage}
\caption{Left: Illustration of potential charging infrastructure in a selected scenario in Hof, Germany. Icons encode stationary charging stations and dynamic charging segments are highlighted in red. Right: Vehicle routes.}
\label{fig:hof_layout}
\end{figure}

\noindent \textbf{Case study:} We conduct a case study in the city Hof, Germany. Figure~\ref{fig:hof_layout} shows the spatial distribution of candidate charging stations in an example scenario, and visualizes the vehicle routes without detours on this network. We summarize the base parameters for our case study in Table~\ref{tbl:case-study} and elaborate on the details in the following. We observed different consumption levels per kilometer traveled on our test track. While the geospatial conditions such as the incline of the the roads certainly play a role, much of the variation comes from external influences such as the temperature and the passenger load. We therefore evaluate different consumption patterns in the range $[0.5 - 2.0]$ representing scenarios between cold winter days with high passenger demand and low-demand spring days that do not require to run the heating or air conditioning. Furthermore, to extend battery life, we enforce the \gls{acr:soc} limits to $10-90\%$ of the vehicles' battery capacity.

\begin{table}[!b]
	\begin{threeparttable}
		\caption{Parameter Case-Study}
		\label{tbl:case-study}
		\centering
		\begin{tabular*}{\textwidth}{@{\extracolsep{\fill}}lrrr@{}}
			\textbf{Parameter} & \textbf{Unit} & \textbf{Base case} & \textbf{Reasoning} \\
			\midrule
            Max. \gls{acr:soc} & kWh & 29.7 & $90\%$ of battery capacity \\
            Min. \gls{acr:soc} & kWh & 3.3 & $10\%$ of battery capacity \\
            Stationary fix cost & EUR & 5.04 & Test track \\
            Dynamic fix cost & EUR & 4.11 & Test track \\
            Variable cost & EUR\slash m & 0.19 & Test track \\
            Energy cost $\EnergyPrice_{\texttt{c}}$ & EUR\slash kWh & 0.3 & Assumption \\
            Charging rate & kW & 30 & Assumption \\
			\bottomrule
		\end{tabular*}
	\end{threeparttable}
\end{table}

Moreover, we conducted a total-cost analysis on the infrastructure cost in our test track. The total-cost analysis yields the cost for constructing stationary charging stations $\InvestmentCost_{\Station}, \medspace \Station \in \SetOfStaticStations$, the cost for dynamic charging stations $\InvestmentCost_{\Station}, \medspace \Station \in \SetOfDynamicStations$, and the cost per meter of dynamic charging segment. We disaggregated the values to reflect $8$ hours of operation assuming a linear depreciation over $25$ years, $365$ days of operation per year, and $2$ shifts per day. 

Note that in \gls{acr:milas}, shuttles are charging inductively at $7$kW. However, the hurdles for a real-life implementation of charging rates up to $30$kW are rather regulatory than technological which is why we assume the higher charging rate when evaluating the technology's future potential. This charging rate  assumption is still conservative compared to other studies \citep[cf.][]{ChenYinEtAl2018}. 

Finally, we extracted the vehicle routes of $12$ vehicles servicing the city bus lines in Hof, Germany --- which is a municipality that is part of the broader project region --- by sensible concatenation of trips published by the operator. We extended the vehicles' timetables such that their routes are feasible given the generally lower speed and battery capacities of shuttles compared to buses, and considered $8$ hours of operations. From the set of stops, we randomly sample $10$ different scenarios by determining random subsets of stops that allow the construction of stationary charging stations and serve as candidate stationary charging stations. We enrich this set of candidate charging stations by $11$ stand-alone stationary charging stations and $6.09$ km of dynamic charging segments belonging to $22$ dynamic charging stations. 

\section{Results} \label{sec:results}
We evaluate our algorithmic framework's computational performance in a computational study in Section~\ref{sec:comp-results} and derive managerial insights from the case study in the subsequent Section~\ref{sec:managerial-insights}.

\subsection{Computational Study} \label{sec:comp-results}
We set a runtime limit of 60 minutes and warmstart the solver at pre-computed solutions that enforce every vehicle to operate on their shortest route with respect to distance. If no such solution exists, we do not pass a warmstart solution to the solver. We conducted our experiments on an Intel i9-9900 processor using a single logical core with 16 GB of memory. 

In order to demonstrate the ability of our algorithmic framework to solve various instances, we compare its performance on the adapted \gls{acr:evrptw-pr} instances with the \gls{acr:bab} approach of the commercial solver CPLEX. To enable a fair comparison, we initialize the \gls{acr:ils} with the same configuration as the benchmark where possible. Table~\ref{tbl:comp-benchmark} yields descriptive statistics about the considered instance classes, the number of instances within each class, and initial insights into the outcome of our experiments. 

\begin{table}[!b] 
\centering
\begin{tabular}{c c c c c c c }
\textbf{} & \textbf{C1} & \textbf{C2} & \textbf{R1} & \textbf{R2} & \textbf{RC1} & \textbf{RC2} \\
\hline
\textbf{\# instances} & 9  & 8  & 12  & 11 & 8 & 8 \\
\textbf{\# warmstart available} & 4  & 8  & 12  & 11 & 8 & 8 \\
\hline
\textbf{Primal solution obtained (B\&B)} & 4  & 8  & 12  & 11 & 8 & 8  \\
\textbf{Primal solution obtained (ILS)} & 9  & 8  & 12  & 11 & 8 & 8  \\
\hline
\end{tabular}
\caption{Number of instances per instance group and number of instances solved by \gls{acr:ils} versus \gls{acr:bab}.}
\label{tbl:comp-benchmark}
\end{table}

While the \gls{acr:ils} framework finds feasible solutions to all instances, the benchmark cannot solve $5$ out of $9$ instances from group C1. For these instances, no warmstart solution is available and the \gls{acr:bab} algorithm does not find any feasible solutions.

\begin{figure}[!t]
\centering
\resizebox{0.5\textwidth}{!}{
\begin{tikzpicture}

\definecolor{darkgray}{RGB}{169,169,169}
\definecolor{darkgray176}{RGB}{176,176,176}
\definecolor{steelblue31119180}{RGB}{31,119,180}

\begin{axis}[
tick align=outside,
tick pos=left,
x grid style={darkgray176},
xmin=0.5, xmax=6.5,
xtick style={color=black},
xtick={1,2,3,4,5,6},
xticklabels={C1,C2,R1,R2,RC1,RC2},
y grid style={darkgray},
ylabel={Primal improvement ILS vs. B\&B [\%]},
ymajorgrids,
ymin=-100, ymax=10,
ytick style={color=black}
]
\path [draw=black, fill=steelblue31119180]
(axis cs:0.6,-53.9455442251537)
--(axis cs:1.4,-53.9455442251537)
--(axis cs:1.4,-40.8381247716671)
--(axis cs:0.6,-40.8381247716671)
--(axis cs:0.6,-53.9455442251537)
--cycle;
\addplot [black]
table {
1 -53.9455442251537
1 -57.3332461435743
};
\addplot [black]
table {
1 -40.8381247716671
1 -40.8381247716671
};
\addplot [black]
table {
0.8 -57.3332461435743
1.2 -57.3332461435743
};
\addplot [black]
table {
0.8 -40.8381247716671
1.2 -40.8381247716671
};
\addplot [black, opacity=0.5, mark=o, mark size=3, mark options={solid,fill opacity=0}, only marks]
table {
1 -20.7922525243078
};
\path [draw=black, fill=steelblue31119180]
(axis cs:1.6,-44.9380489537407)
--(axis cs:2.4,-44.9380489537407)
--(axis cs:2.4,-41.247242674455)
--(axis cs:1.6,-41.247242674455)
--(axis cs:1.6,-44.9380489537407)
--cycle;
\addplot [black]
table {
2 -44.9380489537407
2 -47.3153689634982
};
\addplot [black]
table {
2 -41.247242674455
2 -38.9883678251522
};
\addplot [black]
table {
1.8 -47.3153689634982
2.2 -47.3153689634982
};
\addplot [black]
table {
1.8 -38.9883678251522
2.2 -38.9883678251522
};
\path [draw=black, fill=steelblue31119180]
(axis cs:2.6,-20.8990749392198)
--(axis cs:3.4,-20.8990749392198)
--(axis cs:3.4,-5.73861191539586)
--(axis cs:2.6,-5.73861191539586)
--(axis cs:2.6,-20.8990749392198)
--cycle;
\addplot [black]
table {
3 -20.8990749392198
3 -26.4040336970633
};
\addplot [black]
table {
3 -5.73861191539586
3 -1.93437986537498
};
\addplot [black]
table {
2.8 -26.4040336970633
3.2 -26.4040336970633
};
\addplot [black]
table {
2.8 -1.93437986537498
3.2 -1.93437986537498
};
\path [draw=black, fill=steelblue31119180]
(axis cs:3.6,-4.70751548611178)
--(axis cs:4.4,-4.70751548611178)
--(axis cs:4.4,0.171295861017625)
--(axis cs:3.6,0.171295861017625)
--(axis cs:3.6,-4.70751548611178)
--cycle;
\addplot [black]
table {
4 -4.70751548611178
4 -8.56336406217553
};
\addplot [black]
table {
4 0.171295861017625
4 2.08869346258206
};
\addplot [black]
table {
3.8 -8.56336406217553
4.2 -8.56336406217553
};
\addplot [black]
table {
3.8 2.08869346258206
4.2 2.08869346258206
};
\path [draw=black, fill=steelblue31119180]
(axis cs:4.6,-17.2283790336649)
--(axis cs:5.4,-17.2283790336649)
--(axis cs:5.4,-6.54885768523907)
--(axis cs:4.6,-6.54885768523907)
--(axis cs:4.6,-17.2283790336649)
--cycle;
\addplot [black]
table {
5 -17.2283790336649
5 -25.2958706215113
};
\addplot [black]
table {
5 -6.54885768523907
5 -1.43727264804131
};
\addplot [black]
table {
4.8 -25.2958706215113
5.2 -25.2958706215113
};
\addplot [black]
table {
4.8 -1.43727264804131
5.2 -1.43727264804131
};
\path [draw=black, fill=steelblue31119180]
(axis cs:5.6,-8.59260041703242)
--(axis cs:6.4,-8.59260041703242)
--(axis cs:6.4,-2.48563155450922)
--(axis cs:5.6,-2.48563155450922)
--(axis cs:5.6,-8.59260041703242)
--cycle;
\addplot [black]
table {
6 -8.59260041703242
6 -9.85177891916508
};
\addplot [black]
table {
6 -2.48563155450922
6 0.883029636978736
};
\addplot [black]
table {
5.8 -9.85177891916508
6.2 -9.85177891916508
};
\addplot [black]
table {
5.8 0.883029636978736
6.2 0.883029636978736
};
\addplot [black]
table {
0.6 -50.1681962199002
1.4 -50.1681962199002
};
\addplot [black]
table {
1.6 -43.0256965424752
2.4 -43.0256965424752
};
\addplot [black]
table {
2.6 -12.724457210514
3.4 -12.724457210514
};
\addplot [black]
table {
3.6 -1.5303039044175
4.4 -1.5303039044175
};
\addplot [black]
table {
4.6 -13.0677031148792
5.4 -13.0677031148792
};
\addplot [black]
table {
5.6 -6.57415695975114
6.4 -6.57415695975114
};
\end{axis}

\end{tikzpicture}}
\caption{Relative primal solution value decrease of ILS vs. \gls{acr:bab} if both methods yield a feasible solution}
\label{fig:solution_quality}
\end{figure}
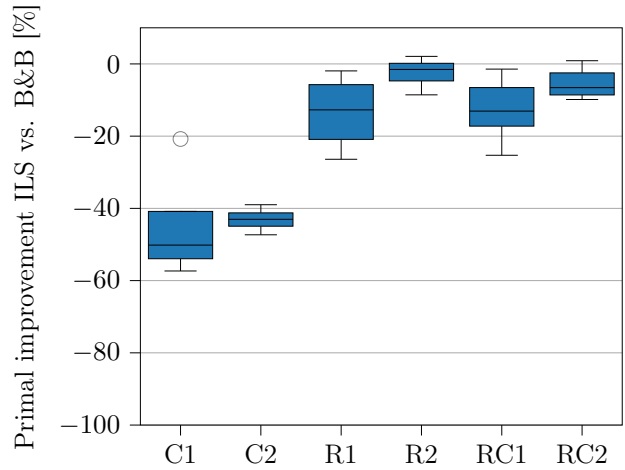

Figure~\ref{fig:solution_quality} shows the relative change in solution value when comparing our \gls{acr:ils} with the \gls{acr:bab}. Negative values correspond to the \gls{acr:ils} finding better solutions than the \gls{acr:bab} algorithm.

\begin{result}
    The \gls{acr:ils} framework finds primal feasible solutions to all instances within an hour of runtime and improves the solution values compared to its benchmark by up to $60\%$. The benchmarking approach is sensitive to being initialized and fails to solve instances where no initial solution is available.
\end{result}

Note that the \gls{acr:ils} may yield slightly higher solution values even when warmstarted at the same configuration because it leverages heuristic solutions to the subproblem and additionally does not allow to partially charge at dynamic charging stations.

\begin{figure}[!t]
\centering
\resizebox{\textwidth}{!}{\input{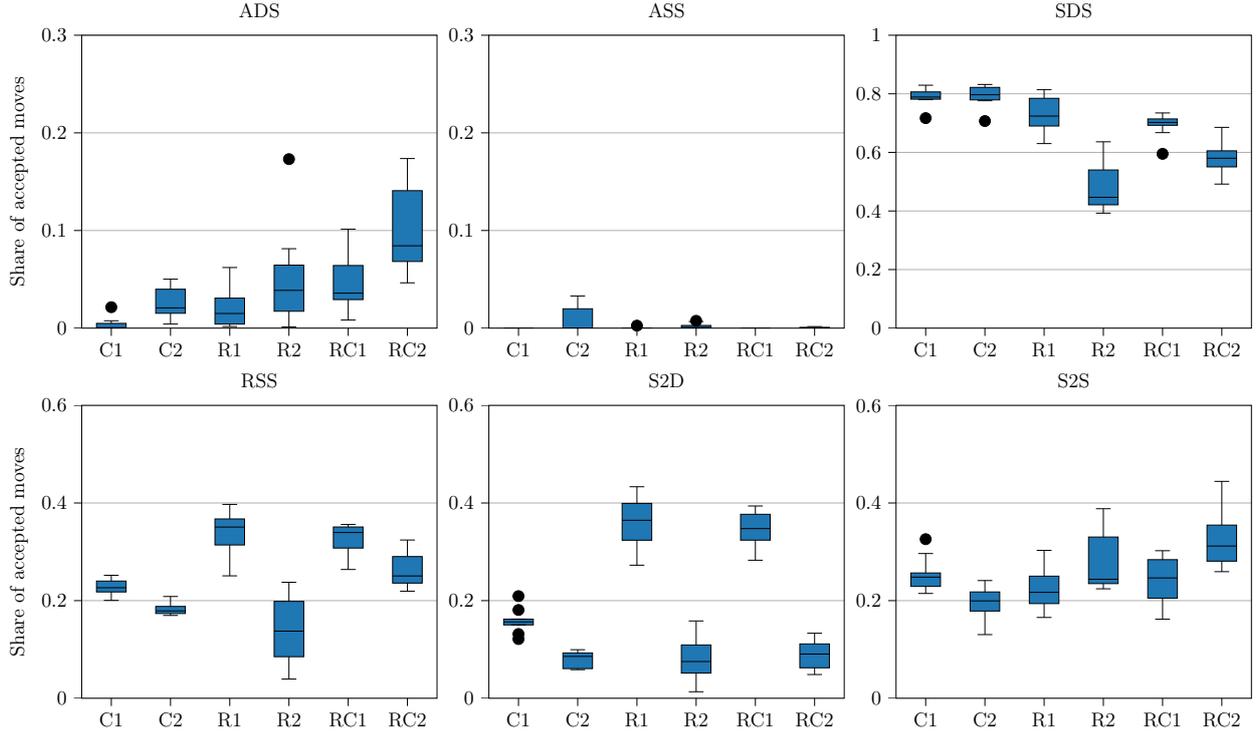}}
\caption{Share of accepted moves per sample-based local search operator and per instance group}
\label{fig:operator_efficiency}
\end{figure}

Thus, our \gls{acr:ils} framework is capable of finding better solutions than the benchmark and is furthermore capable of solving a variety of instances without relying on initial information such as a warmstart solution. Additionally, we demonstrate that the sample-based operators of our local search are suitable to intensify a given feasible solution. Figure ~\ref{fig:operator_efficiency} yields the share of accepted moves per operator and instance class. In this context, we count a move by an operator as accepted, if the operator provides an $\epsilon$-improved solution. Figure~\ref{fig:operator_efficiency} demonstrates that all local search operators provide a benefit to the sample-based local search by improving given solutions. In particular, \gls{acr:strip-dyn} yields a solution value improvement in up to $80\%$ of moves. On the other hand, \gls{acr:add-stat} only very rarely improves a given solution. This observation is due to the trade-off between energy consumption cost and infrastructure cost in which travel distance reductions at the cost of higher infrastructure costs are rarely worthwhile. Moreover, \gls{acr:add-dyn} is particularly effective on the instances which are partially distributed at random and partially clustered. Here, the operator yields $\epsilon$-improved solutions in up to $17\%$ of the executions. The remaining operators yield accepted moves in up to $40\%$ of executions across instance classes.    

\subsection{Managerial Study} \label{sec:managerial-insights}
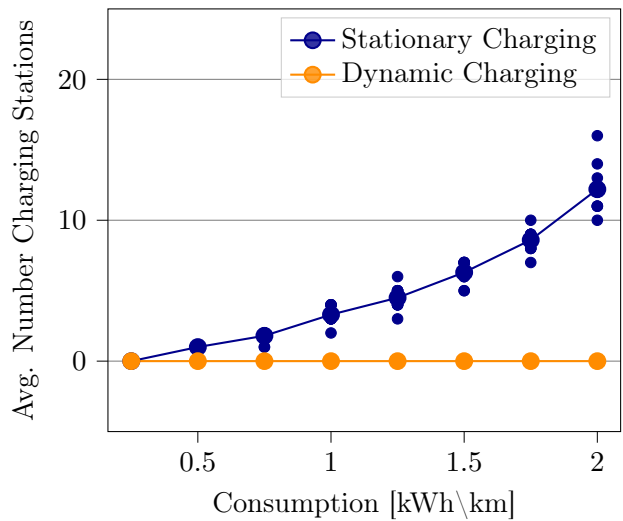
\begin{figure}[!tb]
\centering
\resizebox{0.5\textwidth}{!}{
\begin{tikzpicture}

\definecolor{darkblue}{RGB}{0,0,139}
\definecolor{darkgray176}{RGB}{176,176,176}
\definecolor{darkorange}{RGB}{255,140,0}
\definecolor{gray}{RGB}{128,128,128}
\definecolor{lightgray204}{RGB}{204,204,204}

\begin{axis}[
legend cell align={left},
legend style={fill opacity=0.8, draw opacity=1, text opacity=1, draw=lightgray204},
tick align=outside,
tick pos=left,
x grid style={darkgray176},
xlabel={Consumption [kWh\textbackslash km]},
xmin=0.1625, xmax=2.0875,
xtick style={color=black},
y grid style={gray},
ylabel={Avg. Number Charging Stations},
ymajorgrids,
ymin=-5, ymax=25,
ytick style={color=black}
]
\addplot [draw=darkblue, fill=darkblue, forget plot, mark=*, only marks]
table{%
x  y
0.25 0
0.25 0
0.25 0
0.25 0
0.25 0
0.25 0
0.25 0
0.25 0
0.25 0
0.25 0
};
\addplot [draw=darkblue, fill=darkblue, forget plot, mark=*, only marks]
table{%
x  y
0.5 1
0.5 1
0.5 1
0.5 1
0.5 1
0.5 1
0.5 1
0.5 1
0.5 1
0.5 1
};
\addplot [draw=darkblue, fill=darkblue, forget plot, mark=*, only marks]
table{%
x  y
0.75 2
0.75 2
0.75 1
0.75 2
0.75 1
0.75 2
0.75 2
0.75 2
0.75 2
0.75 2
};
\addplot [draw=darkblue, fill=darkblue, forget plot, mark=*, only marks]
table{%
x  y
1 4
1 3
1 2
1 4
1 3
1 3
1 3
1 4
1 3
1 4
};
\addplot [draw=darkblue, fill=darkblue, forget plot, mark=*, only marks]
table{%
x  y
1.25 5
1.25 4
1.25 3
1.25 5
1.25 4
1.25 4
1.25 4
1.25 5
1.25 5
1.25 6
};
\addplot [draw=darkblue, fill=darkblue, forget plot, mark=*, only marks]
table{%
x  y
1.5 7
1.5 7
1.5 5
1.5 6
1.5 5
1.5 6
1.5 6
1.5 7
1.5 7
1.5 7
};
\addplot [draw=darkblue, fill=darkblue, forget plot, mark=*, only marks]
table{%
x  y
1.75 8
1.75 9
1.75 7
1.75 8
1.75 8
1.75 9
1.75 9
1.75 10
1.75 9
1.75 9
};
\addplot [draw=darkblue, fill=darkblue, forget plot, mark=*, only marks]
table{%
x  y
2 11
2 16
2 11
2 12
2 10
2 12
2 11
2 13
2 12
2 14
};
\addplot [thick, darkblue, mark=*, mark size=3, mark options={solid}]
table {%
0.25 0
0.5 1
0.75 1.8
1 3.3
1.25 4.5
1.5 6.3
1.75 8.6
2 12.2
};
\addlegendentry{Stationary Charging}
\addplot [thick, darkorange, mark=*, mark size=3, mark options={solid}]
table {%
0.25 0
0.5 0
0.75 0
1 0
1.25 0
1.5 0
1.75 0
2 0
};
\addlegendentry{Dynamic Charging}
\end{axis}

\end{tikzpicture}}
\caption{Average number of stationary and dynamic charging stations in the solutions across the different scenarios and for different consumption patterns.}
\label{fig:dyn_in_base_case}
\end{figure}

Given the current cost structure as it is reflected in the base parameters, we investigate the competitiveness of dynamic inductive charging versus stationary inductive charging. Figure~\ref{fig:dyn_in_base_case} shows the average number of charging stations differentiated by stationary charging and dynamic charging with the base cost structure for various given fleet consumption rates.

The required number of charging stations rises with increasing demand due to the vehicles' more frequent need to recharge and the increased relative cost of deviating from the shortest route to utilize nearby chargers. In this context, our experiments show that satisfying the increased demand for charging infrastructure with additional stationary charging stations dominates solutions in which additional dynamic charging stations are constructed. The reason for this observation are the high infrastructure cost for dynamic charging stations compared to stationary charging stations and the necessity to deploy dynamic charging stations at large-scale in order to exchange singular stationary charging stations while maintaining feasibility. 

\begin{result}
    Given the cost associated to deploying inductive charging infrastructure as of now, dynamic inductive charging is not cost-competitive compared to stationary inductive charging.    
\end{result}

Moreover, we allow the vehicles of the shuttle-fleet to deviate from their shortest route between servicing consecutive stops in their timetables. Amongst others, this feature sets us apart from most comparable microscopic studies such as e.g., \cite{KoJang2013} and \cite{LiHeEtAl2024}. In Figure~\ref{fig:value_of_deviation}, we demonstrate that allowing such deviations yields a net cost benefit.

Given the cost for infrastructure components, i.e., stationary and dynamic charging stations, and the cost for consuming and recharging energy of $0.3$\euro, we observe an average reduction of the total cost of up to $3.5\%$.

\begin{result}
    Allowing the vehicles to deviate from the shortest route within their given timetables yields a cost reduction of up to $3.5\%$. While the increased travel distances leads to elevated costs for energy consumption and recharging of up to $5\%$, the decrease in infrastructure cost of up to $8\%$ outweighs this effect and leads to net savings. 
\end{result} 

The application of decentralized energy storage connected to solar panels as they are deployed in \gls{acr:milas}, has the potential to significantly decrease the average cost per kilowatthour and therefore increase the potential of such detours.

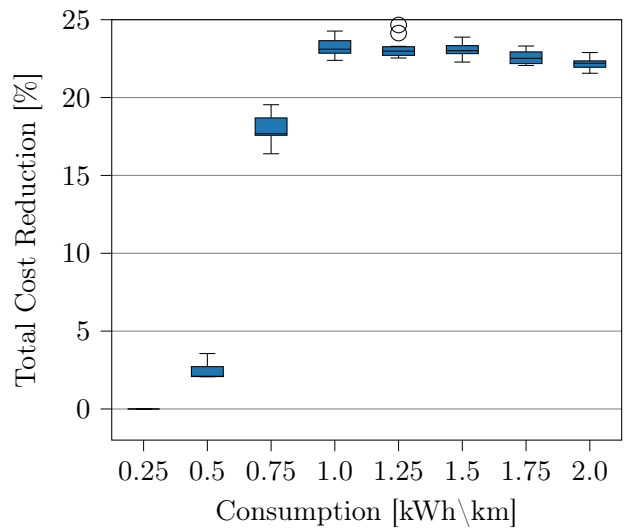
\begin{figure}[!b]
\centering
\resizebox{0.5\textwidth}{!}{
\begin{tikzpicture}

\definecolor{darkgray176}{RGB}{176,176,176}
\definecolor{gray}{RGB}{128,128,128}
\definecolor{steelblue31119180}{RGB}{31,119,180}

\begin{axis}[
tick align=outside,
tick pos=left,
x grid style={darkgray176},
xlabel={Consumption [kWh\textbackslash km]},
xmin=0.5, xmax=8.5,
xtick style={color=black},
xtick={1,2,3,4,5,6,7,8},
xticklabels={0.25,0.5,0.75,1.0,1.25,1.5,1.75,2.0},
y grid style={gray},
ylabel={Total Cost Reduction [\%]},
ymajorgrids,
ymin=-2, ymax=25,
ytick style={color=black}
]
\path [draw=black, fill=steelblue31119180]
(axis cs:0.75,0)
--(axis cs:1.25,0)
--(axis cs:1.25,0)
--(axis cs:0.75,0)
--(axis cs:0.75,0)
--cycle;
\addplot [black]
table {%
1 0
1 0
};
\addplot [black]
table {%
1 0
1 0
};
\addplot [black]
table {%
0.875 0
1.125 0
};
\addplot [black]
table {%
0.875 0
1.125 0
};
\path [draw=black, fill=steelblue31119180]
(axis cs:1.75,2.08)
--(axis cs:2.25,2.08)
--(axis cs:2.25,2.7175)
--(axis cs:1.75,2.7175)
--(axis cs:1.75,2.08)
--cycle;
\addplot [black]
table {%
2 2.08
2 2.07
};
\addplot [black]
table {%
2 2.7175
2 3.56
};
\addplot [black]
table {%
1.875 2.07
2.125 2.07
};
\addplot [black]
table {%
1.875 3.56
2.125 3.56
};
\path [draw=black, fill=steelblue31119180]
(axis cs:2.75,17.5775)
--(axis cs:3.25,17.5775)
--(axis cs:3.25,18.6925)
--(axis cs:2.75,18.6925)
--(axis cs:2.75,17.5775)
--cycle;
\addplot [black]
table {%
3 17.5775
3 16.39
};
\addplot [black]
table {%
3 18.6925
3 19.54
};
\addplot [black]
table {%
2.875 16.39
3.125 16.39
};
\addplot [black]
table {%
2.875 19.54
3.125 19.54
};
\path [draw=black, fill=steelblue31119180]
(axis cs:3.75,22.8475)
--(axis cs:4.25,22.8475)
--(axis cs:4.25,23.6525)
--(axis cs:3.75,23.6525)
--(axis cs:3.75,22.8475)
--cycle;
\addplot [black]
table {%
4 22.8475
4 22.39
};
\addplot [black]
table {%
4 23.6525
4 24.27
};
\addplot [black]
table {%
3.875 22.39
4.125 22.39
};
\addplot [black]
table {%
3.875 24.27
4.125 24.27
};
\path [draw=black, fill=steelblue31119180]
(axis cs:4.75,22.71)
--(axis cs:5.25,22.71)
--(axis cs:5.25,23.265)
--(axis cs:4.75,23.265)
--(axis cs:4.75,22.71)
--cycle;
\addplot [black]
table {%
5 22.71
5 22.54
};
\addplot [black]
table {%
5 23.265
5 23.29
};
\addplot [black]
table {%
4.875 22.54
5.125 22.54
};
\addplot [black]
table {%
4.875 23.29
5.125 23.29
};
\addplot [black, mark=o, mark size=3, mark options={solid,fill opacity=0}, only marks]
table {%
5 24.14
5 24.66
};
\path [draw=black, fill=steelblue31119180]
(axis cs:5.75,22.82)
--(axis cs:6.25,22.82)
--(axis cs:6.25,23.335)
--(axis cs:5.75,23.335)
--(axis cs:5.75,22.82)
--cycle;
\addplot [black]
table {%
6 22.82
6 22.28
};
\addplot [black]
table {%
6 23.335
6 23.88
};
\addplot [black]
table {%
5.875 22.28
6.125 22.28
};
\addplot [black]
table {%
5.875 23.88
6.125 23.88
};
\path [draw=black, fill=steelblue31119180]
(axis cs:6.75,22.1825)
--(axis cs:7.25,22.1825)
--(axis cs:7.25,22.9325)
--(axis cs:6.75,22.9325)
--(axis cs:6.75,22.1825)
--cycle;
\addplot [black]
table {%
7 22.1825
7 22.06
};
\addplot [black]
table {%
7 22.9325
7 23.31
};
\addplot [black]
table {%
6.875 22.06
7.125 22.06
};
\addplot [black]
table {%
6.875 23.31
7.125 23.31
};
\path [draw=black, fill=steelblue31119180]
(axis cs:7.75,21.945)
--(axis cs:8.25,21.945)
--(axis cs:8.25,22.355)
--(axis cs:7.75,22.355)
--(axis cs:7.75,21.945)
--cycle;
\addplot [black]
table {%
8 21.945
8 21.56
};
\addplot [black]
table {%
8 22.355
8 22.89
};
\addplot [black]
table {%
7.875 21.56
8.125 21.56
};
\addplot [black]
table {%
7.875 22.89
8.125 22.89
};
\addplot [black]
table {%
0.75 0
1.25 0
};
\addplot [black]
table {%
1.75 2.1
2.25 2.1
};
\addplot [black]
table {%
2.75 17.68
3.25 17.68
};
\addplot [black]
table {%
3.75 23.105
4.25 23.105
};
\addplot [black]
table {%
4.75 22.98
5.25 22.98
};
\addplot [black]
table {%
5.75 23.01
6.25 23.01
};
\addplot [black]
table {%
6.75 22.52
7.25 22.52
};
\addplot [black]
table {%
7.75 22.2
8.25 22.2
};
\end{axis}

\end{tikzpicture}}
\caption{Total cost reduction with variable energy price curve compared to constant energy price}
\label{fig:value_of_energy_storage}
\end{figure}

Equipping the system with solar panels and storing the generated energy in a decentral battery helps to avoid purchasing energy via the power grid at full exposure to market prices that include the margin for the providers. Instead, we can charge from the solar panels or the connected battery. We model this effect of decentralized energy storage as a piecewise constant function over hourly intervals in which energy costs $0.3$\euro \, per kWh until 8am when the storage is empty and recharging relies on the grid, then the price drops to $0.04$\euro \, when the solar panels take over, and finally at 1pm the energy is provided by the energy storage, the price reaches $0.2067$\euro \, and remains at this level throughout the rest of the day. Figure~\ref{fig:value_of_energy_storage} shows the effects of considering such decentral energy systems, represented by a piecewise constant energy price function over hourly intervals, on the total cost in our study. Specifically, the figure shows the total cost reduction if we leverage a decentralized energy storage with solar panels that has the described effect on the energy price curve versus the case in which energy prices are constant and given by the power grid provider. 

\begin{result}
    The value of photovoltaic systems with decentralized energy storage solutions exceeds $20\%$ of the total system cost for fleet consumptions of $1.0$ kWh\slash km or more. In some scenarios the value reaches $25\%$.
\end{result}

Note that we derived this result based on an early shift from $5$AM to $1$PM. This shift contains the significant decrease in the energy cost function such that the value of the system in relation to the total cost may be smaller when shifting the operating time of the fleet during the day. However, the value can only increase further when evaluating an additional shift and measuring the value of the decentralized system in absolute terms.

\begin{figure}[!t]
\centering
\resizebox{0.5\textwidth}{!}{
\begin{tikzpicture}

\definecolor{darkblue}{RGB}{0,0,139}
\definecolor{darkgray176}{RGB}{176,176,176}
\definecolor{deepskyblue}{RGB}{0,191,255}
\definecolor{gray}{RGB}{128,128,128}
\definecolor{lightgray204}{RGB}{204,204,204}
\definecolor{royalblue}{RGB}{65,105,225}

\begin{axis}[
legend cell align={left},
legend style={fill opacity=0.8, draw opacity=1, text opacity=1, draw=lightgray204},
tick align=outside,
tick pos=left,
x grid style={darkgray176},
xlabel={Consumption [kWh\textbackslash km]},
xmin=0.1625, xmax=2.0875,
xtick style={color=black},
y grid style={gray},
ylabel={Value of Deviating [\% of total cost]},
ymajorgrids,
ymin=-10, ymax=10,
ytick style={color=black}
]
\addplot [draw=deepskyblue, fill=deepskyblue, mark=*, only marks]
table{%
x  y
0.25 0
0.25 0
0.25 0
0.25 0
0.25 0
0.25 0
0.25 0
0.25 0
0.25 0
0.25 0
};
\addlegendentry{Infrastructure Cost}
\addplot [draw=royalblue, fill=royalblue, mark=*, only marks]
table{%
x  y
0.25 0.02
0.25 0.01
0.25 0.01
0.25 0.01
0.25 0.01
0.25 0.01
0.25 0.01
0.25 0.01
0.25 0.01
0.25 0.01
};
\addlegendentry{Routing Cost}
\addplot [draw=deepskyblue, fill=deepskyblue, forget plot, mark=*, only marks]
table{%
x  y
0.5 0
0.5 0
0.5 0
0.5 0
0.5 0
0.5 0
0.5 0
0.5 -4.98
0.5 0
0.5 0
};
\addplot [draw=royalblue, fill=royalblue, forget plot, mark=*, only marks]
table{%
x  y
0.5 1.25
0.5 1.94
0.5 1.58
0.5 1.23
0.5 1.02
0.5 1.22
0.5 1.22
0.5 1.89
0.5 1.69
0.5 1.22
};
\addplot [draw=deepskyblue, fill=deepskyblue, forget plot, mark=*, only marks]
table{%
x  y
0.75 -8.15
0.75 -2.74
0.75 -5.5
0.75 -5.48
0.75 -8.23
0.75 -2.81
0.75 -5.44
0.75 -7.67
0.75 -5.3
0.75 -5.43
};
\addplot [draw=royalblue, fill=royalblue, forget plot, mark=*, only marks]
table{%
x  y
0.75 3.86
0.75 2.78
0.75 3.9
0.75 3.83
0.75 4.77
0.75 2.42
0.75 2.7
0.75 5
0.75 2.44
0.75 4.09
};
\addplot [draw=deepskyblue, fill=deepskyblue, forget plot, mark=*, only marks]
table{%
x  y
1 -3.65
1 -5.27
1 -5.38
1 -3.6
1 -3.67
1 -5.38
1 -3.64
1 -5.11
1 -6.87
1 -5.26
};
\addplot [draw=royalblue, fill=royalblue, forget plot, mark=*, only marks]
table{%
x  y
1 2.31
1 2.52
1 3.4
1 1.9
1 3
1 2.28
1 2.18
1 1.92
1 3.26
1 3.29
};
\addplot [draw=deepskyblue, fill=deepskyblue, forget plot, mark=*, only marks]
table{%
x  y
1.25 -4.07
1.25 -5.24
1.25 -5.32
1.25 -5.28
1.25 -5.36
1.25 -5.34
1.25 -4.04
1.25 -5.1
1.25 -5.15
1.25 -3.93
};
\addplot [draw=royalblue, fill=royalblue, forget plot, mark=*, only marks]
table{%
x  y
1.25 2.1
1.25 2.16
1.25 3.34
1.25 1.63
1.25 1.65
1.25 2.36
1.25 1.36
1.25 1.44
1.25 2.16
1.25 1.34
};
\addplot [draw=deepskyblue, fill=deepskyblue, forget plot, mark=*, only marks]
table{%
x  y
1.5 -4.28
1.5 -4.13
1.5 -6.21
1.5 -5.27
1.5 -6.33
1.5 -5.26
1.5 -4.24
1.5 -6
1.5 -4.11
1.5 -5.18
};
\addplot [draw=royalblue, fill=royalblue, forget plot, mark=*, only marks]
table{%
x  y
1.5 1.54
1.5 2.12
1.5 2.1
1.5 2.28
1.5 2.94
1.5 1.4
1.5 1.59
1.5 2.43
1.5 1.93
1.5 1.91
};
\addplot [draw=deepskyblue, fill=deepskyblue, forget plot, mark=*, only marks]
table{%
x  y
1.75 -6.12
1.75 -5.07
1.75 -7.57
1.75 -6.04
1.75 -5.22
1.75 -4.34
1.75 -3.49
1.75 -4.97
1.75 -5.05
1.75 -5.95
};
\addplot [draw=royalblue, fill=royalblue, forget plot, mark=*, only marks]
table{%
x  y
1.75 2.75
1.75 2.74
1.75 2.88
1.75 1.9
1.75 2.03
1.75 1.29
1.75 1.07
1.75 1.55
1.75 1.53
1.75 2.82
};
\addplot [draw=deepskyblue, fill=deepskyblue, forget plot, mark=*, only marks]
table{%
x  y
2 -6.6
2 -2.15
2 -5.74
2 -5.8
2 -5.88
2 -5.09
2 -5.12
2 -6.24
2 -5.66
2 -5
};
\addplot [draw=royalblue, fill=royalblue, forget plot, mark=*, only marks]
table{%
x  y
2 2.11
2 1.64
2 1.98
2 1.17
2 2.4
2 1.66
2 1.25
2 2.44
2 2.32
2 1.79
};
\addplot [thick, darkblue, mark=diamond*, mark size=3, mark options={solid}]
table {%
0.25 0.011
0.5 0.928
0.75 -2.096
1 -2.177
1.25 -2.929
1.5 -3.077
1.75 -3.326
2 -3.452
};
\addlegendentry{Avg. Net Value}
\end{axis}

\end{tikzpicture}}
\caption{Average value of deviating from the shortest route for total shuttle fleet as a proportion of total cost across scenarios and for different consumption patterns.}
\label{fig:value_of_deviation}
\end{figure}

We furthermore investigate under which conditions dynamic charging becomes an economically sustainable addition to stationary inductive charging. For this purpose, we assume that the cost for stationary charging stations $\InvestmentCost_{\Station}, \medspace \Station \in \SetOfStaticStations$ decreases by scaling effects such that it only contains material- and construction cost. The same applies to the fix cost for dynamic charging stations $\InvestmentCost_{\Station}, \medspace \Station \in \SetOfDynamicStations$. We then scale the variable base cost per meter of segment of \euro{0.19} between $10\%$ and $1\%$. 

\begin{figure}[!t]
\centering
\resizebox{0.5\textwidth}{!}{
\begin{tikzpicture}

\definecolor{cornflowerblue84158205}{RGB}{84,158,205}
\definecolor{darkgray176}{RGB}{176,176,176}
\definecolor{lightgray204}{RGB}{204,204,204}
\definecolor{skyblue147196222}{RGB}{147,196,222}
\definecolor{steelblue37117183}{RGB}{37,117,183}
\definecolor{teal873145}{RGB}{8,73,145}

\begin{axis}[
legend cell align={left},
legend style={
  fill opacity=0.8,
  draw opacity=1,
  text opacity=1,
  at={(0.03,0.97)},
  anchor=north west,
  draw=lightgray204,
},
tick align=outside,
tick pos=left,
x dir=reverse,
x grid style={darkgray176},
xlabel={Cost per Segment Meter [\% of Base Case]},
xmin=0.1045, xmax=1.9855,
xtick style={color=black},
xtick={0.19,0.38,0.57,0.76,0.95,1.14,1.33,1.52,1.71,1.9},
xticklabels={1\%,2\%,3\%,4\%,5\%,6\%,7\%,8\%,9\%,10\%},
y grid style={darkgray176},
ylabel={Avg. Total Segment Length [km]},
ymin=-0.017595, ymax=0.369495,
ytick style={color=black}
]
\addplot [semithick, skyblue147196222, mark=*, mark size=3, mark options={solid}]
table {%
0.19 0.3519
0.38 0.3309
0.57 0.2306
0.76 0.1296
0.95 0
1.14 0
1.33 0
1.52 0
1.71 0
1.9 0
};
\addlegendentry{$0.5\,\mathrm{kWh/km}$}
\addplot [semithick, cornflowerblue84158205, mark=*, mark size=3, mark options={solid}]
table {%
0.19 0.1346
0.38 0.0569
0.57 0.0502
0.76 0.0247
0.95 0
1.14 0.0146
1.33 0
1.52 0
1.71 0
1.9 0
};
\addlegendentry{$1.0\,\mathrm{kWh/km}$}
\addplot [semithick, steelblue37117183, mark=*, mark size=3, mark options={solid}]
table {%
0.19 0.176
0.38 0.0501
0.57 0.0543
0.76 0.0293
0.95 0.0167
1.14 0.0116
1.33 0
1.52 0
1.71 0
1.9 0
};
\addlegendentry{$1.5\,\mathrm{kWh/km}$}
\addplot [semithick, teal873145, mark=*, mark size=3, mark options={solid}]
table {%
0.19 0.2306
0.38 0.1151
0.57 0.06
0.76 0.0617
0.95 0.0296
1.14 0.015
1.33 0
1.52 0.015
1.71 0
1.9 0
};
\addlegendentry{$2.0\,\mathrm{kWh/km}$}
\end{axis}

\end{tikzpicture}}
\caption{Average total deployed segment length for decreasing variable cost for different vehicle energy consumptions}
\label{fig:penetration_threshold}
\end{figure}
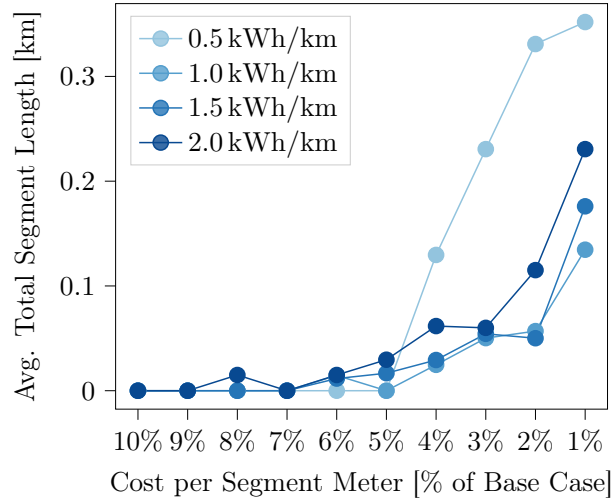

Figure~\ref{fig:penetration_threshold} shows the average total deployed segment length in the $10$ evaluated scenarios given different consumption patterns. The figure shows that for variable prices of less than $5\%$ of the base cost per meter segment, dynamic charging becomes an economically viable addition to stationary inductive charging stations and our algorithm starts to deploy dynamic charging segments at an increasing scale. 

The determining factor for the threshold of penetration is the relation of costs for stationary versus dynamic charging and the capability of dynamic charging stations to satisfy the charging need of the fleet. Hence, we observe that the penetration for the scenario with a low fleet consumption reaches the highest penetration level. The lower overall re-charging requirement of the fleet in this scenario leads to a significantly increased competitiveness of dynamic charging stations. Thus, the lower the consumption compared to the charging rates is, and the lower the variable costs for dynamic charging segments are, the better is their competitiveness. However, even in these competitive scenarios, dynamic inductive charging complements but does not dominate stationary inductive charging because the solutions of the respective scenarios comprise both kinds of charging infrastructure.

\begin{result}
    The penetration of dynamic charging segments increases when the cost per meter segment is smaller than $5\%$ of the base cost of \euro{0.19}. The magnitude of this effect is independent of the fleet consumption.
\end{result}

\section{Conclusion} \label{sec:conclusion}
We introduced a novel microscopic planning problem to determine the optimal deployment of heterogeneous stationary and dynamic inductive charging infrastructure for electric shuttle fleets. Our formulation accounts for fixed stop sequences and timetables, while allowing vehicles to make controlled detours for charging to improve infrastructure utilization and flexibility. The proposed model explicitly incorporates variable energy prices and is readily extensible to heterogeneous vehicle fleets and charging technologies.

To address the computational complexity of realistic instances, we provided a \gls{acr:mip} formulation and developed a tailored metaheuristic based on \gls{acr:ils}. This framework solves vehicle-level subproblems through dynamic programming using a state-of-the-art bi-directional label-setting algorithm with lazy dominance checks. We generated benchmark instances based on adapted \gls{acr:evrptw-pr} datasets and demonstrated that our approach outperforms a warm-started commercial solver by up to 60\% in terms of solution quality.

We also conducted a real-world case study inspired by the \gls{acr:milas} project, in which autonomous, inductively charged shuttles were deployed on a live test track. Insights from operations and consultations with industry partners informed our parameter calibration and scenario design. Results indicate that, under current technological and cost conditions, dynamic inductive charging is not economically viable for shuttle fleets in smaller municipalities when compared to stationary charging alternatives. Nevertheless, allowing moderate detours for charging can reduce total system costs by up to 3.5\%, and integrating decentralized photovoltaic energy storage systems yields savings of more than 20\% --- highlighting the potential of combining infrastructure flexibility with renewable energy integration.

\singlespacing{
\bibliographystyle{apalike}
\bibliography{./references}}

\begin{thebibliography}{}

\bibitem[Alwesabi et~al., 2020]{AlwesabiWangEtAl2020}
Alwesabi, Y., Wang, Y., Avalos, R., and Liu, Z. (2020).
\newblock Electric bus scheduling under single depot dynamic wireless charging
  infrastructure planning.
\newblock {\em Energy}, 213:118855.

\bibitem[Broihan et~al., 2022]{BroihanNozinskiEtAl2022}
Broihan, J., Nozinski, I., Pöch, N., and Helber, S. (2022).
\newblock Designing {D}ynamic {I}nductive {C}harging {I}nfrastructures for
  {A}irport {A}prons with {M}ultiple {V}ehicle {T}ypes.
\newblock {\em Energies}, 15(11):4085.

\bibitem[Chen et~al., 2016]{ChenHeEtAl2016}
Chen, Z., He, F., and Yin, Y. (2016).
\newblock Optimal deployment of charging lanes for electric vehicles in
  transportation networks.
\newblock {\em Transportation Research Part B: Methodological}, 91:344--–365.

\bibitem[Chen et~al., 2018]{ChenYinEtAl2018}
Chen, Z., Yin, Y., and Song, Z. (2018).
\newblock A cost-competitiveness analysis of charging infrastructure for
  electric bus operations.
\newblock {\em Transportation Research Part C: Emerging Technologies},
  93:351--366.

\bibitem[{EU Commission}, 2023]{EUCommission2023}
{EU Commission} (2023).
\newblock
  \url{https://ec.europa.eu/commission/presscorner/detail/en/ip_23_762}.
\newblock Accessed 12-January-2024.

\bibitem[Helber et~al., 2018]{HelberBroihanEtAl2018}
Helber, S., Broihan, J., Jang, Y.~J., Hecker, P., and Feuerle, T. (2018).
\newblock {L}ocation {P}lanning for {D}ynamic {W}ireless {C}harging {S}ystems
  for {E}lectric {A}irport {P}assenger {B}uses.
\newblock {\em Energies}, 11(2):258.

\bibitem[Hodgson, 1990]{Hodgson1990}
Hodgson, M.~J. (1990).
\newblock A {F}low-{C}apturing {L}ocation-{A}llocation {M}odel.
\newblock {\em Geographical Analysis}, 22(3):270--279.

\bibitem[Hwang et~al., 2018]{HwangJangEtAl2018}
Hwang, I., Jang, Y.~J., Ko, Y.~D., and Lee, M.~S. (2018).
\newblock System {O}ptimization for {D}ynamic {W}ireless {C}harging {E}lectric
  {V}ehicles {O}perating in a {M}ultiple-{R}oute {E}nvironment.
\newblock {\em IEEE Transactions on Intelligent Transportation Systems},
  19(6):1709--1726.

\bibitem[Iliopoulou and Kepaptsoglou, 2019]{IliopoulouKepaptsoglou2019}
Iliopoulou, C. and Kepaptsoglou, K. (2019).
\newblock Integrated transit route network design and infrastructure planning
  for on-line electric vehicles.
\newblock {\em Transportation Research Part D: Transport and Environment},
  77:178--197.

\bibitem[Kim and Kuby, 2012]{KimKuby2012}
Kim, J.-G. and Kuby, M. (2012).
\newblock The deviation-flow refueling location model for optimizing a network
  of refueling stations.
\newblock {\em Optimization Approaches to Hydrogen Logistics},
  37(6):5406--5420.

\bibitem[Ko and Jang, 2013]{KoJang2013}
Ko, Y.~D. and Jang, Y.~J. (2013).
\newblock The {O}ptimal {S}ystem {D}esign of the {O}nline {E}lectric {V}ehicle
  {U}tilizing {W}ireless {P}ower {T}ransmission {T}echnology.
\newblock {\em IEEE Transactions on Intelligent Transportation Systems},
  14(3):1255--1265.

\bibitem[Ko et~al., 2015]{KoJangEtAl2015}
Ko, Y.~D., Jang, Y.~J., and Lee, M.~S. (2015).
\newblock The optimal economic design of the wireless powered intelligent
  transportation system using genetic algorithm considering nonlinear cost
  function.
\newblock {\em Computers \& Industrial Engineering}, 89:67--79.

\bibitem[Kuby and Lim, 2005]{KubyLim2005}
Kuby, M. and Lim, S. (2005).
\newblock The flow-refueling location problem for alternative-fuel vehicles.
\newblock {\em Socio-Economic Planning Sciences}, 39(2):125--145.

\bibitem[Kuby and Lim, 2007]{KubyLim2007}
Kuby, M. and Lim, S. (2007).
\newblock {L}ocation of {A}lternative-{F}uel {S}tations {U}sing the
  {F}low-{R}efueling {L}ocation {M}odel and {D}ispersion of {C}andidate {S}ites
  on {A}rcs.
\newblock {\em Networks and Spatial Economics}, 7:129--152.

\bibitem[Li and Huang, 2014]{LiHuang2014}
Li, S. and Huang, Y. (2014).
\newblock Heuristic approaches for the flow-based set covering problem with
  deviation paths.
\newblock {\em Transportation Research Part E: Logistics and Transportation
  Review}, 72:144--158.

\bibitem[Li et~al., 2024]{LiHeEtAl2024}
Li, W., He, Y., Hu, S., He, Z., and Ratti, C. (2024).
\newblock Planning dynamic wireless charging infrastructure for battery
  electric bus systems with the joint optimization of charging scheduling.
\newblock {\em Transportation Research Part C: Emerging Technologies},
  159:104469.

\bibitem[Liu and Wang, 2017]{LiuWangEtAl2017}
Liu, H. and Wang, D. Z.~W. (2017).
\newblock Locating multiple types of charging facilities for battery electric
  vehicles.
\newblock {\em Transportation Research Part B: Methodological}, 103:30--55.

\bibitem[Liu and Song, 2017]{LiuSong2017}
Liu, Z. and Song, Z. (2017).
\newblock Robust planning of dynamic wireless charging infrastructure for
  battery electric buses.
\newblock {\em Transportation Research Part C: Emerging Technologies},
  83:77--103.

\bibitem[Liu et~al., 2017]{LiuSongEtAl2017}
Liu, Z., Song, Z., and He, Y. (2017).
\newblock Optimal {D}eployment of {D}ynamic {W}ireless {C}harging {F}acilities
  for an {E}lectric {B}us {S}ystem.
\newblock {\em Transportation Research Record}, 2647(1):100--108.

\bibitem[Majhi et~al., 2022]{MajhiRanjitkarEtAl2022}
Majhi, R.~C., Ranjitkar, P., and Sheng, M. (2022).
\newblock Optimal allocation of dynamic wireless charging facility for electric
  vehicles.
\newblock {\em Transportation Research Part D: Transport and Environment},
  111:103461.

\bibitem[{MIT Climate Portal}, 2023]{MITClimatePortal2023}
{MIT Climate Portal} (2023).
\newblock
  \url{https://climate.mit.edu/ask-mit/why-have-electric-vehicles-won-out-over-hydrogen-cars-so-far}.
\newblock Accessed 12-January-2024.

\bibitem[Mubarak et~al., 2021]{MubarakUesterEtAl2021}
Mubarak, M., Üster, H., Abdelghany, K., and Khodayar, M. (2021).
\newblock Strategic network design and analysis for in-motion wireless charging
  of electric vehicles.
\newblock {\em Transportation Research Part E: Logistics and Transportation
  Review}, 145:102179.

\bibitem[{Münchner Verkehrsgesellschaft}, 2023]{MVG2023}
{Münchner Verkehrsgesellschaft} (2023).
\newblock
  \url{https://www.mvg.de/ueber/das-unternehmen/unternehmensprofil.html}.
\newblock Accessed 12-January-2024.

\bibitem[Riemann et~al., 2015]{RiemannWangEtAl2015}
Riemann, R., Wang, D. Z.~W., and Busch, F. (2015).
\newblock Optimal location of wireless charging facilities for electric
  vehicles: {F}low-capturing location model with stochastic user equilibrium.
\newblock {\em Transportation Research Part C: Emerging Technologies},
  58(A):1--12.

\bibitem[Schiffer and Walther, 2017]{SchifferWalter2017A}
Schiffer, M. and Walther, G. (2017).
\newblock The electric location routing problem with time windows and partial
  recharging.
\newblock {\em European Journal of Operational Research}, 260(3):995--1013.

\bibitem[Schiffer and Walther, 2018]{SchifferWalter2018}
Schiffer, M. and Walther, G. (2018).
\newblock {A}n {A}daptive {L}arge {N}eighborhood {S}earch for the
  {L}ocation-routing {P}roblem with {I}ntra-route {F}acilities.
\newblock {\em Transportation Science}, 52(2):331--352.

\bibitem[Schwerdfeger et~al., 2022]{SchwerdfegerBockEtAl2022}
Schwerdfeger, S., Bock, S., Boysen, N., and Briskorn, D. (2022).
\newblock Optimizing the electrification of roads with charge-while-drive
  technology.
\newblock {\em European Journal of Operational Research}, 299(3):1111--1127.

\bibitem[Song and Cheng, 2024]{SongCheng2024}
Song, M. and Cheng, L. (2024).
\newblock Full cover wireless charging segment location problem with routing in
  space-time-electricity network.
\newblock {\em Sustainable Cities and Society}, 112:105611.

\bibitem[Song et~al., 2023]{SongChengEtAl2023}
Song, M., Cheng, L., and Zhang, Y. (2023).
\newblock Joint location optimization of charging stations and segments in the
  space-time-electricity network: {A}n augmented {L}agrangian relaxation and
  {ADMM}-based decomposition scheme.
\newblock {\em Computers \& Industrial Engineering}, 183:109517.

\bibitem[Stützle and Ruiz, 2018]{StützleRuiz2018}
Stützle, T. and Ruiz, R. (2018).
\newblock Iterated {L}ocal {S}earch.
\newblock In Mart{í}, R., Pardalos, P.~M., and Resende, M. G.~C., editors,
  {\em Handbook of Heuristics}, pages 579--605. Springer, Cham.

\bibitem[Sun et~al., 2020]{SunChenEtAl2020}
Sun, X., Chen, Z., and Yin, Y. (2020).
\newblock Integrated planning of static and dynamic charging infrastructure for
  electric vehicles.
\newblock {\em Transportation Research Part D: Transport and Environment},
  83:102331.

\bibitem[Thomas et~al., 2019]{ThomasEtAl2019}
Thomas, B.~W., Calogiuri, T., and Hewitt, M. (2019).
\newblock An exact bidirectional ${A}^{\star}$ approach for solving
  resource-constrained shortest path problems.
\newblock {\em Networks}, 73(2):187--205.

\bibitem[Tu et~al., 2016]{TuLiEtAl2016}
Tu, W., Li, Q., Fang, Z., Shaw, S., Zhou, B., and Chang, X. (2016).
\newblock Optimizing the locations of electric taxi charging stations: A
  spatial-temporal demand coverage approach.
\newblock {\em Transportation Research Part C: Emerging Technologies},
  65:172--189.

\bibitem[Upchurch et~al., 2009]{UpchurchKubyEtAl2009}
Upchurch, C., Kuby, M., and Lim, S. (2009).
\newblock {A} {M}odel for {L}ocation of {C}apacitated {A}lternative-{F}uel
  {S}tations.
\newblock {\em Geographical Analysis}, 41(1):85--106.

\bibitem[Yıldız et~al., 2016]{YildizArslanEtAl2016}
Yıldız, B., Arslan, O., and Karaşan, O.~E. (2016).
\newblock A branch and price approach for routing and refueling station
  location model.
\newblock {\em European Journal of Operational Research}, 248(3):815--826.

\bibitem[Zeng et~al., 2024]{ZengXieEtAl2024}
Zeng, X., Xie, C., Xu, M., and Chen, Z. (2024).
\newblock Optimal en-route charging station locations for electric vehicles
  with heterogeneous range anxiety.
\newblock {\em Transportation Research Part C: Emerging Technologies},
  158:104459.

\end{thebibliography}
\newpage
\onehalfspacing
\appendix
\begin{appendices}
	\normalsize
        \section{Notation} \label{app:notation} 
        This table summarizes the notation that we introduced to define the problem setting and formulate a \gls{acr:mip} to derive solutions to our problem.

\begin{table}[!h]
	\begin{threeparttable}
		\caption{Notation}
		\label{tbl:mip-notation}
		\centering
		\begin{tabular}{lr}
			\toprule
			Symbol & Meaning \\
			\midrule
			\multicolumn{2}{c}{Sets} \\
			${\SetOfLines}$ & Set of lines \\
            $\SetOfStops$ & Set of stops \\
            $\SetOfSegments$ & Set of segments \\
			${\StartDepot, \EndDepot}$ & Central depot (start and end representation) \\
			$\StopsOfLine$ & Ordered sequence of locations assigned to line $\Line \in \SetOfLines$ \\
            $\SetOfSegments_{\Station}$ & Ordered sequence of segments belonging to $\Station \in \SetOfDynamicStations$ \\
			${\SetOfStations} := \SetOfStaticStations \cup \SetOfDynamicStations$ & Stationary and dynamic candidate charging stations \\
            \midrule
			\multicolumn{2}{c}{Parameters} \\
            ${\InitialSoC}$ & Initial SoC \\
            ${\MinSoC}$ & Minimum SoC \\
			${\MaxSoC}$ & Maximum SoC \\
            $\EnergyPrice_{\texttt{c}}, \EnergyPrice_{\texttt{r}}$ & Energy price per kWh for consumption (c) and recharging (r) \\
            ${\ArcTravelTime}_{\UnbracedArc}$ & Travel time along arc $\Arc$ \\
			${\ArcConsumption}_{\UnbracedArc}$ & Consumption on arc $\Arc$ \\
			${\TwBegin}^{\Vehicle}_{\Vertex}$ & Earliest departure time of vehicle $\Line$ at $\Vertex \in \SetOfVertices$ \\
			${\TwEnd}^{\Vehicle}_{\Vertex}$ & Latest departure time of vehicle $\Line$ at $\Vertex \in \SetOfVertices$ \\
			${\StationChargingRate}_{\Station}$ & Charging rate at station $\Station \in \SetOfStations$ \\
            $\InvestmentCost_{\Station}, \InvestmentCost_{\UnbracedSegment}$ & Cost for one station $\Station \in \SetOfStations$, or segment $\Segment \in \SetOfSegments$ respectively \\
			\midrule
            \multicolumn{2}{c}{Variables} \\
			$\ArcVar^{\Line}_{\UnbracedArc}$ & Line $\Line$ traverses arc $\Arc$ \\
			$\StaticChargerVar_{\Station}$ & Stationary station $\Station \in \SetOfStaticStations$ constructed \\
            $\TransformerVar_{\Station}$ & Dynamic station $\Station \in \SetOfDynamicStations$ is constructed \\
			$\DynamicChargerVar_{\UnbracedSegment}$ & Dynamic segment $\Segment \in \SetOfSegments$ constructed \\
			$\ArrivalTimeVar^{\Line}_{\Vertex}$ & Departure time of line $\Line$ at vertex $\Vertex$ \\
			$\StateOfChargeVar^{\Line}_{\Vertex}$ & State of charge of line $\Line$ departing at vertex $\Vertex$ \\
			$\ChargingTimeVar^{\Line}_{\Station}$ & Charging time of line $\Line$ at station $\Station \in \SetOfStaticStations$ \\
			$\RechargeVar^{\Line}_{\Station}$ & SoC replenished by line $\Line$ at $\Station \in \SetOfStaticStations$ \\
            $\RechargeVar^{\Line}_{\UnbracedSegment}$ & SoC replenished by line $\Line$ at $\Segment \in \SetOfSegments$ \\
			\midrule
			\multicolumn{2}{c}{Other}\\
			$\IncomingArcs(\Vertex)$ & Set of incoming arcs of vertex $\Vertex$\\
			$\OutgoingArcs(\Vertex)$ & Set of outgoing arcs of vertex $\Vertex$\\
			\bottomrule
		\end{tabular}
	\end{threeparttable}
\end{table}
\newpage

        \section{Time-discrete MIP Formulation} \label{app:mip-formulation}     
        A time-discrete \gls{acr:mip} formulation of our problem relying on the notation from Section~\ref{sec:prob-description} and utilizing some additional decision variables summarized in Appendix~\ref{app:notation} is as follows:
\begin{subequations}
    \label{mips:deterministic}
    \setlength{\abovedisplayskip}{0pt}
    \setlength{\belowdisplayskip}{0pt}
    \setlength{\abovedisplayshortskip}{0pt}
    \setlength{\belowdisplayshortskip}{0pt}
    \begin{multline}
        \hfill 
        \min \sum_{\Station \in \SetOfStaticStations} \InvestmentCost_{\Station} \StaticChargerVar_{\Station} + \sum_{\Segment \in \SetOfSegments} \InvestmentCost_{\UnbracedSegment} \DynamicChargerVar_{\UnbracedSegment} + 
        \sum_{\Station \in \SetOfDynamicStations} \InvestmentCost_{\Station} \TransformerVar_{\Station} + \EnergyPrice \medspace \sum_{\Line \in \SetOfLines} \left(\sum_{\Arc \in \SetOfArcs} \ArcConsumption_{\UnbracedArc} \ArcVar^{\Line}_{\UnbracedArc} + \RechargeVar_{\UnbracedArc}^\Line + \sum_{\Vertex \in \SetOfVertices} \RechargeVar^\Line_\Vertex \right) 
        \hfill
        \label{mips:objective}
    \end{multline}
    \begin{multline}
        \quad \quad \sum_{\Arc \in \IncomingArcs(\Vertex)} \ArcVar^{\Line}_{\UnbracedArc} - \sum_{\Arc \in \OutgoingArcs(\Vertex)} \ArcVar^{\Line}_{\UnbracedArc} = 0
        \hfill \AllLines, \AllVertices \setminus \{\StartDepot, \EndDepot\}
        \label{cons:flow-conservation}
    \end{multline}
    \begin{multline}
        \quad \quad \sum_{\Arc \in \IncomingArcs(\Vertex)} \ArcVar^{\Line}_{\UnbracedArc} \geq 1
        \hfill \AllLines, \forall \Vertex \in \StopsOfLine \setminus \{\EndDepot\}
        \label{cons:mandatory-visits}
    \end{multline}
    \begin{multline}
        \quad \quad \sum_{\Arc \in \OutgoingArcs(\StartDepot)} \ArcVar^{\Line}_{\UnbracedArc} \geq 1
        \hfill \AllLines
        \label{cons:start-depot}
    \end{multline}
    \begin{multline}
        \quad \quad \sum_{\Segment \in \SetOfSegments_{\Station}} \DynamicChargerVar_{\UnbracedSegment}
        \leq |\SetOfSegments_{\Station}| \medspace \TransformerVar_{\Station} \hfill \forall \Station \in \SetOfDynamicStations 
        \label{cons:prevent-dynamic-station-without-transformer}
    \end{multline}
    \begin{multline}
        \quad \quad \ArrivalTimeVar^{\Line}_{\StopOfLine} \leq  \ArrivalTimeVar^{\Line}_{\StopOfLine + 1}
        \hfill \AllLines, \forall \StopOfLine \in \StopsOfLine \setminus \{\StartDepot\}
        \label{cons:stop-order}
    \end{multline}
    \begin{multline}
        \quad \quad \ArrivalTimeVar^{\Line}_{i} + \TravelTime_{\UnbracedArc} \leq \ArrivalTimeVar^{\Line}_{j} - \TimeStep \  \ChargingTimeVar^{\Line}_{j} + (1 - \ArcVar^{\Line}_{\UnbracedArc}) \TwEnd_{\EndDepot}
        \hfill \AllArcs, \AllLines
        \label{cons:discrete-time-propagation}
    \end{multline}
    \begin{multline}
        \quad \quad \StateOfChargeVar^{\Line}_{i} - \Consumption_{\UnbracedArc} + \RechargeVar^{\Line}_{\UnbracedArc} \geq \StateOfChargeVar^{\Line}_{j} - \RechargeVar^{\Line}_{j} - (1 - \ArcVar^{\Line}_{\UnbracedArc})  (\MaxSoC - \MinSoC + \Consumption_{\UnbracedArc})
        \hfill \AllArcs, \AllLines
        \label{cons:charge-propagation-greater}
    \end{multline}
    \begin{multline}
        \quad \quad \StateOfChargeVar^{\Line}_{i} - \Consumption_{\UnbracedArc} + \RechargeVar^{\Line}_{\UnbracedArc} \leq \StateOfChargeVar^{\Line}_{j} - \RechargeVar^{\Line}_{j} + (1 - \ArcVar^{\Line}_{\UnbracedArc})  (\MaxSoC - \MinSoC + \Consumption_{\UnbracedArc})
        \hfill \AllArcs, \AllLines
        \label{cons:charge-propagation-smaller}
    \end{multline}
    \begin{multline}
        \quad \quad \TwBegin_{\Vertex} \leq \ArrivalTimeVar^{\Line}_{\Vertex} \leq \TwEnd_{\Vertex}
        \hfill \AllVertices, \AllLines
        \label{cons:time-limit}
    \end{multline}
    \begin{multline}
        \quad \quad \MinSoC \leq \StateOfChargeVar^{\Line}_{\Vertex} - \RechargeVar^{\Line}_{\Vertex}
        \hfill \AllVertices, \AllLines
        \label{cons:lower-charge-limit}
    \end{multline}
    \begin{multline}
        \quad \quad \StateOfChargeVar^{\Line}_{\Vertex} \leq \MaxSoC
        \hfill \AllVertices, \AllLines
        \label{cons:upper-charge-limit}
    \end{multline}
    \begin{multline}
        \quad \quad \StateOfChargeVar^{\Line}_{\StartDepot} = \InitialSoC
        \hfill \AllLines
        \label{cons:initial-soc}
    \end{multline}
    \begin{multline}
        \quad \quad \RechargeVar^{\Line}_{\Station} = \TimeStep \ \ChargingTimeVar^{\Line}_{\Station} \StationChargingRate_{\Station}
        \hfill \forall\Station \in \SetOfStaticStations, \AllLines
        \label{cons:discrete-static-charging-rate-limit}
    \end{multline}
    \begin{multline}
        \quad \quad \RechargeVar^{\Line}_{\Station} \leq \StaticChargerVar_{\Station} (\MaxSoC - \MinSoC)
        \hfill \forall\Station \in \SetOfStaticStations, \AllLines
        \label{cons:prevent-static-charging-if-not-open}
    \end{multline}
    \begin{multline}
        \quad \quad \RechargeVar^{\Line}_{\UnbracedSegment} \leq \ArcVar_{\UnbracedSegment} \TravelTime_{\UnbracedSegment} \StationChargingRate_{\Station}
        \hfill \forall \Station \in \SetOfDynamicStations, (\SegmentStart, \SegmentEnd) \in \SetOfSegments_{\Station}, \AllLines
        \label{cons:dynamic-charging-rate}
    \end{multline}
    \begin{multline}
        \quad \quad \RechargeVar^{\Line}_{\UnbracedSegment} \leq \DynamicChargerVar_{\UnbracedSegment} (\MaxSoC - \MinSoC)
        \hfill \forall \Segment \in \SetOfSegments, \AllLines
        \label{cons:prevent-dynamic-charging-if-not-open}
    \end{multline}
    \begin{multline}
        \quad \quad \ArcVar^{\Line}_{\UnbracedArc} \in \{0, 1\}
        \hfill \forall\Arc \in \SetOfArcs, \AllLines
        \label{domain:routing}
    \end{multline}
    \begin{multline}
        \quad \quad \StaticChargerVar_\Station \in \{0, 1\}				
        \hfill \forall\Station \in \SetOfStaticStations
        \label{domain:static-charging}
    \end{multline}
    \begin{multline}
        \quad \quad \DynamicChargerVar_{\UnbracedSegment} \in \{0, 1\}
        \hfill \forall \Segment \in \SetOfSegments
        \label{domain:dynamic-charging}
    \end{multline}
    \begin{multline}
        \quad \quad \TransformerVar_{\Station} \in \{0, 1\}
        \hfill \forall \Station \in \SetOfDynamicStations
        \label{domain:transformer-variable}
    \end{multline} 
    \begin{multline}
        \quad \quad \StateOfChargeVar^{\Line}_{\Vertex} \in \mathbb{R}^{+}
        \hfill \forall \Vertex \in \SetOfVertices, \AllLines
        \label{domain:state-of-charge}
    \end{multline}
    \begin{multline}
        \quad \quad \RechargeVar^{\Line}_{\Station} \in \mathbb{R}^{+}
        \hfill \forall\Station \in \SetOfStaticStations, \AllLines
        \label{domain:recharge}
    \end{multline}
    \begin{multline}
        \quad \quad \ArrivalTimeVar^{\Line}_{\Vertex} \in \mathbb{R}^{+}
        \hfill \forall \Vertex \in \SetOfVertices, \AllLines
        \label{domain:arrival-time}
    \end{multline}
    \begin{multline}
        \quad \quad \ChargingTimeVar^{\Line}_{\Station} \in \mathbb{N}
        \hfill \forall \Station \in \SetOfStaticStations, \AllLines
        \label{cons:discrete-domain-charging-time}
    \end{multline}        
\end{subequations}

In the objective~\eqref{mips:objective}, we minimize the cost for charging infrastructure and operating the fleet $\SetOfLines$ which consumes energy during operations and recharges at both stationary charging stations represented by vertices and dynamic charging stations represented by arcs. The energy price is constant and equal for both consumption and recharging. Note that we can use different energy prices for consumption and recharging by separating the respective terms in the objective. We design the feasible space as follows. 

Constraints~\eqref{cons:flow-conservation} are classical flow conservation constraints. Constraints~\eqref{cons:mandatory-visits} ensure that all locations in the vehicle's sequence $\StopsOfLine$ are visited. In Constraints~\eqref{cons:start-depot}, we ensure that all vehicles leave the depot represented by its start representation $\StartDepot$. Constraints~\eqref{cons:prevent-dynamic-station-without-transformer} prevent the construction of dynamic charging segments without the construction of the respective station $\Station \in \SetOfDynamicStations$. We ensure that the stop order on every line is reflected in the departure times by Constraints~\eqref{cons:stop-order}. Accordingly, in Constraints~\eqref{cons:discrete-time-propagation} we propagate the departure times of every line $\Line \in \SetOfLines$ depending on the routing decision represented by the decision variables $\ArcVar^{\Line}_{\UnbracedArc}$. When propagating the departure times, we allow idling on arcs. In Constraints~\eqref{cons:charge-propagation-greater} and~\eqref{cons:charge-propagation-smaller}, we propagate \gls{acr:soc} for all vehicles operating on lines $\Line \in \SetOfLines$. Here, we enforce equality between the consumption and recharge profile and the variables $\StateOfChargeVar$ reflecting the \gls{acr:soc} at discrete points in order to make sure that the upper limits on the battery \gls{acr:soc} are respected. Constraints~\eqref{cons:time-limit},~\eqref{cons:lower-charge-limit},~and \eqref{cons:upper-charge-limit} enforce the departure times at stops in the given time window and the \gls{acr:soc} between the given limits for every vehicle. Constraints~\eqref{cons:initial-soc} set the \gls{acr:soc} at the depot to the exogenously given parameter $\InitialSoC$. Constraints~\eqref{cons:discrete-static-charging-rate-limit} and Constraints~\eqref{cons:dynamic-charging-rate} set the \gls{acr:soc} of every vehicle at any charging station $\Station \in \SetOfStations$ according to the time spent at the stationary station or traversing the arc representing the dynamic segment respectively. In Constraints~\eqref{cons:prevent-static-charging-if-not-open} and \eqref{cons:prevent-dynamic-charging-if-not-open}, the correspondence between the potential time spent charging at a station and the resulting \gls{acr:soc} is conditioned on the construction of the charger or the segment respectively. 

The remaining Constraints~\eqref{domain:routing} - \eqref{cons:discrete-domain-charging-time} define and enforce the domains of the decision variables. Note that a continuous-time formulation based on the provided \gls{acr:mip} is straightforward by assuming real-valued decision variables $\ChargingTimeVar$ and adapting Constraints~\eqref{cons:discrete-time-propagation} and \eqref{cons:discrete-static-charging-rate-limit} accordingly. In order to achieve a fair comparison between the \gls{acr:bab} algorithm applied to Problem~\ref{mips:deterministic} and our algorithmic framework, we allow repeated visits to the same physical charger on a single vehicle's route without breaking propagation constraints~\ref{cons:discrete-time-propagation} and \ref{cons:charge-propagation-greater} by inserting copies for every candidate charging station $\Station \in \SetOfStations$. The copies carry a zero investment cost and can only be constructed if their associated originals are constructed as well.

        \section{Sample-based Local Search Operators} \label{app:sb-ls-operators}
        We elaborate on the local search operators in the following. For the description of the local search operators, we denote by $\VehicleBound$ a lower bound on the routing cost per vehicle. 

\noindent \textbf{\glsentrylong{acr:strip-dyn}: } This operator tightens the dynamic charging station configuration such that the \gls{acr:ils} does not construct unnecessary segments. In an extreme case, this includes removing dynamic charging stations from a given solution completely. Algorithm~\ref{alg:strip-dyn} yields a pseudocode description of the operator. 

\begin{algorithm}[!tb]
    \caption{\Gls{acr:strip-dyn}$(\Solution=(\mathbf{\StaticChargerVar},\mathbf{\TransformerVar}, \mathbf{\DynamicChargerVar}, \mathbf{X}, \mathbf{T}, \mathbf{R}))$}
    \label{alg:strip-dyn}
    \begin{algorithmic}[1]
        \State $\Solution^{\star} \gets \Solution$ \label{l:sds-init-best-solution} \Comment{Initialize best solution at initial solution}
        \State Sample $\tilde{\mathcal{F}} \sim \mathcal{U}\left\{\mathcal{I} \subseteq \SetOfDynamicStations_{\mathbf{\TransformerVar}=1} : |\mathcal{I}| = \min \{ \NeighborhoodSize, |\SetOfDynamicStations_{\mathbf{\TransformerVar}=1}|\} \right\}$ \label{l:sds-sample-random-subset} \Comment{Determine random subset}
        \For{$\Station \in \tilde{\mathcal{F}}$}
        \State $\Solution^{\prime} \gets \texttt{TightenConfiguration}(\Solution, \Station)$ \label{l:sds-tighten-config} \Comment{Find a tighter configuration with Algorithm~\ref{alg:strip}}
        \If {$\texttt{Cost}(\Solution^{\prime}) + \epsilon \leq \texttt{Cost}(\Solution^{\star})$}
            \State $\Solution^{\star} \gets \Solution^{\prime}$ \label{l:sds-accept-new-solution} \Comment{Accept new solution}
        \EndIf
        \EndFor
        \State \Return $\Solution^{\star}$ \label{l:sds-return}
    \end{algorithmic}
\end{algorithm}

The operator starts by initializing the best solution~(l.~\ref{l:sds-init-best-solution}) and sampling a maximum of $\NeighborhoodSize$ dynamic charging stations from the dynamic charging stations that are part of the current configuration. For every dynamic station in the sampled subset, the operator tightens the configuration while retaining feasibility by executing Algorithm~\ref{alg:strip}~(l.~\ref{l:sds-tighten-config}), and accepts the resulting new solution if it improves the current solution by at least a small $\epsilon \ll 1$~(l.~\ref{l:sds-accept-new-solution}).

\begin{algorithm}[!b]
    \caption{$\textsc{TightenConfiguration}(\Solution=(\mathbf{\StaticChargerVar},\mathbf{\TransformerVar}, \mathbf{\DynamicChargerVar}, \mathbf{X}, \mathbf{T}, \mathbf{R}), \Station \in \SetOfDynamicStations \text{ with } |\SetOfSegments_{\Station}|=m)$}
    \label{alg:strip}
    \begin{algorithmic}[1]
        \State $lb = -1, \medspace ub=m, \medspace r=m$ \label{l:init-index-bounds} \Comment{Initialize index bounds}
        \State $\DynamicChargerVar_{\UnbracedSegment} = \mathbb{1}_{\{\SegmentIndex \le r\}}, (\SegmentStart_{\SegmentIndex}, \SegmentEnd_{\SegmentIndex}) \in \SetOfSegments_{\Station}$ \label{l:setting-init-config} 
 \Comment{Feasibility is guaranteed}
        \While{$ub - lb > 1$} \label{l:feasibility-check}
            \If {$\Solution\textsc{ is feasible}$}
                \State $ub = r, \medspace r = ub - \lceil \frac{ub}{2} \rceil$ \label{l:decrease-ub} \Comment{Reduce upper index bound}
            \Else 
                \State $lb = r, \medspace r = lb + \lfloor \frac{ub -lb}{2} \rfloor$ \label{l:increase-lb} \Comment{Increase lower index bound}
            \EndIf
            \State $\DynamicChargerVar_{\UnbracedSegment} = \mathbb{1}_{\{\SegmentIndex \le r\}}, (\SegmentStart_{\SegmentIndex}, \SegmentEnd_{\SegmentIndex}) \in \SetOfSegments_{\Station}$ \label{l:update-config} \Comment{Update configuration}
        \EndWhile 
        \State \Return $\Solution$
    \end{algorithmic}
\end{algorithm}

Algorithm~\ref{alg:strip} relies on the principle of a binary search over the dynamic charging station's index set to determine a tight configuration that retains feasibility. In Line~(l.~\ref{l:init-index-bounds}), we initialize the lower and upper index bounds. The lower index bound refers to the greatest dynamic segment index $lb \in \{1,\dots,m\}$ such that $\DynamicChargerVar_{\UnbracedSegment} = \mathbb{1}_{\{\SegmentIndex \le r\}}, (\SegmentStart_{\SegmentIndex}, \SegmentEnd_{\SegmentIndex}) \in \SetOfSegments_{\Station}$ and the algorithm has proven infeasibility for the resulting configuration (given the remainder of the configuration ($\mathbf{\StaticChargerVar}, \mathbf{\TransformerVar}, \mathbf{\DynamicChargerVar}$) is constant). Vice versa, the upper index bound refers to the lowest such index for which the resulting configuration is feasible. Then, the algorithm initializes the configuration determined by the initial upper index bound~(l.~\ref{l:setting-init-config}). The algorithm evaluates feasibility for the resulting configuration and adapts the lower~(l.~\ref{l:increase-lb}) and upper index bounds~(l.~\ref{l:decrease-ub}) accordingly in a binary iterative fashion. After every update of the index bounds, the algorithm updates the current configuration~(l.~\ref{l:update-config}). Finally, the algorithm terminates when neither the lower nor the upper index bound can be tightened further. Note that the solution returned by the operator \gls{acr:strip-dyn} may depend on the order in which the sampled dynamic stations are tightened. Moreover, there is a minor heuristic element involved as the operators enforces solutions in which segments that are part of the configuration are always associated with an index that is smaller than all segments that are not part of the configuration. 

\noindent \textbf{\glsentrylong{acr:remove-stat}: } This operator removes stationary charging stations from a solution's configuration. Algorithm~\ref{alg:remove-stat} provides a description of the operator in pseudocode. The operator randomly samples $\min \{ \NeighborhoodSize, |\SetOfStaticStations_{\mathbf{y}=1}|\}$ stationary charging stations from the set of stationary charging stations in the current configuration~(l.~\ref{l:rss-sample-random-subset}). Then, the operator removes every sampled station from the current configuration individually, evaluates the resulting solution, and accepts moves that yield an epsilon improvement with $\epsilon \ll 1$~(l.~\ref{l:rss-remove-static-station}-\ref{l:rss-accept-new-solution}). Finally, the operator returns the best evaluated solution~(l.~\ref{l:rss-return}).

\begin{algorithm}[!tb]
    \caption{\Gls{acr:remove-stat}$(\Solution=(\mathbf{\StaticChargerVar},\mathbf{\TransformerVar}, \mathbf{\DynamicChargerVar}, \mathbf{X}, \mathbf{T}, \mathbf{R}))$}
    \label{alg:remove-stat}
    \begin{algorithmic}[1]
        \State $\Solution^{\star} \gets \Solution$ \label{l:rss-init-best-solution} \Comment{Initialize best solution at initial solution}
        \State Sample $\tilde{\mathcal{F}} \sim \mathcal{U}\left\{\mathcal{I} \subseteq \SetOfStaticStations_{\mathbf{y}=1} : |\mathcal{I}| = \min \{ \NeighborhoodSize, |\SetOfStaticStations_{\mathbf{y}=1}|\} \right\}$ \label{l:rss-sample-random-subset} \Comment{Determine random subset}
        \For{$\Station \in \tilde{\mathcal{F}}$}
        \State $\Solution^{\prime} \gets \texttt{RemoveAndEvaluate}(\Solution, \Station)$ \label{l:rss-remove-static-station} \Comment{Remove $\Station$ from $\Solution$ and evaluate to new solution}
        \If {$\texttt{Cost}(\Solution^{\prime}) + \epsilon \leq \texttt{Cost}(\Solution^{\star})$}
            \State $\Solution^{\star} \gets \Solution^{\prime}$ \label{l:rss-accept-new-solution} \Comment{Accept new solution}
        \EndIf
        \EndFor
        \State \Return $\Solution^{\star}$ \label{l:rss-return}
    \end{algorithmic}
\end{algorithm}

\noindent \textbf{\glsentrylong{acr:swap-to-dyn}: } This operator exchanges stationary charging stations with dynamic charging stations. Algorithm~\ref{alg:swap-to-dyn} yields a description of the operator.

\begin{algorithm}[!b]
    \caption{\Gls{acr:swap-to-dyn}$(\Solution=(\mathbf{\StaticChargerVar},\mathbf{\TransformerVar}, \mathbf{\DynamicChargerVar}, \mathbf{X}, \mathbf{T}, \mathbf{R}))$}
    \label{alg:swap-to-dyn}
    \begin{algorithmic}[1]
        \State $\Solution^{\star} \gets \Solution$ \label{l:s2d-init-best-solution} \Comment{Initialize best solution at initial solution}
        \State Sample $\tilde{\mathcal{F}} \sim \mathcal{U}\left\{\mathcal{I} \subseteq \SetOfStaticStations_{\mathbf{y}=1} : |\mathcal{I}| = \min \{ \NeighborhoodSize, |\SetOfStaticStations_{\mathbf{y}=1}|\} \right\}$ \label{l:s2d-sample-random-subset} \Comment{Determine random subset}
        \For{$\Station \in \tilde{\mathcal{F}}$}
        \State $\Station^{\prime} \gets \texttt{GetClosestStation}(\SetOfDynamicStations_{\mathbf{\TransformerVar}=0}, \Station)$ \label{l:s2d-get-closest-station} \Comment{By geometric distance to first segment start}
        \State $\Solution^{\prime} \gets \texttt{Remove}(\Solution, \Station)$ \label{l:s2d-remove} 
        \Comment{Remove stationary station from configuration}
        \If{$\eqref{eq:infrastructure-cost} \texttt{ evaluated for }\Solution^{\prime} + \sum_{\Vehicle \in \SetOfVehicles} \VehicleBound + \min_{\Segment \in \SetOfSegments_{\Station^{\prime}}} \{ \Cost_{\UnbracedSegment} + \Cost_{\Station^{\prime}} \} \geq \texttt{Cost}(\Solution^{\star})$}
            \State $\texttt{continue}$ \label{s2d:evaluate-min-total-cost} \Comment{Check if move can yield improvement}
        \EndIf
        \State $\Solution^{\prime} \gets \texttt{Add}(\Solution^{\prime}, \Station^{\prime})$ \label{l:s2d-add} 
        \State $\Solution^{\prime} \gets \texttt{TightenConfiguration}(\Solution^{\prime}, \Station^{\prime})$ \label{l:s2d-tighten-config} \Comment{Find a tighter configuration with Algorithm~\ref{alg:strip}}
        \If {$\texttt{Cost}(\Solution^{\prime}) + \epsilon \leq \texttt{Cost}(\Solution^{\star})$}
            \State $\Solution^{\star} \gets \Solution^{\prime}$ \label{l:s2d-accept-new-solution} \Comment{Accept new solution}
        \EndIf
        \EndFor
        \State \Return $\Solution^{\star}$ \label{l:s2d-return}
    \end{algorithmic}
\end{algorithm}

First, we initialize the best solution at $\Solution$~(l.~\ref{l:s2d-init-best-solution}). Second, we sample $\min \{ \NeighborhoodSize, |\SetOfStaticStations_{\mathbf{y}=1}|\}$ stationary charging stations. For every station in the sample, we determine the closest dynamic station that is not in the current configuration~(l.~\ref{l:s2d-get-closest-station}) and afterwards remove the sampled station from the current configuration~(l.~\ref{l:s2d-remove}). In this context, we use the geometric distance between the stationary charging station and the first segment start of the dynamic charging stations as distance measure. The operator checks if replacing the sampled station $\Station$ by its closest counterpart $\Station^{\prime}$ can improve the objective. The check adds up the minimum associated infrastructure cost when only the shortest segment of $\Station^{\prime}$ belongs to the solution and the minium routing cost of the fleet (cf. Section~\ref{sec:subproblem}). If the operator confirms this check~(l.~\ref{s2d:evaluate-min-total-cost}), it continues by adding $\Station^{\prime}$ to the solution~(l.~\ref{l:s2d-add}) and then tightening the dynamic charging stations by calling Algorithm~\ref{alg:strip}~(l.~\ref{l:s2d-tighten-config}). We accept new solutions if they provide a decrease in objective value by at least a small $\epsilon \ll 1$~(l.~\ref{l:s2d-accept-new-solution}). 

\noindent \textbf{\glsentrylong{acr:swap-to-stat}: } This operator follows a similar approach as \gls{acr:swap-to-dyn}. Algorithm~\ref{alg:swap-to-stat} provides the respective pseudocode description.

\begin{algorithm}[!tb]
    \caption{\Gls{acr:swap-to-stat}$(\Solution=(\mathbf{\StaticChargerVar},\mathbf{\TransformerVar}, \mathbf{\DynamicChargerVar}, \mathbf{X}, \mathbf{T}, \mathbf{R}))$}
    \label{alg:swap-to-stat}
    \begin{algorithmic}[1]
        \State $\Solution^{\star} \gets \Solution$ \label{l:s2s-init-best-solution} \Comment{Initialize best solution at initial solution}
        \State Sample $\tilde{\mathcal{F}} \sim \mathcal{U}\left\{\mathcal{I} \subseteq \SetOfStaticStations_{\mathbf{y}=1} \cup \SetOfDynamicStations_{\mathbf{\TransformerVar}=1}: |\mathcal{I}| = \min \{ \NeighborhoodSize, |\SetOfStaticStations_{\mathbf{y}=1}| + |\SetOfDynamicStations_{\mathbf{\TransformerVar}=1}|\} \right\}$ \label{l:s2s-sample-random-subset} \Comment{Determine random subset}
        \For{$\Station \in \tilde{\mathcal{F}}$}
        \State $\Station^{\prime} \gets \texttt{GetClosestStation}(\SetOfStaticStations_{\mathbf{y}=0}, \Station)$ \label{l:s2s-get-closest-station} \Comment{By geometric distance}
        \State $\Solution^{\prime} \gets \texttt{Remove}(\Solution, \Station)$ \label{l:s2s-remove} 
        \Comment{Remove station from configuration}
        \If{$\eqref{eq:infrastructure-cost} \texttt{ evaluated for }\Solution^{\prime} + \sum_{\Vehicle \in \SetOfVehicles} \VehicleBound + \Cost_{\Station^{\prime}} \geq \texttt{Cost}(\Solution^{\star})$}
            \State $\texttt{continue}$ \label{l:s2s-evaluate-min-total-cost} \Comment{Check if move can yield improvement}
        \EndIf
        \State $\Solution^{\prime} \gets \texttt{Add}(\Solution^{\prime}, \Station^{\prime})$ \label{l:s2s-add} 
        \If {$\texttt{Cost}(\Solution^{\prime}) + \epsilon \leq \texttt{Cost}(\Solution^{\star})$}
            \State $\Solution^{\star} \gets \Solution^{\prime}$ \label{l:s2s-accept-new-solution} \Comment{Accept new solution}
        \EndIf
        \EndFor
        \State \Return $\Solution^{\star}$ \label{l:s2s-return}
    \end{algorithmic}
\end{algorithm}

This operator samples $\min \{ \NeighborhoodSize, |\SetOfStaticStations_{\mathbf{y}=1}| + |\SetOfDynamicStations_{\mathbf{\TransformerVar}=1}|\}$ from the stations in the current configuration~(l.~\ref{l:s2s-sample-random-subset}). In every move, the operator removes one station from the sample from the configuration~(l.~\ref{l:s2s-remove}) and determines their closest stationary counterpart that is not part of the current configuration~(l.~\ref{l:s2s-get-closest-station}). Subsequently, the operator evaluates if the resulting configuration could lead to an improved solution~(l.~\ref{l:s2s-evaluate-min-total-cost}). If this is not the case, the operator continues evaluating the next swap. Otherwise, the operator adds station $\Station^{\prime}$ to the configuration~(l.~\ref{l:s2s-add}) and evaluates the configuration. If the evaluation yields an improvement by at least a small $\epsilon \ll 1$ in solution value, the operator accepts the corresponding solution~(l.~\ref{l:s2s-accept-new-solution}).

\noindent \textbf{\glsentrylong{acr:add-dyn}: } Algorithm~\ref{alg:add-dyn} yields a description of this operator.

\begin{algorithm}[!tb]
    \caption{\Gls{acr:add-dyn}$(\Solution=(\mathbf{\StaticChargerVar},\mathbf{\TransformerVar}, \mathbf{\DynamicChargerVar}, \mathbf{X}, \mathbf{T}, \mathbf{R}))$}
    \label{alg:add-dyn}
    \begin{algorithmic}[1]
        \State $\Solution^{\star} \gets \Solution$ \label{l:ad-init-best-solution} \Comment{Initialize best solution at initial solution}
        \State \State Sample $\tilde{\mathcal{F}} \sim \mathcal{U}\left\{\mathcal{I} \subseteq \SetOfStaticStations_{\mathbf{y}=1} \cup \SetOfDynamicStations_{\mathbf{\TransformerVar}=1}: |\mathcal{I}| = \min \{ \NeighborhoodSize, |\SetOfStaticStations_{\mathbf{y}=1}| + |\SetOfDynamicStations_{\mathbf{\TransformerVar}=1}|\} \right\}$ \label{l:ad-sample-random-subset} \Comment{Determine random subset}
        \For{$\Station \in \tilde{\mathcal{F}}$}
        \State $\Station^{\prime} \gets \texttt{GetFurthestStation}(\SetOfDynamicStations_{\mathbf{\TransformerVar}=0}, \Station)$ \label{l:ad-get-furthest-station} \Comment{By geometric distance}
        \If{$\eqref{eq:infrastructure-cost} \texttt{ evaluated for }\Solution + \sum_{\Vehicle \in \SetOfVehicles} \VehicleBound + \min_{\Segment \in \SetOfSegments_{\Station^{\prime}}} \{ \Cost_{\UnbracedSegment} + \Cost_{\Station^{\prime}} \} \geq \texttt{Cost}(\Solution^{\star})$}
            \State $\texttt{continue}$ \label{l:ad-evaluate-min-total-cost} \Comment{Check if move can yield improvement}
        \EndIf
        \State $\Solution^{\prime} \gets \texttt{Add}(\Solution, \Station^{\prime})$ \label{l:ad-add} 
        \State $\Solution^{\prime} \gets \texttt{TightenConfiguration}(\Solution^{\prime}, \Station^{\prime})$ \label{l:ad-tighten-config} \Comment{Find a tighter configuration with Algorithm~\ref{alg:strip}}
        \If {$\texttt{Cost}(\Solution^{\prime}) + \epsilon \leq \texttt{Cost}(\Solution^{\star})$}
            \State $\Solution^{\star} \gets \Solution^{\prime}$ \label{l:ad-accept-new-solution} \Comment{Accept new solution}
        \EndIf
        \EndFor
        \State \Return $\Solution^{\star}$ \label{l:ad-return}
    \end{algorithmic}
\end{algorithm}

First, the operator initializes $\Solution^{\star}$ and samples a maximum of $\min \{ \NeighborhoodSize, |\SetOfDynamicStations_{\mathbf{\TransformerVar}=0}|\}$ from the stations in the current configuration~(l.~\ref{l:ad-sample-random-subset}). Then, the operator loops over the sample, determines the dynamic charging station that is not part of the configuration and bears the greatest distance to the current station from the sample~(l.~\ref{l:ad-get-furthest-station}) and adds this dynamic charging station to the configuration~(l.~\ref{l:ad-add}). The evaluation of the move includes applying Algorithm~\ref{alg:strip} such that the added dynamic charging station is not excessively long~(l.~\ref{l:ad-tighten-config}). It waives the explicit evaluation of the move if the resulting infrastructure cost summed with the lower bound on the routing cost per vehicle exceeds the current best solution value~(l.~\ref{l:ad-evaluate-min-total-cost}). Finally, the operator accepts new solutions if they improve the solution value of $\Solution^{\star}$ by more than a small $\epsilon \ll 1$. 

\noindent \textbf{\glsentrylong{acr:add-stat}: } This operator follows the equivalent approach to \gls{acr:add-dyn} with the only difference that it evaluates adding a stationary charging station to the given configuration. Algorithm~\ref{alg:add-stat} describes the operator analogously to the description of the other operators.

\begin{algorithm}[!b]
    \caption{\Gls{acr:add-stat}$(\Solution=(\mathbf{\StaticChargerVar},\mathbf{\TransformerVar}, \mathbf{\DynamicChargerVar}, \mathbf{X}, \mathbf{T}, \mathbf{R}))$}
    \label{alg:add-stat}
    \begin{algorithmic}[1]
        \State $\Solution^{\star} \gets \Solution$ \label{l:as-init-best-solution} \Comment{Initialize best solution at initial solution}
        \State \State Sample $\tilde{\mathcal{F}} \sim \mathcal{U}\left\{\mathcal{I} \subseteq \SetOfStaticStations_{\mathbf{y}=1} \cup \SetOfDynamicStations_{\mathbf{\TransformerVar}=1}: |\mathcal{I}| = \min \{ \NeighborhoodSize, |\SetOfStaticStations_{\mathbf{y}=1}| + |\SetOfDynamicStations_{\mathbf{\TransformerVar}=1}|\} \right\}$ \label{l:as-sample-random-subset} \Comment{Determine random subset}
        \For{$\Station \in \tilde{\mathcal{F}}$}
        \State $\Station^{\prime} \gets \texttt{GetFurthestStation}(\SetOfStaticStations_{\mathbf{y}=0}, \Station)$ \label{l:as-get-furthest-station} \Comment{By geometric distance}
        \If{$\eqref{eq:infrastructure-cost} \texttt{ evaluated for }\Solution + \sum_{\Vehicle \in \SetOfVehicles} \VehicleBound + \Cost_{\Station^{\prime}} \geq \texttt{Cost}(\Solution^{\star})$}
            \State $\texttt{continue}$ \label{l:as-evaluate-min-total-cost} \Comment{Check if move can yield improvement}
        \EndIf
        \State $\Solution^{\prime} \gets \texttt{Add}(\Solution, \Station^{\prime})$ \label{l:as-add} 
        \If {$\texttt{Cost}(\Solution^{\prime}) + \epsilon \leq \texttt{Cost}(\Solution^{\star})$}
            \State $\Solution^{\star} \gets \Solution^{\prime}$ \label{l:as-accept-new-solution} \Comment{Accept new solution}
        \EndIf
        \EndFor
        \State \Return $\Solution^{\star}$ \label{l:as-return}
    \end{algorithmic}
\end{algorithm}

First, the operator initializes $\Solution^{\star}$. Then, it samples a subset of the stations in the current configuration~(l.~\ref{l:as-sample-random-subset}). For every station in the sample, the operator identifies the stationary charging station that is not in the current configuration and bears the greatest distance to the sampled one~(l.~\ref{l:as-sample-random-subset}). The operator evaluates the configuration after adding this station if the resulting infrastructure cost added to the lower bound on the routing cost for all vehicles remains below the current best solution value~(l.~\ref{l:as-evaluate-min-total-cost}). In line with the other operators, Algorithm~\ref{alg:add-stat} accepts new solutions if they improve the solution value of $\Solution^{\star}$ by more than a small $\epsilon \ll 1$~(l.~\ref{l:as-accept-new-solution}). 
\end{appendices}

\end{document}